\documentclass[10pt,leqno]{amsart}
\ifx\relax \undefined
  \usepackage[T1]{fontenc}
  \usepackage[utf8]{inputenc}
\fi
\usepackage{lmodern,amsmath,amssymb,amsthm,mathtools}
\usepackage[a4paper,margin=25mm]{geometry}
\usepackage{enumerate,xcolor}
\usepackage[pagebackref,colorlinks=true,linkcolor=blue,citecolor=blue,urlcolor=blue]{hyperref}
\ifdefined\backrefalt
  \renewcommand*{\backref}[1]{}%
  \renewcommand*{\backrefalt}[4]{\ifcase #1\or $\uparrow$#2\else $\uparrow$#2\fi}%
\fi
\hypersetup{pdftitle={Schubert positivity for filtered bundles}}

\numberwithin{equation}{section}
\newtheorem{theorem}{Theorem}[section]
\newtheorem{proposition}[theorem]{Proposition}
\newtheorem{lemma}[theorem]{Lemma}
\newtheorem{corollary}[theorem]{Corollary}
\theoremstyle{definition}
\newtheorem{definition}[theorem]{Definition}

\newtheorem{remark}[theorem]{Remark}

\newcommand{\C}{\mathbb C}
\newcommand{\R}{\mathbb R}
\newcommand{\PP}{\mathbb P}
\newcommand{\OO}{\mathcal O}
\newcommand{\ii}{\mathrm i}
\DeclareMathOperator{\Gr}{Gr}
\DeclareMathOperator{\Pol}{Pol}
\DeclareMathOperator{\IH}{IH}
\DeclareMathOperator{\IC}{IC}
\DeclareMathOperator{\rk}{rk}
\DeclareMathOperator{\Ann}{Ann}
\newcommand{\clpol}[1]{\overline{\Pol(#1)}}

\allowdisplaybreaks[2]
\usepackage{microtype}
\title[Schubert positivity for filtered bundles]{Schubert positivity for filtered bundles}
\author{Enhan Li}
\address{Department of Mathematics, University of Science and Technology of China}
\email{lablehustc0719@mail.ustc.edu.cn}
\date{}

\begin{document}
\begin{abstract}
We characterize formal filtered nefness of a formal quotient flag on a polarized Hodge--Lefschetz module by the formal nefness of its positive-rank total quotients. Using this criterion, we prove that quotient Schubert functors preserve formal nefness and that Schubert incidence modules inherit every commuting formally filtered nef flag. As geometric and numerical consequences, we recover the positivity criteria of Fulton--Lazarsfeld and Fulton and extend the Lorentzian polynomial theorem of Lu and Zheng to derived Schubert classes of filtered bundles. In combinatorics, our framework applies to tautological flags of flag matroids on ordinary and augmented matroid Chow rings as well as on the intersection cohomology of rational projective fans equipped with a common weight map. We deduce concordance-sensitive intersection inequalities, establish the Schubert positivity of the leading Riemann--Roch Chow class of a flag matroid, and show that the nonvanishing supports of derived Schubert numbers on an augmented Chow ring form the Higgs lifts of an intrinsically determined matroid quotient.
\end{abstract}
\maketitle
\tableofcontents

\section{Introduction}
Let $F_s\twoheadrightarrow\cdots\twoheadrightarrow F_1\twoheadrightarrow F_0=0$ be a flag of quotient vector bundles on a projective variety. When the largest bundle $F_s$ is nef, every quotient $F_i$ remains nef, even though a successive kernel $G_i=\ker(F_i\to F_{i-1})$ may fail to be nef. For example, the kernel of a general surjection $\OO^{\oplus2}\twoheadrightarrow\OO(1)$ on $\PP^1$ is $\OO(-1)$. Nevertheless, any Schubert polynomial evaluated at the Chern roots of $G_1,\ldots,G_s$ in the order of the flag has nonnegative degree in top dimension. For ample filtered bundles, Fulton proved that a nonzero homogeneous polynomial has positive degree on all flags of prescribed ranks if and only if its expansion in the Schubert basis has nonnegative coefficients \cite[\S1, Theorem]{Ful95}, extending the single-bundle Schur positivity theorem of Fulton and Lazarsfeld \cite[Theorem~I]{FL83}.

An analogous positivity problem arises in combinatorics for the tautological classes of a flag matroid $M_1\preceq\cdots\preceq M_s$, where $N\preceq M$ means that $N$ is a matroid quotient of $M$. By recent work of Larson and Partida \cite[Theorem~1.4]{LP}, the projective bundle ring associated with each individual constituent $M_i$ satisfies the hard Lefschetz property and the Hodge--Riemann relations. However, a Schubert polynomial with multiple descents involves successive differences of the tautological classes of consecutive constituents, which need not come from nef vector bundles even when the matroids are realizable. Understanding why such multi-descent Schubert polynomials remain positive---and how their positivity depends on the quotient relations $M_{i-1}\preceq M_i$---requires a filtered Hodge--Lefschetz formalism that operates directly on Chern operators without assuming a common linear realization. We develop this formalism for formal quotient flags over general polarized Hodge--Lefschetz modules and apply it both to algebraic geometry and to matroid theory.

\subsection*{Formal quotient flags}
Let $(V,S)$ be a Hodge--Lefschetz module of weight $d$ equipped with a polarizing cone $\mathcal C$. A formal bundle $F$ of rank $e$ over $V$ is specified by commuting Chern operators $c_j(F)$ of bidegree $(j,j)$ for $1\le j\le e$. When $e>0$, the associated projective bundle module $W(F)$ carries a relative hyperplane operator $\xi$ satisfying Grothendieck's relation, and $F$ is called \emph{formally nef} if $h+t\xi$ polarizes $W(F)$ for every $h\in\mathcal C$ and $t>0$ \cite[Definition~8.1]{LZ}. A formal quotient flag $F_\bullet$ with marked ranks $0=r_0<r_1<\cdots<r_s$ consists of formal graded pieces $G_i$ of ranks $r_i-r_{i-1}$, whose total quotients $F_j$ are determined by the Whitney relation $c(F_j)=\prod_{i\le j}c(G_i)$.

To formulate filtered nefness for $F_\bullet$, let $B$ be a smooth complex projective variety carrying a flag of quotient vector bundles $W_1\twoheadrightarrow\cdots\twoheadrightarrow W_s$ (ordered from largest rank to smallest, allowing repeated stages and trailing zero bundles) with $W_1$ globally generated, and form the formal incidence bundle
\[
K=\bigoplus_{i=1}^s G_i\boxtimes W_i
\]
over $V\otimes H^*(B,\R)$. We say that $F_\bullet$ is \emph{formally filtered nef} if every such incidence bundle $K$ of positive rank is formally nef with respect to the cone $\mathcal C+\operatorname{Amp}(B)$. Taking $B$ to be a point with $W_i=\C$ for $i\le j$ and $W_i=0$ for $i>j$ recovers the formal nefness of the $j$th total quotient $F_j$. Our first main result shows that these $s$ conditions over a point already imply the incidence condition over arbitrary smooth projective bases.

\begin{theorem}[Total-quotient criterion]\label{intro:criterion}
A formal quotient flag on $(V,\mathcal C)$ is formally filtered nef if and only if each of its total quotients of positive rank is formally nef.
\end{theorem}

We establish sufficiency first in the geometric case where $G_1,\ldots,G_s$ are vector bundles on a smooth projective variety $Y$, testing against quotient flags $W_\bullet$ on $Y$ via $\bigoplus_i G_i\otimes W_i$ (Theorem~\ref{thm:filtered}). The proof compares two resolutions of a determinantal locus: on each perverse support of the first resolution, the second resolution expresses the pullback of the hyperplane class in terms of a shorter flag. For a general Hodge--Lefschetz module $(V,\mathcal C)$, we carry out this comparison universally over a product of Grassmannians and specialize the universal Chern classes to operators on $V$ (Theorem~\ref{thm:abstract-stages}).

\subsection*{Schubert functors and Schubert incidence modules}
Fix a formally filtered nef flag $F_\bullet$ on $(V,\mathcal C)$ of marked ranks $r_1<\cdots<r_s$. On the partial flag variety $T=\operatorname{Fl}(r_1,\ldots,r_s;n)$ with universal quotient bundles $Q_i$ of rank $r_i$, choosing $n$ sufficiently large yields a specialization homomorphism $\varphi$ from $H^{2*}(T,\R)$ to the algebra of operators on $V$ sending $c_j(Q_i)$ to $c_j(F_i)$. For any vector bundle $U$ on $T$, we write $\varphi(U)$ for the formal bundle on $V$ of rank $\rk U$ with Chern operators $\varphi(c_j(U))$.

\begin{theorem}[Specialization and quotient Schubert functors]\label{intro:functors}
Every nef vector bundle $U$ of positive rank on $T$ specializes to a formally nef bundle $\varphi(U)$ on $V$. More generally, if $U_\ell\twoheadrightarrow\cdots\twoheadrightarrow U_1$ is a flag of quotient bundles on $T$ with $U_\ell$ nef of positive rank, then the specializations $\varphi(U_j)$ form the total quotients of a formally filtered nef flag after omitting zero bundles and repeated ranks.

In particular, the Schur bundles $\mathbf S_\lambda(Q_i)$ and the quotient Schubert bundles $\mathbf S'_w(Q_\bullet)$ of Kra\'skiewicz and Pragacz specialize to formally nef bundles whenever they have positive rank, provided in the partial-flag case that the Schubert bundle on a complete refinement descends to $T$ (which holds whenever the character of $\mathbf S'_w$ is symmetric within each block between consecutive marked ranks).
\end{theorem}

The specialization statement follows by realizing the diagonal class on $T\times T$ as the top Chern class of an incidence bundle and applying top Chern descent (Theorem~\ref{thm:flag-specialization}). Since universal quotient Schubert bundles are globally generated, the assertion for Schubert functors follows immediately (Theorem~\ref{thm:schubert-preservation} and Corollary~\ref{cor:schur-formal}), extending Fulton's observation that Schubert functors preserve ampleness for filtered vector bundles \cite[\S4]{Ful95}.

Beyond bundle functors, Schubert polynomials act directly as commuting endomorphisms on $V$. Let $\mathcal F_1,\ldots,\mathcal F_r$ be formally filtered nef flags on $V$ whose Chern operators commute with one another and with the polarizing cone $\mathcal C$. For each $1\le i\le r$, choose a permutation $w_i$ whose descents lie among the marked ranks of $\mathcal F_i$, and set
\[
g_i=\mathfrak S_{w_i},\qquad k_i=\ell(w_i),\qquad \Gamma=\prod_{i=1}^r g_i(\mathcal F_i),\qquad k=\sum_{i=1}^r k_i.
\]
Although the quotient module $\Gamma V$ itself need not satisfy higher-degree Hodge--Riemann relations---as demonstrated by an explicit counterexample on $(\PP^2)^3$ in Subsection~\ref{subsec:pairing-base}---we construct an auxiliary polarized Hodge--Lefschetz module that lifts the pairing against $\Gamma$.

\begin{proposition}[Schubert incidence module]\label{intro:incidence}
There exist a smooth projective variety $Z$, a Hodge--Lefschetz module $(W,S_W)$ of weight $d-k$ polarized by $\mathcal C+\operatorname{Amp}(Z)$, and a linear map $j\colon V\to W$ commuting with operators in $\operatorname{span}(\mathcal C)$ such that
\[
S_W(jx,jy)=S(x,\Gamma y)\qquad(x,y\in V).
\]
Moreover, every formally filtered nef flag on $V$ whose Chern operators commute with those of $\mathcal F_1,\ldots,\mathcal F_r$ induces a formally filtered nef flag on $W$.
\end{proposition}

Here $Z$ is a product of Bott--Samelson resolutions of Schubert varieties, and $W$ is the image of a product of top Chern operators on $V\otimes H^*(Z,\R)$ (Proposition~\ref{prop:schubert-incidence-module} and Theorem~\ref{thm:schubert-incidence-preservation}). Restricting the mixed Hodge--Riemann form on $W$ along $j$ yields a signature bound on $V$, which in bidegree $(1,1)$ specializes to a Hodge-index inequality for $\Gamma$.

\subsection*{Numerical inequalities and Lorentzian polynomials}
Assume now that $V$ is effective, meaning that $V^{p,q}=0$ whenever $\min(p,q)<0$, and that $V^{0,0}=\C e$ for a chosen real generator $e$. For an operator $A$ of bidegree $(d,d)$, we write $\int_V A=S(e,Ae)$, and we call a $(1,1)$-operator \emph{nef} if it lies in $\overline{\mathcal C}$.

\begin{corollary}[Schubert volume polynomial]\label{intro:volume}
If $m=d-k\ge0$ and $D_1,\ldots,D_s\in\overline{\mathcal C}$ are commuting nef divisors, then the homogeneous polynomial
\[
\frac1{m!}\int_V\Gamma(t_1D_1+\cdots+t_sD_s)^m
\]
is Lorentzian or zero. In particular, $\int_V\Gamma D_1\cdots D_m\ge0$ for all nef divisors $D_1,\ldots,D_m$.
\end{corollary}

Indeed, Proposition~\ref{intro:incidence} identifies these intersection numbers with mixed divisor pairings against $j(e)$ on the polarized module $W$ (Theorem~\ref{thm:ineq-schubert-volume}). When $m=0$, Corollary~\ref{intro:volume} recovers the nonnegativity of top-degree Schubert numbers for nef flags. Testing against dual Schubert cycles yields the characterization of numerically positive polynomials due to Fulton--Lazarsfeld and Fulton (Corollary~\ref{cor:classical-pos}). When $m\ge2$, taking $L$ to be a product of $m-2$ nef divisors yields the Khovanskii--Teissier type inequality
\[
\left(\int_V\Gamma LHK\right)^2 \ge \left(\int_V\Gamma LH^2\right)\left(\int_V\Gamma LK^2\right)
\]
for all nef divisors $H$ and $K$.

While Corollary~\ref{intro:volume} holds the Schubert operator $\Gamma$ fixed, we also compare Schubert classes across different codimensions by translating all Chern roots of a flag by an indeterminate $u$ and expanding $g(x+u\mathbf 1)=\sum_{j\ge0}g^{[j]}(x)u^j$. Allowing an independent translation parameter for each flag leads to the following Lorentzianity theorem for derived Schubert classes.

\begin{theorem}[Derived Schubert classes]\label{intro:derived}
Let $q,n_1,\ldots,n_r$ be nonnegative integers satisfying $d=q+\sum_{i=1}^r(k_i-n_i)$. Then the polynomial
\[
\sum_{\substack{|a|=q\\0\le a_i\le n_i}} \frac{x_1^{a_1}\cdots x_r^{a_r}}{a_1!\cdots a_r!}\int_V\prod_{i=1}^r g_i^{[n_i-a_i]}(\mathcal F_i)
\]
is Lorentzian or zero. If $r\ge1$, every positive-rank total quotient of each $\mathcal F_i$ is formally ample, $q\le n_i\le k_i$, and $g_i(1,\ldots,1)>0$ for all $1\le i\le r$, then this polynomial is strictly Lorentzian.
\end{theorem}

For one-stage flags of nef vector bundles, Theorem~\ref{intro:derived} recovers the projective case of Lu and Zheng's theorem \cite[Theorem~1.9]{LuZhengHR}. For two nef bundles, the resulting log-concavity of derived Schur numbers is due to Ross and Toma \cite[Theorem~10.5]{RTnef}. We prove the nef statement by applying Corollary~\ref{intro:volume} over $V\otimes H^*(\prod_i\PP^{n_i},\R)$ after twisting the $i$th flag by the hyperplane class of $\PP^{n_i}$ (Theorem~\ref{thm:ineq-luzheng}), and we establish strict Lorentzianity by using formal ampleness to perturb all total quotients negatively before differentiating the Hessian along translation parameters (Theorem~\ref{thm:strict-filtered} and Corollary~\ref{cor:stagewise-strict}).

\subsection*{Flag matroids and combinatorial applications}
A central motivation for working over abstract Hodge--Lefschetz modules comes from matroid theory. When a flag matroid $M_1\preceq\cdots\preceq M_s$ on a finite ground set $E$ is realizable over $\C$ by a chain of linear quotients $\C^E\twoheadrightarrow V_s\twoheadrightarrow\cdots\twoheadrightarrow V_1$, it induces a genuine flag of globally generated quotient vector bundles on the permutohedral variety $X_E$, so classical Schubert positivity applies directly. For nonrealizable flag matroids---or even for pairs of realizable matroids admitting no simultaneous realization---no such flag of nef vector bundles exists: the $K$-theoretic tautological classes $F_i=\mathcal S_{M_i}^\vee\in K^0(X_E)$ of Berget, Eur, Spink, and Tseng \cite[Definition~3.9]{BEST} have virtual differences $F_i-F_{i-1}$ whose line-bundle representatives typically fail ordinary nefness. At the same time, positivity results for single matroids or arbitrary tuples of matroids \cite{LP,LS} apply to products of Schur classes of the individual constituents $F_i$, which cannot detect whether $M_{i-1}$ is a quotient of $M_i$.

Our total-quotient criterion bridges this gap. By Larson and Partida \cite[Theorem~1.4]{LP}, each positive-rank constituent $F_i$ is formally nef on $X_E$, while the quotient relation $M_{i-1}\preceq M_i$ ensures that the ratios $c(F_i)/c(F_{i-1})$ have the prescribed successive ranks. Theorem~\ref{intro:criterion} therefore upgrades the chain $(F_1,\ldots,F_s)$ into a formally filtered nef flag, and a geometric resolution-and-descent argument extends this property to pullbacks along arbitrary morphisms $Y\to X_E$ from smooth projective varieties (Theorem~\ref{thm:mat}).

More generally, let $\Sigma$ be a rational fan and let $f$ be an integral linear map from its cocharacter lattice to $\mathbb Z^E$ carrying each cone of $\Sigma$ into a cone of the normal fan of every base polytope $P(M_i)$. On any cone where $B_i$ is a $f$-minimum-weight basis of $M_i$, put
\[
C_i(u)=\prod_{e\in B_i}(1-f^*t_e\,u),
\]
where $t_e$ denotes the $e$th coordinate function on $\R^E$. Because the same weight map $f$ governs all constituents, ties among weights can be broken consistently so that the greedy bases satisfy $B_1\subseteq\cdots\subseteq B_s$.

\begin{theorem}[Flag matroids on fans]\label{intro:fans}
The polynomials $C_i(u)$ are the total Chern polynomials of a formally filtered nef flag on $A^*(\Sigma)$ whenever $\Sigma$ is a simplicial Lefschetz fan admitting a strictly convex class, and on the intersection cohomology $\IH^*(\Sigma)$ of every rational projective fan $\Sigma$ (where the coefficients of $C_i(u)$ act by conewise polynomials).
\end{theorem}

We deduce the intersection cohomology statement by passing to a projective simplicial subdivision and descending along the decomposition of minimal extension sheaves (Theorem~\ref{thm:fan-flags}). Specializing $\Sigma$ to the ordinary and augmented Bergman fans of an arbitrary loopless base matroid $N$ on $E$---chosen independently of the flag $M_1\preceq\cdots\preceq M_s$---yields formally filtered nef flags formed by the restrictions of $\mathcal S_{M_i}^\vee$ on $A_{\rm ord}(N)$ and of $\mathcal Q^{\rm aug}_{M_i^*}$ on $A_{\rm aug}(N)$, where $M_i^*$ is the dual matroid (Corollary~\ref{cor:both-chow-rings}).

In Section~\ref{sec:applications}, we develop several combinatorial consequences that clarify the geometric and matroid-theoretic meaning of this Schubert positivity:
\begin{enumerate}
\item \emph{Concordance-sensitive inequalities.} Already for two marked ranks $(1,3)$ on $n\ge4$ elements, a rank-one quotient $N\preceq M$ of a rank-three matroid with $a=c_1(F_N)$ and $b=c_1(F_M)$ has $\mathfrak S_{3142}(N,M)=a^2(b-a)$, so Corollary~\ref{intro:volume} yields
\begin{equation}\label{eq:intro-concordance}
\int_{X_E}a^2 b\,D_1\cdots D_{n-4}\ge\int_{X_E}a^3 D_1\cdots D_{n-4}
\end{equation}
for all nef divisors $D_i$ on $X_E$, as well as nonnegativity and log-concavity for iterated finite differences of the mixed degrees $\int a^{d-j}b^j$. Here the quotient hypothesis is essential: for $N=U_{1,4}$ and the matroid $M$ on $\{1,2,3,4\}$ with sole basis $\{1,2,3\}$, the element $4$ is a loop in $M$ but not in $N$, so $N$ is not a quotient of $M$, and the left-hand side of \eqref{eq:intro-concordance} minus the right-hand side equals $-1$ even though both matroids are realizable. This separates multi-descent Schubert positivity from results valid for arbitrary tuples of matroids \cite{LP,LS} (Subsection~\ref{subsec:comparison}).

\item \emph{Virtual orbit closures on flag varieties.} When a flag matroid $\mathcal M$ is realizable, the closure $Z$ of its torus orbit in $T=\operatorname{Fl}(r_1,\ldots,r_s;|E|)$ has a Schubert-positive fundamental cycle $[Z]$, recovered as the leading term of the Riemann--Roch class of $[\OO_Z]\in K^0(T)$. For a general flag matroid $\mathcal M$, no subvariety $Z\subseteq T$ exists, yet Cameron, Dinu, Micha\l ek, and Seynnaeve \cite[Definition~84]{CDMS} constructed a combinatorial $K$-class $y(\mathcal M)\in K^0(T)$. We show that its Riemann--Roch image behaves in top dimension like a genuine orbit closure.
\end{enumerate}

\begin{proposition}[Leading Chow class of a flag matroid]\label{intro:chow}
The highest-dimensional nonzero homogeneous component of $\operatorname{ch}(y(\mathcal M))\operatorname{td}(T_T)\cap[T]$ lies in dimension $\dim\bigl(\sum_{i=1}^s P(M_i)\bigr)$, and its coefficients in the Schubert basis are nonnegative integers.
\end{proposition}

For a single matroid on a Grassmannian, Proposition~\ref{intro:chow} recovers \cite[Corollary~1.10]{LS}. Furthermore, pairing this leading cycle (and the subleading component) with a fixed Schubert class and powers of the Pl\"ucker divisors produces Lorentzian multidegree polynomials (Proposition~\ref{prop:orbit-rr} and Corollary~\ref{cor:orbit-multidegrees}).

Finally, on the augmented Chow ring $A_{\rm aug}(N)$ of a loopless rank-$d$ matroid $N$ with coordinate divisors $y_e$ ($e\in E$), the supports of derived Schubert intersection numbers have a direct structural meaning in terms of matroid morphisms. Let $\Gamma$ be a product of compatible Schubert classes of the flags on $A_{\rm aug}(N)$ from Corollary~\ref{cor:both-chow-rings}, translate all Chern roots by $u$ to write $\Gamma(u)=\sum_{j=0}^k G_j u^{k-j}$ with $k=\deg\Gamma$, and for $0\le j\le\min(k,d)$ set
\[
\mathcal B_j=\left\{I\subseteq E:\ |I|=d-j,\ \deg_N\!\left(G_j\prod_{e\in I}y_e\right)>0\right\},
\]
so that $\mathcal B_0$ is the set of bases of $N$.

\begin{theorem}[Intrinsic Higgs lifts]\label{intro:higgs}
Let $q$ be the largest index with $\mathcal B_q\ne\varnothing$. For every $0\le j\le q$, the family $\mathcal B_j$ is the set of bases of a rank-$(d-j)$ matroid $K_j$ on $E$. Moreover, $K_q$ is a quotient of $N$, and $K_0,K_1,\ldots,K_q$ are the Higgs lifts of $K_q$ relative to $N$, with rank functions
\[
\rk_{K_j}(A)=\min\{\rk_N(A),\,\rk_{K_q}(A)+q-j\}\qquad(A\subseteq E).
\]
\end{theorem}

Because $y_e^2=0$ in $A_{\rm aug}(N)$, each layer generates a multiaffine Lorentzian polynomial, and polarizing the auxiliary translation variable combines all layers into the basis generating polynomial of a single matroid on an enlarged ground set (Theorem~\ref{thm:derived-higgs}). Whereas Eur and Huh \cite{EH20} proved log-concavity for the Higgs lifts of a prescribed matroid morphism, Theorem~\ref{intro:higgs} shows that the terminal support $\mathcal B_q$ of a Schubert class intrinsically singles out a quotient $N\to K_q$ and determines all intermediate supports $\mathcal B_j$. As a concrete specialization on ordinary Chow rings and on the augmented Chow ring of a Boolean matroid, we also show that for any loopless matroid $M$ of rank $r\ge2$, weighting subsets $T\subseteq E$ by $\rk_M(T)-1$ produces Lorentzian reduced rank polynomials (Corollary~\ref{cor:reduced-rank-polynomial}).

\subsection*{Acknowledgments}
The author thanks Matt Larson and Weizhe Zheng for enlightening discussions.

\subsection*{AI disclosure}
The main ideas of the paper are due to the author, although AI was used in the exploration of proofs, literature search, and proofreading. The author takes full responsibility for all mathematical content in the paper.

\section{Preliminaries}\label{sec:prelim}

\subsection{Hodge--Lefschetz modules}
Following Cattani \cite[Definition~2.3]{Cattani} in the grading convention of Li and Zheng \cite[Definition~2.7]{LZ}, a \emph{Hodge--Lefschetz module} $(V,S)$ of weight $d$ consists of a finite-dimensional $\mathbb Z$-graded real vector space $V=\bigoplus_{l\in\mathbb Z}V_l$, a real bilinear form $S$ on $V$, and a bigrading
\[
V_l\otimes_{\R}\C=\bigoplus_{p+q=d-l}V^{p,q},\qquad \overline{V^{p,q}}=V^{q,p},
\]
such that $V^{p,q}$ and $V^{p',q'}$ are orthogonal under the $\C$-bilinear extension of $S$ unless $p+p'=q+q'=d$. We say that $V$ is \emph{effective} if $V^{p,q}=0$ whenever $\min(p,q)<0$. Effectivity will be imposed only where explicitly stated. For $i\ge0$, we write $\mathfrak o^{i,i}(V,S)$ for the space of real endomorphisms $L$ of $V$ of bidegree $(i,i)$ satisfying $S(Lu,v)=(-1)^iS(u,Lv)$ for all $u,v\in V$.

An operator $N\in\mathfrak o^{1,1}(V,S)$ \emph{polarizes} $(V,S)$ if, for every $l\ge0$, the power $N^l\colon V_l\to V_{-l}$ is an isomorphism and the Hermitian form
\[
(u,v)\longmapsto \ii^{\,q-p}S(u,N^l\overline v)
\]
is positive definite on the primitive subspace $V^{p,q}\cap\ker N^{l+1}$ whenever $p+q=d-l$. In particular, $S$ is nondegenerate. Throughout what follows, we fix a nonempty relatively open convex cone $\mathcal C\subseteq\mathfrak o^{1,1}(V,S)$ of pairwise commuting polarizing operators. By Cattani's mixed hard Lefschetz and mixed Hodge--Riemann theorems \cite[Theorems~4.2 and~4.4]{Cattani}, translated via the grading shift of \cite[Remark~2.8]{LZ} (after a Tate twist if necessary to place the finitely many nonzero bidegrees in the effective range), any tuple $N_0,\ldots,N_l\in\mathcal C$ with $\Omega=N_1\cdots N_l$ induces an isomorphism $\Omega\colon V_l\to V_{-l}$ such that $\ii^{\,q-p}S(u,\Omega\overline v)$ is positive definite on $V^{p,q}\cap\ker(N_0\Omega)$ for $p+q=d-l$.

When $Y$ is a smooth complex projective variety of dimension $d$, its cohomology carries a canonical Hodge--Lefschetz structure with $V_l=H^{d-l}(Y,\R)$, $V^{p,q}=H^{p,q}(Y)$, and bilinear pairing
\[
S(u,v)=(-1)^{j(j+1)/2}\int_Y uv\qquad\bigl(u\in H^j(Y,\R)\bigr),
\]
polarized by the real ample cone $\mathcal C=\operatorname{Amp}(Y)$ \cite[Example~2.10]{LZ}. In combinatorial and algebraic settings where only even degrees occur, one equips a graded real vector space $A^*=\bigoplus_{k\ge0}A^k$ with the diagonal bigrading $V^{k,k}=A^k\otimes_{\R}\C$, as in the Lefschetz modules of \cite[\S1]{AHL} and the projective bundle rings of \cite[\S3.1]{LP}. On any compact K\"ahler manifold, the Lefschetz $\mathfrak{sl}_2$-triples attached to K\"ahler classes generate a semisimple Lie algebra by Looijenga and Lunts \cite[Proposition~1.6]{LL97}. More generally, on an irreducible complex projective variety, Saito constructed a canonical polarizable pure Hodge structure on its intersection cohomology \cite[Corollary~5.4.7]{Saito88}, for which de Cataldo and Migliorini established a geometric proof of the hard Lefschetz theorem and Hodge--Riemann relations \cite[Theorem~2.2.3]{dCM}.

Polarized Hodge--Lefschetz modules likewise govern numerous structures in convex geometry, representation theory, and combinatorics (compare \cite[\S2]{LZ} and \cite[Example~1.1]{AHL}):
\begin{itemize}
\item In convex geometry, McMullen established hard Lefschetz and the Hodge--Riemann--Minkowski inequalities on the polytope algebra of a simple polytope \cite[Theorems~7.3 and~8.2]{McM}, recovered analytically via volume polynomials by Timorin \cite[Theorems~5.1.1 and~5.3.1]{Tim}. For arbitrary convex polytopes, Karu proved hard Lefschetz and the Hodge--Riemann relations on the intersection cohomology of normal fans \cite[Theorems~0.1 and~0.2]{Karu}, as well as their relative counterparts on subdivision multiplicity spaces \cite[Theorems~1.1 and~1.2]{Karu19}.
\item In Coxeter theory, Elias and Williamson proved hard Lefschetz and the Hodge--Riemann relations on reduced Soergel bimodules $B_w\otimes_R\R$ \cite[Theorems~1.3 and~1.4]{EW14}, together with relative versions on multiplicity spaces of $B_x\otimes_R B_y$ \cite[Theorems~1.2 and~3.4]{EW21}.
\item In matroid theory, the ordinary Chow ring of a matroid originates in the wonderful arrangement models of De Concini and Procesi \cite[\S\S4--5]{DCP} and the atomic-lattice Chow rings of Feichtner and Yuzvinsky \cite[Definition~3 and Theorem~3]{FY}, realized on the order complex subdivision of the Bergman fan due to Ardila and Klivans \cite[\S2]{AK} and related to nested-set fans by Feichtner and Sturmfels \cite[Theorem~4.1 and Corollary~4.3]{FeSt}. Adiprasito, Huh, and Katz proved Poincar\'e duality, hard Lefschetz, and the Hodge--Riemann relations for ordinary Chow rings of loopless matroids \cite[Theorems~1.4, 6.19, and~8.8]{AHK}. Braden, Huh, Matherne, Proudfoot, and Wang introduced the augmented Chow ring and used semi-small deletion decompositions to establish the K\"ahler package both on ordinary and augmented Chow rings \cite[\S1.2 and Theorems~1.2, 1.5, and~1.6]{BHMPW22} and on the intersection cohomology module of a matroid over its graded M\"obius algebra \cite[Theorem~1.6]{BHMPW26}. Related K\"ahler packages hold on conormal fans of matroids by Ardila, Denham, and Huh \cite[Theorem~5.28]{ADH} and on unimodular quasiprojective quasilinear tropical fans by Amini and Piquerez \cite{AP}.
\end{itemize}
Moreover, Amini, Huh, and Larson proved a general decomposition theorem for Lefschetz modules: upon restriction to any subalgebra generated by degree-one classes in $\overline{\mathcal C}$, each indecomposable summand remains polarizable up to a grading shift \cite[Theorems~1.4 and~1.5]{AHL}.

\subsection{Projective bundle modules}
Fix a commutative graded real subalgebra $A=\bigoplus_{i\ge0}A^i$ of endomorphisms of $V$ containing the identity $1_V$ in $A^0$ and the polarizing cone $\mathcal C$ in $A^1$, with $A^i\subseteq\mathfrak o^{i,i}(V,S)$ for all $i\ge0$ (it always suffices to take the subalgebra generated by $\mathcal C$ and the Chern operators under consideration). A \emph{formal bundle} $F$ of rank $e\ge0$ in $A$ consists of $c_0(F)=1_V$ together with operators $c_j(F)\in A^j$ for $1\le j\le e$. We set $c_j(F)=0$ for $j<0$ or $j>e$ and write $c(F;u)=\sum_{j=0}^e c_j(F)u^j$.

When $e>0$, let $\xi$ be an indeterminate of bidegree $(1,1)$ representing the relative hyperplane operator for rank-one quotients, and set
\[
P_F(z)=z^e-c_1(F)z^{e-1}+\cdots+(-1)^ec_e(F),\qquad W(F)=V[\xi]/P_F(\xi)V[\xi],
\]
where $V[\xi]$ denotes the $A[\xi]$-module $V\otimes_{\R}\R[\xi]$ (no ring structure on $V$ itself is needed). Endowed with the bigrading $W(F)^{p,q}=\bigoplus_{a=0}^{e-1}\xi^a V^{p-a,q-a}$ and the bilinear form
\[
T(\xi^a v,\xi^b w)=(-1)^a S(v,s_{a+b-e+1}w),
\]
where $s_j$ is the coefficient of $u^j$ in $c(F;-u)^{-1}$ for $j\ge0$ (equivalently, the coefficient of $\xi^{e-1}$ in the reduction of $\xi^{e-1+j}$ modulo $P_F(\xi)$) and $s_j=0$ for $j<0$, the pair $(W(F),T)$ forms a Hodge--Lefschetz module of weight $d+e-1$ \cite[Definition~8.1]{LZ}.

We call a positive-rank formal bundle $F$ \emph{formally nef} over $(V,\mathcal C)$ if the cone $\mathcal C+\R_{>0}\xi$ polarizes $W(F)$. Direct sums, tensor products, exterior and Schur powers, and real twists $F\langle H\rangle$ by $H\in A^1$ (which shifts formal Chern roots by $x\mapsto x+H$) are defined via the universal symmetric-polynomial identities in the Chern operators and require no splitting into rank-one operators.

\subsection{Formal quotient flags}
Fix positive integers $b_1,\ldots,b_s$ and set $r_j=b_1+\cdots+b_j$, with $r_0=0$. A \emph{formal quotient flag} $F_\bullet$ of marked ranks $r_1<\cdots<r_s$ in $A$ consists of formal bundles $G_1,\ldots,G_s$ of ranks $b_1,\ldots,b_s$, called its graded pieces, together with the total quotients $F_j$ of rank $r_j$ defined by
\[
c(F_j;u)=\prod_{i=1}^j c(G_i;u)\qquad(1\le j\le s).
\]
Equivalently, $F_\bullet$ is determined by its total quotients $F_1,\ldots,F_s$, provided that every successive ratio $c(F_i;u)/\allowbreak c(F_{i-1};u)$ has degree at most $b_i$. For a flag of quotient vector bundles $F_s\twoheadrightarrow\cdots\twoheadrightarrow F_1$ with $F_0=0$, the graded pieces are $G_i=\ker(F_i\to F_{i-1})$.

\subsection{Schubert classes and specialization}\label{subsec:specialization}
Let $T=\operatorname{Fl}(r_1,\ldots,r_s;n)$ be the partial flag variety parametrizing quotient flags of $\C^n$ of ranks $r_1<\cdots<r_s$. On a complete flag refinement, write $x_a$ for the first Chern class of the $a$th successive line quotient, so that the Schubert polynomial of Lascoux and Sch\"utzenberger \cite[Definition~(1.3)]{LS82} associated with the simple transposition $s_j=(j,j+1)$ is $\mathfrak S_{s_j}=x_1+\cdots+x_j$. A permutation $w$ is \emph{compatible} with the marked ranks $r_1,\ldots,r_s$ if every descent of $w$ belongs to $\{r_1,\ldots,r_s\}$. By the basis lemma in \cite[\S2]{Ful95}, as $w$ ranges over compatible permutations of length $k$, the polynomials $\mathfrak S_w$ form a basis for the space of homogeneous degree-$k$ polynomials symmetric in each variable block $(x_{r_{i-1}+1},\allowbreak\ldots,\allowbreak x_{r_i})$. We evaluate any compatible Schubert polynomial at a formal quotient flag $F_\bullet$ by substituting $c_j(G_i)$ for the $j$th elementary symmetric polynomial in block $i$. When $s=1$, this specialization recovers the Schur classes $s_\lambda(F_1)=\det(c_{\lambda_i+j-i}(F_1))$ in the convention of Fulton--Lazarsfeld \cite[Theorem~I]{FL83} and Ross--Toma \cite[Theorem~1.1]{RT}. By the dual Jacobi--Trudi identity, $s_\lambda(F_1)$ is the Schur polynomial of conjugate shape $\lambda'$ in the formal Chern roots.

Let $Q_s\twoheadrightarrow\cdots\twoheadrightarrow Q_1\twoheadrightarrow Q_0=0$ denote the universal quotient flag on $T$. By Borel's presentation \cite[Theorem~26.1(b)]{Borel}, written in the block form of \cite[\S\S2.3--2.4]{Buch}, eliminating the Chern classes of the rank-$(n-r_s)$ kernel $\ker(\C^n\otimes\OO_T\to Q_s)$ from the Whitney relation imposes only the vanishing of the coefficients of $u^j$ in $c(Q_s;u)^{-1}$ for $j>n-r_s$. Since $V$ is finite dimensional, $A^j=0$ for all sufficiently large $j$. Choosing $n$ large enough that $A^j=0$ for $j>n-r_s$ and that all permutations under consideration lie in $S_n$, we obtain a well-defined graded algebra homomorphism
\begin{equation}\label{eq:specialization}
\varphi\colon H^{2*}(T,\R)\longrightarrow A,\qquad c_j\bigl(\ker(Q_i\to Q_{i-1})\bigr)\longmapsto c_j(G_i),
\end{equation}
under which the geometric Schubert classes of the universal quotient flag \cite[Propositions~7.5 and~8.1]{Ful92} map to the specialized Schubert operators $\mathfrak S_w(F_\bullet)$.

\section{Formal filtered nefness and the total-quotient criterion}\label{sec:definition}

\subsection{Definition and statement of the criterion}
Let $F_\bullet$ be a formal quotient flag on $(V,S)$ with graded pieces $G_1,\ldots,G_s$. To test filtered positivity of $F_\bullet$, we pair it with auxiliary quotient flags on smooth complex projective varieties $B$, equipping $V_B=V\otimes H^*(B,\R)$ with the tensor-product Hodge--Lefschetz structure and polarizing cone $\mathcal C+\operatorname{Amp}(B)$ \cite[Proposition~2.11]{LZ} (when $B$ is a point, $\operatorname{Amp}(B)=\{0\}$). Let $W_1\twoheadrightarrow\cdots\twoheadrightarrow W_s$ be a flag of quotient vector bundles on $B$ with $W_1$ globally generated, where repeated stages and zero bundles at the end of the flag are permitted, and define the formal incidence bundle
\begin{equation}\label{eq:inc-test}
K(F_\bullet,W_\bullet)=\bigoplus_{i=1}^s G_i\boxtimes W_i
\end{equation}
on $V_B$. Explicitly, if $x_{ia}$ and $y_{ib}$ denote the formal Chern roots of $G_i$ and $W_i$, the elementary symmetric functions in the pairwise sums $x_{ia}+y_{ib}$ are separately symmetric in each alphabet and hence define commuting operators in $A\otimes H^{2*}(B,\R)$.

\begin{definition}\label{def:filtered}
A formal quotient flag $F_\bullet$ is \emph{formally filtered nef} over $(V,\mathcal C)$ if, for every smooth complex projective variety $B$ and every flag $W_1\twoheadrightarrow\cdots\twoheadrightarrow W_s$ on $B$ with $W_1$ globally generated such that $K=K(F_\bullet,W_\bullet)$ has positive rank, the cone
\[
\mathcal C+\operatorname{Amp}(B)+\R_{>0}\xi
\]
polarizes $W(K)$ over $V_B$.
\end{definition}

\begin{remark}[Necessity of all total quotients]\label{rem:filtered-quotients}
Taking $B$ to be a point with $W_i=\C$ for $i\le j$ and $W_i=0$ for $i>j$ yields $c(K;u)=c(F_j;u)$, so every positive-rank total quotient $F_j$ of a formally filtered nef flag is necessarily formally nef.

Formal nefness of the top bundle $F_s$ alone does not suffice, even on $\PP^1$. Recall that a vector bundle is called \emph{lef} if its tautological quotient line bundle is lef in the sense of de Cataldo and Migliorini \cite[Definitions~2.1.1 and~2.1.3]{dCM02}. For $E=\OO\oplus\OO(1)$ on $\PP^1$, the global sections of $\OO_{\PP(E)}(1)$ define the blowdown $\mathbb F_1\to\PP^2$, whose sole positive-dimensional fiber---the section corresponding to $E\twoheadrightarrow\OO$---lies over a point of codimension two. Hence the map is semismall and $E$ is lef (though its nef quotient $\OO$ is not lef, since $c_1(\OO)=0$ fails hard Lefschetz on $\PP^1$). Now $E'=\OO(2)\oplus\OO(-1)$ shares the same rank and Chern classes as the globally generated bundle $E$, so $W(E')=W(E)$ and $E'$ is formally nef on $\PP^1$. Nevertheless, the two-stage quotient flag $E'\twoheadrightarrow\OO(-1)$ fails the condition over a point at stage $1$: writing $h=c_1(\OO(1))$, the projective hyperplane operator of $\OO(-1)$ is $-h$, and $h+(-h)=0$ fails hard Lefschetz.
\end{remark}

\begin{theorem}[Total-quotient criterion]\label{thm:abstract-stages}
A formal quotient flag on $(V,\mathcal C)$ is formally filtered nef if and only if every total quotient of positive rank is formally nef.
\end{theorem}

We establish sufficiency first for one-stage flags against arbitrary nef auxiliary bundles, next for genuine vector bundles on a smooth projective variety by comparing two determinantal resolutions, and finally for general Hodge--Lefschetz modules in Subsection~\ref{subsec:abstract-stages}.

\subsection{Flags with one stage}

\begin{proposition}\label{prop:one-stage}
For a flag with a single stage, Definition~\ref{def:filtered} is equivalent to formal nefness of its bundle $F$. More generally, if $F$ is formally nef on $(V,\mathcal C)$ and $U$ is a nef vector bundle of positive rank on a smooth projective variety $B$, then $F\boxtimes U$ is formally nef on $V_B$.
\end{proposition}

\begin{proof}
Necessity follows by taking $B$ to be a point and $U=\C$. Conversely, write $r=\rk F$ and $m=\rk U$, and let $p\colon Y\to B$ be the $m$-fold fiber power of $\PP_B(U)$ over $B$, with nef tautological hyperplane classes $y_1,\ldots,y_m$ pulled back from the factors. Over $V_Y=V\otimes H^*(Y,\R)$, consider the formal direct sum $D=\bigoplus_{i=1}^m F\langle y_i\rangle$. Because translating the projective variable by $y_i$ preserves the top-coefficient Segre pushforward and hence the projective pairing, the projective bundle module of the $i$th summand identifies with $W(F)\otimes H^*(Y,\R)$ with tautological operator $\xi+y_i$. For any $N\in\mathcal C$, $a\in\operatorname{Amp}(Y)$, and $t>0$, the operators $N+t\xi$ and $a+ty_i$ polarize $W(F)$ and $H^*(Y,\R)$, respectively. By the tensor-product and direct-sum theorems \cite[Propositions~2.11 and~8.14]{LZ}, the cone $\mathcal C+\operatorname{Amp}(Y)+\R_{>0}z$ polarizes $W(D)$, where $z$ denotes the projective variable of $D$.

Now let the symmetric group $S_m$ act on $Y$ by permuting the $m$ factors and on $W(D)$ by fixing $z$, and set $\Delta=\det(y_i^{j-1})_{1\le i,j\le m}$ and $k=\binom{m}{2}$. Alternating the projective bundle basis $\prod_i y_i^{a_i}$ ($0\le a_i<m$) leaves only determinants with exponents $0,\ldots,m-1$, identifying the sign-alternating subspace of $H^*(Y,\R)$ with $\Delta\,p^*H^*(B,\R)$. Moreover, $e_j(y_1,\ldots,y_m)$ acts on this subspace as $p^*c_j(U)$: indeed, over the even subring generated by the Chern classes of $U$ (and hence over all of $H^*(B,\R)$ by base change), multiplication by $y$ on $M=H^*(B)[y]/(P_U(y))$ induces multiplication by $\det(u-T)=P_U(u)$ on the top exterior power of the standard basis via $\prod_i(u-y_i)$.

Since the coefficients of $\prod_{i=1}^m P_F(z-y_i)$ are symmetric in $y_1,\ldots,y_m$, they restrict on the sign subspace to those of $P_{F\boxtimes U}(z)$, so the free projective bases yield a linear isomorphism of modules up to a $(k,k)$ bidegree shift:
\[
W(F\boxtimes U)\xrightarrow{\ \sim\ } W(D)^{\rm sign},\qquad z^a v\longmapsto z^a\Delta\,p^*v.
\]
This map multiplies the prescribed bilinear pairing by $m!>0$: expanding $\Delta^2$ shows that $\prod_i y_i^{m-1}$ appears with coefficient $(-1)^k m!$ while every other monomial has some exponent strictly less than $m-1$ and pushes forward to zero, so $p_*(\Delta^2)=(-1)^k m!$. Moving the first factor of $\Delta$ across the signed pairing contributes another factor of $(-1)^k$, and the Segre pushforwards in $z$ match by the symmetric-coefficient identity above. Notice that the $(k,k)$ shift increases the weight by $2k$ while leaving the centered degree grading unchanged.

Finally, fix $N\in\mathcal C$, $h\in\operatorname{Amp}(B)$, and $t>0$, and choose $\epsilon>0$ small enough that $h-\epsilon c_1(U)$ is ample on $B$. Because $\sum_{i=1}^m y_i$ is nef and relatively ample for $p$, the class
\[
b=p^*\bigl(h-\epsilon c_1(U)\bigr)+\epsilon\sum_{i=1}^m y_i
\]
is ample on $Y$. Consequently, the $S_m$-invariant operator $N+b+tz$ polarizes $W(D)$ and restricts on $W(D)^{\rm sign}$ to $N+h+tz$. Since the sign projector commutes with $N+b+tz$ and is orthogonal with respect to the pairing, both hard Lefschetz and the positive definiteness of primitive pairings pass to $W(D)^{\rm sign}$, proving that $F\boxtimes U$ is formally nef. Specializing to globally generated $U$ settles Definition~\ref{def:filtered} for one-stage flags.
\end{proof}

\subsection{The criterion for vector bundles on a variety}\label{sec:stagecriterion}
Let $Y$ be a smooth complex projective variety, whose real ample cone polarizes $H^*(Y,\R)$ on each connected component. By the Whitney formula, a genuine quotient flag on $Y$ shares the total Chern polynomials of the partial direct sums $\bigoplus_{i\le j}G_i$ of its graded pieces. Formulating the geometric criterion directly in terms of bundles $G_1,\ldots,G_s$ will later allow us to apply it to direct sums of line bundles representing $K$-theoretic differences of matroid classes.

\begin{theorem}\label{thm:filtered}
Let $G_1,\ldots,G_s$ be vector bundles on a smooth complex projective variety $Y$, and set $E_j=G_1\oplus\cdots\oplus G_j$ for $1\le j\le s$. The following conditions are equivalent:
\begin{enumerate}
\item Every partial sum $E_j$ of positive rank is formally nef.
\item For every flag $W_1\twoheadrightarrow\cdots\twoheadrightarrow W_s$ of quotient vector bundles on $Y$ with $W_1$ globally generated, the bundle $\bigoplus_{i=1}^s G_i\otimes W_i$ is formally nef whenever it has positive rank.
\end{enumerate}
\end{theorem}

Setting $W_i=\OO_Y$ for $i\le j$ and $W_i=0$ for $i>j$ deduces (a) immediately from (b). For the converse, we work component by component on $Y$ and compare two resolutions of a determinantal locus obtained by forgetting the smallest quotient $W_s$: over each perverse support of the first projection, a second resolution factors the linear map through a Grassmannian quotient of $G_s$, reducing the length of the flag by one.

\subsubsection{Lefschetz operators on perverse supports}
For any smooth projective variety $Z$, let $\Pol(Z)\subseteq N^1(Z)_{\R}$ denote the connected component of the open cone of polarizing real divisor classes containing $\operatorname{Amp}(Z)$. By \cite[Proposition~3.11(a)]{LZ}, $\Pol(Z)$ is an open convex cone on which mixed hard Lefschetz and the mixed Hodge--Riemann relations hold \cite[Theorems~4.2 and~4.4]{Cattani}, and convexity implies $\Pol(Z)+\clpol{Z}\subseteq\Pol(Z)$.

If $F$ is a vector bundle of positive rank on $Y$ with relative hyperplane class $\xi$ on $\PP(F)\to Y$, then $F$ is formally nef if and only if $\xi\in\clpol{\PP(F)}$. Indeed, for any $h\in\operatorname{Amp}(Y)$, the class $h+\delta\xi$ is ample on $\PP(F)$ for sufficiently small $\delta>0$: if $\xi\in\clpol{\PP(F)}$, then $h+t\xi=(h+\delta\xi)+(t-\delta)\xi\in\Pol(\PP(F))$ for all $t>\delta$ and hence for all $t>0$. Conversely, if $F$ is formally nef, the polarizing ray $\{h+t\xi : t>0\}$ meets $\operatorname{Amp}(\PP(F))$ and therefore lies in $\Pol(\PP(F))$, so $\xi=\lim_{t\to\infty}(h+t\xi)/t\in\clpol{\PP(F)}$.

Suppose next that $a$ and $\eta$ belong to the closure of a common polarizing component of commuting $(1,1)$-operators on a Hodge--Lefschetz module and preserve both factors of a graded direct sum $U\oplus U'$. By a theorem of Cattani and Kaplan \cite[Theorem~3.3(i)]{CK}, in the module formulation of \cite[Proposition~2.14]{LZ}, the monodromy filtration of $a+t\eta$ is independent of $t>0$ and restricts to the monodromy filtration of $(a+t\eta)|_U$ on $U$. Since hard Lefschetz for a degree-two operator on $U$ with center $b$ is equivalent to the identity
\[
M_j=\bigoplus_{l\ge b-j}U^l
\]
for its zero-centered monodromy filtration, hard Lefschetz for $(a+t\eta)|_U$ at a single $t>0$ implies hard Lefschetz on $U$ with the same center for every $t>0$.

In what follows, $\IH^*(S)$ denotes the rational/real intersection cohomology of a projective variety $S$ with its standard cohomological grading.

\begin{lemma}[Support criterion]\label{lem:support}
Let $f\colon Z\to P$ be a morphism of projective varieties with $Z$ smooth, and let $a$ and $\eta$ be real divisor classes on $P$ with $a$ ample. Assume that every simple perverse constituent of $Rf_*\R_Z[\dim Z]$ is the intersection complex $\IC(S,\R)$ of its support $S$. Then $f^*\eta\in\clpol{Z}$ if and only if, for each such support $S$ and every $t>0$, the restriction $(a+t\eta)|_S$ satisfies hard Lefschetz on $\IH^*(S)$ with center $\dim S$.
\end{lemma}

\begin{proof}
Choose a decomposition-theorem splitting for $Rf_*\R_Z[\dim Z]$ in the complex projective form of Beilinson, Bernstein, Deligne, and Gabber \cite[Theorem~6.2.5]{BBD} and de Cataldo and Migliorini \cite[Theorem~2.1.1(b),(c)]{dCM}. By naturality of cup products \cite[Remark~4.4.3]{dCM}, any divisor class on $P$ acts on each summand $\IH^*(S)$ via restriction to $S$.

Any real ample class $b$ on $S$ satisfies hard Lefschetz on $\IH^*(S)$. Indeed, write $b=\sum_i c_i b_i$ with $c_i>0$ and each $b_i$ rational ample. Pulling back to a projective resolution $\pi\colon\widetilde S\to S$ realizes $\IH^*(S)$ as a direct summand of $H^*(\widetilde S,\R)$ whose complement is stable under classes from $S$. Hard Lefschetz holds when all $c_i=1$ by \cite[Theorem~2.2.3(c)]{dCM}, and it persists for general $c_i>0$ by constancy of the monodromy filtration on $H^*(\widetilde S,\R)$.

If $f^*\eta\in\clpol{Z}$, then $f^*a\in\clpol{Z}$ as well because $f^*a$ is nef. Since $(a+t\eta)|_S$ is ample for small $t>0$, it satisfies hard Lefschetz on $\IH^*(S)$ for small $t>0$, and hence for all $t>0$ by constancy of the monodromy filtration on the invariant summand $\IH^*(S)$ Notice that this forward implication requires only that $\IC(S,\R)$ occur as a direct summand of $Rf_*\R_Z[\dim Z]$.

Conversely, set $n=\dim Z$, let $M_{S,j}$ be the multiplicity space of $\IC(S,\R)$ in perverse degree $j$ so that $\IH^l(S)\otimes M_{S,j}$ contributes to $H^{l+n-\dim S+j}(Z,\R)$, and fix an integral ample class $v$ on $Z$. The perverse matrix entry of cup product with $v$ from degree $j$ to degree $j'$ is an $\operatorname{Ext}^{2+j-j'}$ morphism, which vanishes for $j'>j+2$ and acts on multiplicity spaces at $j'=j+2$ via a linear map $v_2$. By relative hard Lefschetz \cite[Theorem~6.2.10]{BBD}, \cite[Theorem~2.1.1(a)]{dCM}, the powers
\[
v_2^j\colon M_{S,-j}\xrightarrow{\ \sim\ }M_{S,j}
\]
are isomorphisms for all $j\ge0$. Rescaling basis vectors in perverse degree $j$ by $\epsilon^{j/2}$ multiplies any component of $\epsilon v$ raising perverse degree by $k$ by $\epsilon^{1-k/2}$. As $\epsilon\to0$, the operator $f^*(a+t\eta)+\epsilon v$ therefore converges on the $S$-summands to
\[
(a+t\eta)|_S\otimes 1+1\otimes v_2.
\]
Since $(a+t\eta)|_S$ is Lefschetz with center $\dim S$ for all $t\ge0$ (by hypothesis for $t>0$ and by ampleness of $a|_S$ at $t=0$) while $v_2$ is Lefschetz with center $0$, their diagonal $\mathfrak{sl}_2$-sum is Lefschetz with center $n$ after the degree shift $n-\dim S$. For any fixed $T>0$, continuity of the Lefschetz determinants on the compact interval $t\in[0,T]$ furnishes $0<\epsilon\le1$ such that $f^*(a+t\eta)+\epsilon v$ satisfies hard Lefschetz for all $t\in[0,T]$. Because $f^*a+\epsilon v$ is ample at $t=0$ and the primitive forms remain nondegenerate throughout $[0,T]$, their signatures stay positive definite at $t=T$, placing $f^*(a+T\eta)+\epsilon v$ in $\Pol(Z)$. Dividing by $T$ and letting $T\to\infty$ yields $f^*\eta\in\clpol{Z}$.
\end{proof}

The next descent lemma will also be used in Section~\ref{sec:matroids} to pull back matroid classes.

\begin{lemma}[Descent from a regular zero locus]\label{lem:section}
Let $p\colon X\to Y$ be a morphism of nonempty connected smooth projective varieties, and let $L\subseteq X$ be the smooth zero locus of a regular section of a globally generated vector bundle $J$ on $X$ such that $p_L=p|_L\colon L\to Y$ is birational. If $D$ is a formally nef vector bundle of positive rank on $X$ and $F$ is a vector bundle of positive rank on $Y$ such that $D|_L$ and $p_L^*F$ have the same rank and Chern classes, then $F$ is formally nef.
\end{lemma}

\begin{proof}
Write $\pi\colon Z=\PP(D)\to X$, let $\xi$ be its relative hyperplane class, and put $A=H^*(Z,\R)$. The fundamental class of $\PP(D|_L)$ in $Z$ is $\gamma=\pi^*c_{\rk J}(J)$. We claim first that the image of $\Pol(Z)$ polarizes the quotient module $A/\Ann(\gamma)$ equipped with the trace map $\deg_\gamma(\alpha)=\int_Z\gamma\alpha$. When $\rk J=0$, we have $\gamma=1$ and $A/\Ann(\gamma)=A$. When $\rk J>0$, the pullback $\pi^*J$ is globally generated on $Z$, so its relative hyperplane class $\xi_J$ on $\PP(\pi^*J)$ is nef and relatively ample. Hence $a+t\xi_J$ is ample for every $a\in\operatorname{Amp}(Z)$ and $t>0$. By cone extension \cite[Proposition~8.17]{LZ}, $a+t\xi_J$ polarizes $H^*(\PP(\pi^*J),\R)$ for all $a\in\Pol(Z)$ and $t>0$, and top Chern descent \cite[Theorem~8.10(a)]{LZ} then implies that the image of $\Pol(Z)$ polarizes $A/\Ann(\gamma)$.

Next set $R_F=H^*(\PP(F),\R)$ with monic projective relation $P_F(\xi_F)=0$. Because $D|_L$ and $p_L^*F$ share the same Chern classes, the pullback of $P_F(\xi)$ vanishes on $\PP(D|_L)$, so $\gamma P_F(\xi)=0$ in $A$ by the Gysin formula. Thus we obtain a well-defined graded algebra homomorphism
\[
\varphi\colon R_F\longrightarrow A/\Ann(\gamma),\qquad \sum_j a_j\xi_F^j\longmapsto \sum_j(p\pi)^*a_j\,\xi^j\bmod\Ann(\gamma).
\]
Putting $r=\rk F$, fiber integration along $\PP(D|_L)\to L$ together with birationality of $p_L$ yields
\[
\deg_\gamma\bigl(\varphi(w)\bigr)=\int_L p_L^*a_{r-1}=\int_Y a_{r-1}=\deg_{R_F}(w)
\]
for every reduced element $w=\sum_{i=0}^{r-1}a_i\xi_F^i\in R_F$. Hence $\varphi$ preserves the nondegenerate Poincar\'e pairings and embeds $R_F$ into $A/\Ann(\gamma)$ as an orthogonal direct summand invariant under $(p\pi)^*h$ and $\xi$ for any $h\in\operatorname{Amp}(Y)$. Since $(p\pi)^*h$ is nef on $Z$ and $\xi\in\clpol{Z}$ by formal nefness of $D$, their images in $A/\Ann(\gamma)$ belong to the closure of a common polarizing cone. As $h+t\xi_F$ is ample on $\PP(F)$ for small $t>0$, constancy of the monodromy filtration on the summand $R_F$ propagates hard Lefschetz and the ample primitive signatures to all $t>0$, showing that $h+t\xi_F\in\Pol(\PP(F))$ for every $t>0$.
\end{proof}

\subsubsection{Proof of the geometric criterion}
For a vector bundle $G$ on $Y$, write $\Gr_q(G)\to Y$ for the Grassmann bundle of rank-$q$ quotients of $G$.

\begin{lemma}\label{lem:grass}
Let $A$ and $G$ be vector bundles on a smooth projective variety $Y$, and for $0\le q\le\rk G$ let $\pi\colon\Gr_q(G)\to Y$ be the projection with universal rank-$q$ quotient $V_q$. If $A$ and $A\oplus G$ are formally nef whenever they have positive rank, then $\pi^*A\oplus V_q$ is formally nef whenever its rank is positive. Moreover, formal nefness is preserved under pullback along any Grassmann bundle or product projection.
\end{lemma}

\begin{proof}
Set $e=\rk G$. The cases $q=0$ and $q=e$ follow immediately from the hypotheses once the Grassmann-pullback assertion is checked, so suppose $0<q<e$ and consider the natural morphism
\[
\PP(\pi^*A\oplus V_q)\longrightarrow\PP(A\oplus G),
\]
under which the target hyperplane class $\xi$ pulls back to the source hyperplane class. Over $\PP(A)$, where the functional on $G$ vanishes, the fiber is $\Gr_q(\C^e)$. Over $\PP(A\oplus G)\setminus\PP(A)$, where the kernel of $G\twoheadrightarrow V_q$ must lie inside the hyperplane annihilated by the nonzero functional on $G$, the fiber is $\Gr_{q-1}(\C^{e-1})$. Over each of these two locally trivial strata, the morphism is a Grassmann bundle whose fiber cohomology is generated by restrictions of tautological Chern classes, so the cohomology local systems are constant and the only perverse supports in the decomposition theorem are $\PP(A\oplus G)$ and (when $\rk A>0$) $\PP(A)$. Since $\xi|_S\in\clpol{S}$ on both smooth supports $S$ by hypothesis, $(a+t\xi)|_S\in\Pol(S)$ for every ample class $a$ on $\PP(A\oplus G)$ and $t>0$, so Lemma~\ref{lem:support} implies formal nefness of $\pi^*A\oplus V_q$. Finally, base change along a Grassmann bundle or product projection induces a Grassmann bundle or product projection between projectivizations, to which Lemma~\ref{lem:support} applies with the full target as sole support.
\end{proof}

Choosing generating global sections of $W_1$ classifies any auxiliary quotient flag on $Y$ via a morphism to a partial flag variety. The next graph lemma pulls formal nefness back along that classifying morphism.

\begin{lemma}\label{lem:graph}
Let $T$ be a partial flag variety of quotients of a finite-dimensional vector space $H_0$, and let $b\colon Y\to T$ be a morphism from a smooth projective variety $Y$. If $D$ is a formally nef bundle of positive rank on $Y\times T$, then its pullback along $(\mathrm{id}_Y,b)\colon Y\to Y\times T$ is formally nef.
\end{lemma}

\begin{proof}
Working on each connected component of $Y$, let $J$ be the globally generated quotient of $\operatorname{End}(H_0)\otimes\OO_{T\times T}$ by the subbundle of endomorphisms carrying each kernel of the first universal flag into the corresponding kernel of the second. The section of $J$ induced by $\mathrm{id}_{H_0}$ vanishes precisely along the diagonal $\Delta_T$, where $J|_{(t,t)}$ identifies via the infinitesimal $\operatorname{GL}(H_0)$-action with $\mathfrak{gl}(H_0)/\mathfrak p_t\simeq T_tT$ and the differential along the second factor is an isomorphism. Pulling back along $b\times\mathrm{id}_T\colon Y\times T\to T\times T$ therefore yields a regular section cutting out the graph of $b$ smoothly, and since the projection from the graph to $Y$ is an isomorphism, Lemma~\ref{lem:section} concludes the proof.
\end{proof}

\begin{proof}[Proof of Theorem~\ref{thm:filtered}]
Assume (a), work on a connected component of $Y$, and choose a finite-dimensional vector space $H_0$ generating $W_1$. For $0\le i\le s$, let $T_i$ be the partial flag variety of quotients of $H_0$ having the first $i$ prescribed ranks, denoting its universal quotients again by $W_i$ (with $T_0$ a point and $W_0=H_0\otimes\OO$). To streamline the induction on $s$, we allow an initial summand $G_0\otimes H_0$, setting $E_i=\bigoplus_{a=0}^i G_a$ and specializing to $G_0=0$ at the end. When $s=0$, the bundle $G_0\otimes H_0$ is a direct sum of copies of $G_0$ and is formally nef by \cite[Proposition~8.14]{LZ}.

Assume now that $s\ge1$ and that the universal statement holds for flags of length $s-1$ over every smooth projective base. Suppressing pullbacks on $Y\times T_{s-1}$, write
\[
C=\bigoplus_{i=0}^{s-1} G_i\otimes W_i,\qquad G=G_s,\qquad H=W_{s-1}.
\]
Dualizing the final rank-$j$ universal quotient $H\twoheadrightarrow W_s$ embeds $W_s^\vee\subseteq H^\vee$ as a rank-$j$ subbundle, inducing a projective morphism
\[
f\colon \PP\bigl(C\oplus(G\otimes W_s)\bigr)\longrightarrow P:=\PP\bigl(C\oplus(G\otimes H)\bigr)
\]
that forgets $W_s^\vee$ and pulls the hyperplane class $\xi$ on $P$ back to the hyperplane class on the source.

Viewing a point of $P$ as a nonzero pair consisting of a functional on $C$ and a linear map $\varphi\colon G\to H^\vee$ up to scale, its fiber under $f$ parametrizes subspaces $W_s^\vee\subseteq H^\vee$ containing $\operatorname{im}\varphi$. For $0\le q\le\min(j,\rk G)$, let $D_q\subseteq P$ be the determinantal locus $\{\rk\varphi\le q\}$. Local matrix normal forms show that the fixed-rank strata are smooth and connected. On the open stratum of $D_q$, the image $\operatorname{im}\varphi$ is a rank-$q$ subbundle of $H^\vee$ and $f$ is a Grassmann bundle whose fiber cohomology is generated by tautological Chern classes. By the decomposition theorem, every simple perverse constituent of $Rf_*\R[\dim]$ is therefore $\IC(D_q,\R)$ for some $0\le q\le\min(j,\rk G)$ with $D_q\ne\varnothing$.

For each such $q$, form the alternative resolution by choosing a rank-$q$ quotient $G\twoheadrightarrow V_q$ through which $\varphi$ factors: setting $Y_q=\Gr_q(G)$, we obtain a projective morphism
\[
f_q\colon \PP\bigl(C\oplus(V_q\otimes H)\bigr)\longrightarrow D_q
\]
from a smooth projective bundle over $Y_q\times T_{s-1}$ that is an isomorphism over the dense open stratum $\{\rk\varphi=q\}$ (sending $\varphi$ to $G\twoheadrightarrow G/\ker\varphi$, and equal to $\mathrm{id}_{\PP(C)}$ when $q=0$). Grouping the two $W_{s-1}$-summands rewrites the source bundle of $f_q$ as
\[
\bigoplus_{i=0}^{s-2}(G_i\otimes W_i)\oplus\bigl((G_{s-1}\oplus V_q)\otimes W_{s-1}\bigr),
\]
which is the universal bundle of a flag of length $s-1$ over $Y_q\times T_{s-1}$ (for $s=1$, this leaves the single stage $(G_0\oplus V_q)\otimes H_0$). By Lemma~\ref{lem:grass}, the partial sums $E_0,\ldots,E_{s-2}$ remain formally nef upon pullback to $Y_q$, as does the final partial sum $E_{s-1}\oplus V_q$ since both $E_{s-1}$ and $E_s=E_{s-1}\oplus G_s$ are formally nef whenever their ranks are positive. The induction hypothesis therefore places $f_q^*\xi$ in the closed polarizing cone of the source of $f_q$. Since $f_q$ is birational onto $D_q$, $\IC(D_q,\R)$ appears as a summand of its direct image, and the forward implication of Lemma~\ref{lem:support} shows that $(a+t\xi)|_{D_q}$ satisfies hard Lefschetz on $\IH^*(D_q)$ for every ample $a$ on $P$ and $t>0$. By the converse implication of Lemma~\ref{lem:support} applied to $f$, the universal bundle on $Y\times T_s$ is formally nef. Pulling back along the classifying morphism $b\colon Y\to T_s$ via Lemma~\ref{lem:graph} then yields (b).
\end{proof}

Applying Theorem~\ref{thm:filtered} on $Y\times B$ together with product pullback (Lemma~\ref{lem:grass}) shows that any genuine flag on $Y$ with formally nef positive-rank total quotients is formally filtered nef in the sense of Definition~\ref{def:filtered}.

\subsection{Proof of the total-quotient criterion on a module}\label{subsec:abstract-stages}
We now establish sufficiency in Theorem~\ref{thm:abstract-stages} for an arbitrary Hodge--Lefschetz module $(V,S)$ assuming that every positive-rank total quotient of $F_\bullet$ is formally nef. The strategy is to build the determinantal resolutions of Theorem~\ref{thm:filtered} universally over a product of Grassmannians and specialize universal Chern classes to $A$ (using complex dimensions for fibers and ordinary cohomological degrees for fiber cohomology).

\subsubsection{Specializing universal families}
Let $S=\prod_i\Gr_{b_i}(\C^L)$ be the product of Grassmannians of quotients of ranks $b_i=\rk G_i$ (including a factor for $G_0$ when the induction employs a nonzero initial piece), and set $R=H^*(S,\R)$. Because complementary Chern relations on $\Gr_{b_i}(\C^L)$ start in Chern degree $L-b_i+1$ whereas $A$ vanishes in high degrees, taking $L\gg0$ yields a graded algebra homomorphism $\rho\colon R\to A$ sending the $i$th universal Chern classes to $c_j(G_i)$. All projective bundles, Grassmann bundles, and determinantal varieties in the proof of Theorem~\ref{thm:filtered} then assemble into locally trivial associated families over $S$, which we denote by a subscript $S$ (writing the corresponding letter without subscript for a fiber over a point).

For each smooth projective family $\pi\colon X_S\to S$ arising in the construction---a tower of Grassmann and projective bundles together with a fixed auxiliary flag variety---its fiber cohomology lies in types $(k,k)$ and $H^*(X_S,\R)$ is free over $R$ on homogeneous lifts of a fiber basis. Setting
\[
V[X]=H^*(X_S,\R)\otimes_R V,\qquad S_X(a\otimes u,b\otimes v)=(-1)^{\deg a/2}S\bigl(u,\rho(\pi_*(ab))v\bigr),
\]
produces a balanced Hodge--Lefschetz module of weight $d+\dim X$ whose projective-bundle instances coincide (compatibly in towers) with the formal projective bundle modules $W(F)$ for the specialized Chern operators.

For $N\in\mathcal C$ and any divisor $a$ on $X_S$ that is ample on the fiber $X$, the operator $N+\delta a$ polarizes $V[X]$ for all sufficiently small $\delta>0$: rescaling homogeneous basis lifts of fiber degree $i$ by $\delta^{i/2}$ ensures that cup product with $a$ (which raises fiber degree by at most two because its $R$-matrix entries have nonnegative degree) sends $N+\delta a$ to $N\otimes 1+1\otimes a|_X$ and $\delta^{-\dim X}S_X$ to the signed tensor-product pairing as $\delta\to0$, where polarization holds by \cite[Proposition~2.11]{LZ} and persists for small $\delta>0$ by openness. Varying $N$ and $a$ along compact segments shows that all such small perturbations lie in a single open polarizing component $P(X)$ inside the span of $\mathcal C$ and the real divisor classes of $X_S$, and $P(X)$ is convex by \cite[Proposition~3.11(a)]{LZ}.

\subsubsection{Two-scale rescaling for fiber-ample pullbacks}
Let $f_S\colon X_S\to P_S$ be one of the projective morphisms of associated families with $X_S$ smooth, whose shifted direct image splits $R$-linearly into shifts of constant-coefficient intersection complexes $\IC(D_S)$ on associated supports $D_S$, and set
\[
I_V(D)=\IH^*(D_S,\R)\otimes_R V.
\]
We check first that $\IH^*(D_S,\R)$ is free over $R$ with $R/R^{>0}\otimes_R\IH^*(D_S,\R)\simeq\IH^*(D)$. Applying the decomposition theorem to $p\colon D_S\to S$ with $\IC(D_S)$-coefficients \cite[Theorem~6.2.5]{BBD}, \cite[Theorem~1.6.1]{dCMsurvey}, local triviality and $i_D^*\IC(D_S)[-\dim S]\simeq\IC(D)$ along a fiber inclusion $i_D\colon D\hookrightarrow D_S$ identify the stalks of $\mathcal H^{l-\dim D}(Rp_*\IC(D_S)[-\dim S])$ with $\IH^l(D)$. Since $S$ is simply connected, these local systems are constant, splitting $\IH^*(D_S,\R)$ into a free graded $R$-module reducing to $\IH^*(D)$ modulo $R^{>0}$. Likewise, proper base change along $i_P\colon P\hookrightarrow P_S$ identifies $i_P^*Rf_{S*}\R_{X_S}[\dim X_S][-\dim S]$ with $Rf_*\R_X[\dim X]$ while preserving perverse degrees, and all splittings take place over $R$ prior to tensoring with $V$.

Now choose rational divisors $a$ on $P_S$ (ample on $P$) and $v$ on $X_S$ (for example ample on $X_S$) such that $f^*a+\theta v$ is ample on $X$ for small $\theta>0$. We claim that for fixed $N\in\mathcal C$,
\begin{equation}\label{eq:specialized-two-scales}
N+\delta f^*a+\delta\theta v\in P(X)\qquad(0<\delta,\theta\ll1),
\end{equation}
with independent upper bounds on $\delta$ and $\theta$. Indeed, choose homogeneous $R$-basis lifts of $\IH^l(D)$ in each perverse summand of degree $j$ (contributing in total fiber degree $i=l+\dim X-\dim D+j$) tensor multiplicity vectors, and rescale lifts of degrees $(i,j)$ by $\delta^{i/2}\theta^{j/2}$. Because $f^*a$ preserves $j$, $v$ raises $j$ by at most two, and both raise $i$ by at most two over $R$, entries of $\delta\theta v$ acquire factors $\delta^{1-(i'-i)/2}\theta^{1-(j'-j)/2}$ with nonnegative exponents, so $N+\delta f^*a+\delta\theta v$ extends continuously to $(\delta,\theta)=(0,0)$ with limit $N+a|_D+v_2$ on each support. Since $N$, $a|_D$ on $\IH^*(D)$, and $v_2$ on the multiplicity spaces each satisfy hard Lefschetz \cite[Theorem~6.2.10]{BBD}, \cite[Theorem~2.1.1(a)]{dCM}, \cite[Theorem~1.6.3]{dCMsurvey}, their tensor-product $\mathfrak{sl}_2$-sum satisfies hard Lefschetz with center $d+\dim X$ in an open rectangle $0<\delta<\delta_0$, $0<\theta<\theta_0$. Fixing $\theta$ in that rectangle, our single-scale perturbation places the operator in $P(X)$ for small $\delta>0$, and constancy of primitive signatures across the rectangle proves \eqref{eq:specialized-two-scales}. Letting $\theta\to0$ shows that $N+\delta f^*a\in\overline{P(X)}$ for $0<\delta<\delta_0$, while single-scale rescaling on $I_V(D)$ shows that $N+\delta a$ satisfies hard Lefschetz on $I_V(D)$ with center $d+\dim D$.

\subsubsection{Comparison of supports}
Suppose $\xi$ is a divisor on $P_S$ such that, for every support $D_S$ of $f_S$, there exists a proper associated morphism $g_S\colon X'_S\to P_S$ with $X'_S$ smooth, constant-coefficient associated supports as in \eqref{eq:specialized-two-scales}, $\IC(D_S)$ occurring as a summand of its shifted direct image, and $g^*\xi\in\overline{P(X')}$. We deduce that
\begin{equation}\label{eq:specialized-support}
f^*\xi\in\overline{P(X)}.
\end{equation}
Fix $N\in\mathcal C$ and choose $\delta>0$ small enough for all finitely many two-scale bounds. On $V[X']$, both $N+\delta g^*a$ and $g^*\xi$ lie in $\overline{P(X')}$ and preserve the $R$-linear direct summand $I_V(D)$ together with its complement. Since $N+\delta a+t\xi$ satisfies hard Lefschetz on $I_V(D)$ for small $t>0$, constancy of the monodromy filtration \cite[Theorem~3.3(i)]{CK}, \cite[Proposition~2.14]{LZ} extends hard Lefschetz on $I_V(D)$ to all $t>0$. For any $T>0$, rescaling $V[X]$ by $\epsilon^{j/2}$ in perverse degree $j$ alone sends $N+\delta f^*a+t f^*\xi+\epsilon v$ as $\epsilon\to0$ to $(N+\delta a+t\xi)\otimes1+1\otimes v_2$, which is Lefschetz throughout $t\in[0,T]$. Choosing $0<\epsilon\le1$ small enough that hard Lefschetz holds on $[0,T]$ and \eqref{eq:specialized-two-scales} applies at $t=0$, constancy of primitive signatures yields $N+\delta f^*a+T f^*\xi+\epsilon v\in P(X)$. Dividing by $T$ and letting $T\to\infty$ proves \eqref{eq:specialized-support}.

\subsubsection{Inductive step and graph descent}
Carry out the Grassmann and determinantal constructions of Theorem~\ref{thm:filtered} over $S$ (including an initial summand $G_0\otimes H_0$, with $s=0$ covered by \cite[Proposition~8.14]{LZ}). For $s\ge1$, let $R_q=\Gr_q(\C^{b_s})$ with universal quotient $V_q$, set $V'=V[R_q]$, fix $N\in\mathcal C$, and let $\eta$ be an integral ample divisor on $R_q$ (taking $V'=V$, $\eta=0$, and omitting rank-zero modules when $q\in\{0,b_s\}$). Applying \eqref{eq:specialized-support} to the projections $\PP(E_i)|_{R_q}\to\PP(E_i)$ for $i\le s-1$ (with the identity comparison map) and to $\PP(E_{s-1}\oplus V_q)\to\PP(E_{s-1}\oplus G_s)$ (with the inclusions of the smooth supports $\PP(E_{s-1}\oplus G_s)$ and $\PP(E_{s-1})$ as comparison maps, where each hyperplane class $\xi_i$ lies in the closed distinguished component by formal nefness), and then applying \eqref{eq:specialized-two-scales} along the projection to the Grassmann base, we obtain a common $\delta>0$ such that $N'=N+\delta\eta$ makes $E_0,\ldots,E_{s-2}$ and $E_{s-1}\oplus V_q$ formally nef over $V'$.

By the induction hypothesis over $(V',\R_{>0}N')$, the shorter flag ending in $G_{s-1}\oplus V_q$---whose universal projective module over $T_{s-1}$ agrees with $V[X_q]$ by composition of fiber integrations and universal Grassmann classification---has polarizers $N'+h_{s-1}+t\xi$ ($h_{s-1}\in\operatorname{Amp}(T_{s-1})$, $t>0$). Picking an integral ample $a_{s-1}$ on $T_{s-1}$ and small $\theta_0>0$ so that $v=a_{s-1}+\theta_0\xi$ is fiber-ample for $X_q\to R_q$ (or $v=\xi$ when $T_{s-1}$ is a point), \eqref{eq:specialized-two-scales} places $N+\delta\eta+\epsilon v$ in $P(X_q)$ for small $\epsilon>0$ (shrinking $\delta$ if needed), identifying the induction cone with $P(X_q)$ and yielding $\xi\in\overline{P(X_q)}$. Hence \eqref{eq:specialized-support} for $f_S$ implies $\xi\in\overline{P(X)}$ over $T_s$. Moreover, for any $h\in\operatorname{Amp}(T_s)$ and $u>0$, small projective perturbation \cite[Lemma~8.2]{LZ} polarizes $N+uh+t\xi$ for small $t>0$ uniformly on compact $u$-intervals. Matching with $P(X)$ at small $u$ via \eqref{eq:specialized-two-scales} shows that $N+h+t\xi\in P(X)$ for $u=1$ and all $t>0$.

Finally, for an arbitrary generated flag on $B$ classified by $b\colon B\to T_s$ (when $T_s$ is not already a point), pull the universal incidence bundle $D$ back to $V\otimes H^*(B\times T_s,\R)$. Since the bundle $J$ cutting out the graph of $b$ in Lemma~\ref{lem:graph} is globally generated and hence formally nef over $V\otimes H^*(B\times T_s,\R)$, simultaneous projectivization and top Chern descent \cite[Proposition~8.6 and Theorem~8.10]{LZ} polarize the image of $\gamma=c_{\rk J}(J)$ on $W(D)$. Because multiplication by $\gamma$ factors into surjective graph restriction followed by injective graph Gysin pushforward, $\gamma W(D)$ identifies isometrically with $W(b^*D)$, and writing any $h_B\in\operatorname{Amp}(B)$ as $(h_B-\epsilon b^*h_T)+\epsilon b^*h_T$ for $h_T\in\operatorname{Amp}(T_s)$ and small $\epsilon>0$ exhibits $N+h_B+t\xi$ among the descended polarizers, establishing Theorem~\ref{thm:abstract-stages}.

\subsection{Incidence modules and relative primitive parts}
Definition~\ref{def:filtered} extends naturally to finitely many commuting flags at once.

\begin{proposition}\label{prop:inc-module}
Let $F_\bullet^{(1)},\ldots,F_\bullet^{(a)}$ be formally filtered nef flags in the same commuting operator algebra $A$ on $(V,\mathcal C)$. For each $1\le i\le a$, choose a quotient flag on a common smooth projective variety $B$ with globally generated first stage such that the incidence bundle $K_i$ has rank $e_i>0$, and set $\Delta=\prod_{i=1}^a c_{e_i}(K_i)$. Then the multi-projective bundle module of $K_1,\ldots,K_a$ over $V_B$ is polarized by
\[
\mathcal C+\operatorname{Amp}(B)+\sum_{i=1}^a\R_{>0}\xi_i,
\]
and the image $\Delta V_B$ is a Hodge--Lefschetz module of weight $d+\dim B-\sum_{i=1}^a e_i$ polarized by $\mathcal C+\operatorname{Amp}(B)$ with pairing $S_\Delta(\Delta u,\Delta v)=S_B(u,\Delta v)$.
\end{proposition}

\begin{proof}
Each single projective bundle module $W(K_i)$ over $V_B$ is polarized by $\mathcal C+\operatorname{Amp}(B)+\R_{>0}\xi_i$ by Definition~\ref{def:filtered}, and all Chern operators commute by construction. Simultaneous projectivization and iterated top Chern descent \cite[Proposition~8.6 and Corollary~8.12]{LZ} (generalizing the single-divisor descent of Cattani, Kaplan, and Schmid \cite[Lemma~1.16]{CKS}) therefore polarize both the multi-projective bundle module and $\Delta V_B$, equipped with bigrading $(\Delta V_B)^{p,q}=\Delta V_B^{p,q}$ and nondegenerate pairing $S_\Delta$ (any rank-zero factor may be omitted since $c_0(K_i)=1$).
\end{proof}

Next let $V=\bigoplus_{j=0}^d V^j$ be a Lefschetz module in the diagonal grading $V^j=V^{j,j}_{\R}$, and set $Q(x,y)=(-1)^jS(x,y)$ for $x\in V^j$, so multiplication by elements of $A^1$ is self-adjoint for $Q$. For any nef operator $\ell\in\overline{\mathcal C}$ and $0\le k\le d$, define
\[
P_k(V,\ell)=\frac{\ker\ell^{k+1}}{\ker\ell^k+\ell\ker\ell^{k+2}},\qquad Q_k([x],[y])=Q(x,\ell^k y).
\]
In any homogeneous Jordan chain $b,\ell b,\ldots,\ell^s b$, only the initial vector $b$ can survive in $P_k(V,\ell)$, and it does so precisely when $s=k$. Hence $P_k(V,\ell)$ is the degree-$k$ primitive subquotient for the zero-centered monodromy filtration of $\ell$, and $Q_k$ is well defined by self-adjointness of $\ell$.

\begin{proposition}\label{prop:relative-primitive-flag}
If $F_\bullet$ is a formally filtered nef flag on $V$, then its Chern operators induce a formally filtered nef flag on each nonzero $P_k(V,\ell)$ (of weight $d-k$) with respect to the image of $\mathcal C$.
\end{proposition}

\begin{proof}
Because the Chern operators commute with $\ell$, they preserve the numerator and denominator of $P_k(V,\ell)$ and satisfy the same Whitney relations on the subquotient (when $\ell=0$, only $P_0(V,0)=V$ is nonzero). Working over the subalgebra generated by $\operatorname{span}(\mathcal C)$ (and $\xi$ on projective bundle modules), a chain starting in degree $j$ of length $k+1$ has perverse index $2j+k$ in the convention of Amini, Huh, and Larson, identifying $(P_k(V,\ell),Q_k)$ with $\operatorname{Prim}^{\bullet,[k]}_\ell$ and showing that the image of $\mathcal C$ polarizes $P_k(V,\ell)$ \cite[Corollary~1.11]{AHL}.

For any positive-rank total quotient $F$ of rank $e$ with $W=W(F)$, the operator $\ell$ acts coefficientwise in the free basis $1,\xi,\ldots,\xi^{e-1}$, identifying $P_k(W,\ell)$ with the projective bundle module of the induced Chern operators over $P_k(V,\ell)$ and matching both pairings via $Q_k(\xi^a[x],\xi^b[y])=Q(x,\ell^k s_{a+b-e+1}(F)y)$. Applying \cite[Corollary~1.11]{AHL} to $(W,\mathcal C+\R_{>0}\xi)$ shows that $h+t\xi$ polarizes $P_k(W,\ell)$ for all $h\in\mathcal C$ and $t>0$, so Theorem~\ref{thm:abstract-stages} yields formal filtered nefness of the induced flag.
\end{proof}

Similarly, applying divisor descent \cite[Lemma~2.18]{LZ} to $V$ and coefficientwise to each $W(F)$ (compare Karu \cite[Theorem~1.4]{Karu19} for convex functions on complete fans) shows that $F_\bullet$ induces a formally filtered nef flag on $\ell V$ of weight $d-1$ with pairing $Q_\ell(\ell x,\ell y)=Q(x,\ell y)$, and hence by iteration on any nonzero image of a product of nef divisors.

\section{Schubert functors preserve formal nefness}\label{sec:transfer}

Let $F_\bullet$ be a formally filtered nef flag of marked ranks $r_1<\cdots<r_s$ on $(V,\mathcal C)$, choose $n>r_s$ large enough for the specialization homomorphism $\varphi\colon H^{2*}(T,\R)\to A$ of \eqref{eq:specialization} on $T=\operatorname{Fl}(r_1,\ldots,r_s;n)$ with universal quotient flag $Q_s\twoheadrightarrow\cdots\twoheadrightarrow Q_1$, and write $\varphi(U)$ for the formal bundle on $V$ associated with a vector bundle $U$ on $T$. Denote by $\mathbf S'_w$ the quotient Schubert functor of Kra\'skiewicz and Pragacz \cite[\S7]{PragaczKP}.

\begin{theorem}[Quotient Schubert functors]\label{thm:schubert-preservation}\label{cor:kp-transfer}
For every partition $\lambda$ and every $1\le i\le s$, the specialized Schur bundle $\varphi(\mathbf S_\lambda(Q_i))$ is formally nef whenever it has positive rank. Likewise, $\varphi(\mathbf S'_w(Q_\bullet))$ is formally nef whenever it has positive rank on a complete quotient flag, and on a partial quotient flag whenever the Schubert bundle on a complete refinement descends to $T$.

More generally, if $U_\ell\twoheadrightarrow\cdots\twoheadrightarrow U_1$ is any flag of quotient vector bundles on $T$ with $U_\ell$ nef of positive rank, then $\varphi(U_1),\ldots,\varphi(U_\ell)$ form the total quotients of a formally filtered nef flag on $(V,\mathcal C)$ after omitting zero bundles and repeated ranks.
\end{theorem}

For example, because every surjection $Q_{i+1}\twoheadrightarrow Q_i$ splits locally, applying $\mathbf S_\lambda$ yields surjections $\mathbf S_\lambda(Q_{s})\twoheadrightarrow\cdots\twoheadrightarrow\mathbf S_\lambda(Q_1)$ with $\mathbf S_\lambda(Q_s)$ globally generated, so the positive-rank specializations $\varphi(\mathbf S_\lambda(Q_i))$ form a formally filtered nef flag once repeated stages are omitted.

Specializing to a one-stage flag shows that Schur functors---and more generally arbitrary nonzero polynomial functors---preserve formal nefness for an individual formal bundle. For a formal bundle $F$ of rank $e>0$ on $(V,\mathcal C)$ and a partition $\lambda$, write $\mathbf S_\lambda(F)$ for the formal bundle obtained by evaluating the universal rank-$e$ Chern polynomials of $\mathbf S_\lambda$ at $c_1(F),\ldots,c_e(F)$.

\begin{corollary}[Schur functors of formal bundles]\label{cor:schur-formal}
Let $F$ be a formally nef bundle of rank $e>0$ on $(V,\mathcal C)$. For every partition $\lambda$ such that $\mathbf S_\lambda(\C^e)\ne0$, the formal Schur bundle $\mathbf S_\lambda(F)$ is formally nef on $(V,\mathcal C)$. In particular, all tensor powers $F^{\otimes k}$, symmetric powers $\operatorname{Sym}^k F$ ($k\ge0$), and exterior powers $\bigwedge^k F$ ($0\le k\le e$) of $F$ are formally nef.
\end{corollary}

\begin{proof}
By Proposition~\ref{prop:one-stage}, the one-stage flag consisting of $F$ alone is formally filtered nef on $(V,\mathcal C)$. Taking $s=1$ and $T=\Gr_e(\C^n)$ ($n\gg0$) with universal rank-$e$ quotient $Q$, both $\mathbf S_\lambda(Q)$ and $Q^{\otimes k}$ are globally generated vector bundles on $T$ with specializations $\varphi(\mathbf S_\lambda(Q))=\mathbf S_\lambda(F)$ and $\varphi(Q^{\otimes k})=F^{\otimes k}$. The conclusion therefore follows from Theorem~\ref{thm:schubert-preservation}, using only Proposition~\ref{prop:one-stage} and Grassmannian specialization in Theorem~\ref{thm:flag-specialization}.
\end{proof}

\subsection{Specialization from a partial flag variety}
We prove first that arbitrary nef bundles on $T$ specialize to formally nef bundles over $V$ by cutting out the diagonal $\Delta_T\subseteq T\times T$ of dimension $m=\dim T$. Setting $Q_0=0$ and $S_i=\ker(\C^n\otimes\OO_T\to Q_i)$, the duals $S_1^\vee\twoheadrightarrow\cdots\twoheadrightarrow S_s^\vee$ form a quotient flag on $T$ with $S_1^\vee$ globally generated.

\begin{lemma}\label{lem:partial-flag-diagonal}
On $T\times T$, there exists a globally generated vector bundle $I$ of rank $m$ satisfying
\[
[I]=\sum_{i=1}^s(Q_i-Q_{i-1})\boxtimes S_i^\vee\quad\text{in }K^0(T\times T),\qquad c_m(I)=[\Delta_T]\quad\text{in }H^{2m}(T\times T,\R).
\]
\end{lemma}

\begin{proof}
Writing superscripts $(1)$ and $(2)$ for the two factors, let the fiber of $I$ consist of compatible tuples of linear maps $f_i\colon S_i^{(2)}\to Q_i^{(1)}$ ($1\le i\le s$) such that $f_i|_{S_{i+1}^{(2)}}$ equals $S_{i+1}^{(2)}\xrightarrow{f_{i+1}}Q_{i+1}^{(1)}\to Q_i^{(1)}$. Because quotient maps split locally, forgetting $f_s$ defines a vector-bundle surjection with kernel $\operatorname{Hom}\bigl(S_s^{(2)},\ker(Q_s^{(1)}\to Q_{s-1}^{(1)})\bigr)$, so induction on $s$ shows that $I$ is locally free of rank $\sum_{i=1}^s(r_i-r_{i-1})(n-r_i)=m$ with the stated $K$-class.

Restriction and projection define a surjection $\operatorname{End}(\C^n)\otimes\OO_{T\times T}\twoheadrightarrow I$, so $I$ is globally generated. The section induced by $\mathrm{id}_{\C^n}$ vanishes precisely where $S_i^{(2)}\subseteq S_i^{(1)}$ for all $i$, and hence set-theoretically along $\Delta_T$. At $(x,x)\in\Delta_T$, differentiating perturbations $1+\epsilon D$ of the first flag identifies the differential along the first factor (up to sign) with the isomorphism $\mathfrak{gl}_n/\mathfrak p_x\xrightarrow{\sim} I_{(x,x)}$, where $\mathfrak p_x$ is the parabolic stabilizer of $x$. Thus the section is regular and $c_m(I)=[\Delta_T]$, matching Pragacz's partial-flag diagonal construction \cite[\S1]{Pra08} (with factors exchanged and indices reversed; compare \cite[(7.6)]{Ful92} and \cite[\S2, (7)--(9)]{KaPr}).
\end{proof}

\begin{theorem}\label{thm:flag-specialization}
If $U$ is a nef vector bundle of rank $e>0$ on $T$, then $\varphi(U)$ is formally nef over $(V,\mathcal C)$. More precisely, it suffices that $F_\bullet$ satisfy the incidence condition of Definition~\ref{def:filtered} for the single auxiliary base $B=\PP_T(U)$ and the flag $W_i=p^*S_i^\vee$, where $p\colon B\to T$ is the projection.
\end{theorem}

\begin{proof}
Put $z=c_1(\OO_B(1))$ and $M=V\otimes H^*(B,\R)$. By Lemma~\ref{lem:partial-flag-diagonal}, the top Chern operator of the incidence bundle for $W_i=p^*S_i^\vee$ on $M$ is $\beta=(\varphi\otimes p^*)[\Delta_T]$. Since this rank-$m$ incidence bundle is formally nef by hypothesis, top Chern descent \cite[Theorem~8.10(a) and Example~2.9(b)]{LZ} equips $\beta M$ with a Hodge--Lefschetz structure of weight $d+\dim B-m=d+e-1$, polarized by $\mathcal C+\operatorname{Amp}(B)$, with $(\beta M)^{p,q}=\beta M^{p,q}$ and pairing $S_{\beta M}(\beta u,\beta v)=S_M(u,\beta v)$.

On the other hand, Grothendieck's projective bundle theorem \cite[\S2, axiom~A1, and \S3, (1)]{Gro} identifies the balanced tensor product
\[
R=V\otimes_{H^*(T,\R)} H^*(B,\R)=V[\xi]\big/\bigl(\xi^e-\varphi(c_1(U))\xi^{e-1}+\cdots+(-1)^e\varphi(c_e(U))\bigr)
\]
with $W(\varphi(U))=\bigoplus_{j=0}^{e-1}\xi^j V$ via the quotient map $\rho\colon M\to R$ sending $z$ to $\xi$. Writing $\sigma_j=p_*(z^{e-1+j})$ for $j\ge0$ and $\sigma_j=0$ for $j<0$, the projective bundle pairing $S_R(\xi^i x,\xi^j y)=(-1)^i S\bigl(x,\varphi(\sigma_{i+j-e+1})y\bigr)$ on $R$ is nondegenerate. Choosing dual homogeneous bases $(b_\nu)$ and $(b_\nu^\vee)$ of $H^*(T,\R)$ so that $[\Delta_T]=\sum_\nu b_\nu^\vee\otimes b_\nu$, we compute for $a\in H^{2k}(B,\R)$:
\begin{align*}
S_M\bigl(x\otimes a,\beta(y\otimes b)\bigr)
&=(-1)^k\sum_\nu S\bigl(x,\varphi(b_\nu^\vee)y\bigr)\int_B a\,p^*b_\nu\,b\\
&=(-1)^k S\bigl(x,\varphi(p_*(ab))y\bigr)
=S_R\bigl(\rho(x\otimes a),\rho(y\otimes b)\bigr),
\end{align*}
where the middle step uses the projection formula and the final step follows on monomials $a=z^i p^*c$, $b=z^j p^*c'$ because moving $\varphi(c)$ across $S$ contributes $(-1)^{\deg(c)/2}$. Nondegeneracy of $S_M$ and $S_R$ together with surjectivity of $\rho$ implies $\ker\beta=\ker\rho$, so $\beta u\mapsto\rho(u)$ is a bigrading-preserving isometry compatible with all descended operators.

Fix $H\in\operatorname{Amp}(T)$. Because $z$ is nef and $p$-ample, $tz+\epsilon p^*H\in\operatorname{Amp}(B)$ for all $t,\epsilon>0$, so the descended operators $h+t\xi+\epsilon\varphi(H)$ ($h\in\mathcal C$, $t,\epsilon>0$) polarize $R$ and lie in a single connected polarizing component $P\subseteq\operatorname{span}(\mathcal C,\xi,\varphi(H))$, which is an open convex cone by Kashiwara--Kawai \cite[Proposition~1.2.2]{KK} in the form of \cite[Proposition~3.11]{LZ}. Letting $h\to0$ along a ray and $\epsilon\to0$ places $\xi\in\overline P$. For any $h\in\mathcal C$, small projective perturbation \cite[Lemma~8.2]{LZ} shows that $h+\delta\xi$ polarizes $R$ for small $\delta>0$. Being a limit of points in $P$ inside the open polarizing locus, $h+\delta\xi$ belongs to $P$, and hence for every $t>\delta$ we have
\[
h+t\xi=(h+\delta\xi)+(t-\delta)\xi\in P+\overline P\subseteq P,
\]
proving that $\mathcal C+\R_{>0}\xi$ polarizes $R=W(\varphi(U))$.
\end{proof}

\subsection{Quotient Schubert bundles and character descent}
\begin{proof}[Proof of Theorem~\ref{thm:schubert-preservation}]
Every Schur functor of a globally generated bundle is globally generated, and on a complete flag variety $\mathbf S'_w$ is a quotient of a tensor product of exterior powers of the universal quotient bundles \cite[proof of Lemma~7.4]{PragaczKP} and is therefore globally generated as well. If $q\colon\widehat T\to T$ is a complete flag refinement and the Schubert bundle descends to $U$ on $T$, the identity $q_*\OO_{\widehat T}=\OO_T$ identifies $H^0(T,U)$ with $H^0(\widehat T,q^*U)$, so $U$ remains globally generated and hence nef on $T$. Theorem~\ref{thm:flag-specialization} then shows that $\varphi(U)$ is formally nef.

For a quotient flag $U_\ell\twoheadrightarrow\cdots\twoheadrightarrow U_1$ on $T$ with $U_\ell$ nef of positive rank, omit zero bundles and repeated ranks, reindex, and set $U_0=0$. Each $U_j$ is nef because $\PP_T(U_j)\hookrightarrow\PP_T(U_\ell)$ restricts $\OO_{\PP_T(U_\ell)}(1)$ to $\OO_{\PP_T(U_j)}(1)$ \cite[Theorem~6.2.12]{Laz}, so each $\varphi(U_j)$ is formally nef by Theorem~\ref{thm:flag-specialization}. Setting $D_j=\ker(U_j\to U_{j-1})$, the Whitney identity $c(\varphi(U_j);u)=c(\varphi(D_j);u)c(\varphi(U_{j-1});u)$ exhibits $\varphi(U_1),\ldots,\varphi(U_\ell)$ as the total quotients of the formal flag with graded pieces $\varphi(D_j)$, which is formally filtered nef by Theorem~\ref{thm:abstract-stages}.
\end{proof}

Kra\'skiewicz and Pragacz introduced Schubert functors in \cite{KP87}, defining the submodule version $\mathbf S_w$ in \cite[\S1]{KP} and identifying its torus character with $\mathfrak S_w$ in \cite[Theorem~4.1]{KP}. In the geometric setting, the quotient presentation above shows directly that $\mathbf S'_w$ preserves nefness, extending Fulton's ampleness observation \cite[\S4]{Ful95}.

\begin{remark}[Checking descent via characters]\label{rem:kp-descent}
On a complete refinement $q\colon\widehat T\to T$, each universal quotient is a quotient of $\C^n\otimes\OO_{\widehat T}$, so the $\operatorname{GL}_n$-equivariant bundle $\mathbf S'_w(\widehat Q_\bullet)$ of rank $R$ is a quotient of $E_0\otimes\OO_{\widehat T}$ for a fixed $\operatorname{GL}_n$-representation $E_0$ built from tensor products of exterior powers of $\C^n$. Whenever the torus character of $\mathbf S'_w$ in the successive line quotients of $\widehat T$ is symmetric within each block between consecutive marked ranks (including the final kernel block), the sum of all weights in that character extends to a character of the parabolic subgroup defining $T$, so $\det\mathbf S'_w(\widehat Q_\bullet)$ is the pullback of a line bundle on $T$ and is trivial along every fiber of $q$. Consequently, the classifying morphism $\widehat T\to\Gr_R(E_0)$ pulls the very ample Pl\"ucker line bundle back to a fiberwise trivial bundle and must be constant on the connected projective fibers of $q$. Because $q$ is Zariski-locally trivial, $\mathbf S'_w(\widehat Q_\bullet)$ descends canonically to $T$. Moreover, if $U$ and $U'$ are bundles on $T$, the projection formula and $q_*\OO_{\widehat T}=\OO_T$ yield $\operatorname{Hom}_T(U,U')\simeq\operatorname{Hom}_{\widehat T}(q^*U,q^*U')$, while surjectivity descends by faithful flatness of $q$, so equivariant quotient maps between descending Schubert bundles descend as well.
\end{remark}

\section{Schubert operators and their incidence modules}\label{sec:schubert-descent}

\subsection{Formal filtered nefness on the incidence module}
Let $(V,S)$ be a Hodge--Lefschetz module of weight $d$ polarized by $\mathcal C$, and let $F^{(1)}_\bullet,\ldots,F^{(a)}_\bullet$ be formally filtered nef flags in a common commuting operator algebra $A$ on $V$. Choosing for each $1\le i\le a$ a compatible Schubert polynomial $\mathfrak S_{w_i}$, set
\[
\Gamma=\prod_{i=1}^a\mathfrak S_{w_i}\bigl(F^{(i)}_\bullet\bigr),\qquad k_i=\ell(w_i),\qquad k=\sum_{i=1}^a k_i,\qquad m=d-k.
\]
Each factor acts as a commuting operator of bidegree $(k_i,k_i)$. Following \cite[Example~2.9(b)]{LZ}, we endow $\Gamma V$ with the bigrading $(\Gamma V)^{p,q}=\Gamma(V^{p,q})$ of weight $m$ and the bilinear form
\begin{equation}\label{eq:schubert-induced-form}
S_\Gamma(\Gamma x,\Gamma y)=S(x,\Gamma y),
\end{equation}
which is well defined and nondegenerate by the adjoint property of $\Gamma$ and the nondegeneracy of $S$.

\begin{proposition}\label{prop:schubert-incidence-module}
There exist a smooth complex projective variety $Z$, a product $\beta$ of top Chern operators on $V\otimes H^*(Z,\R)$, and a linear map
\[
j\colon V\longrightarrow W=\beta\bigl(V\otimes H^*(Z,\R)\bigr),\qquad j(x)=\beta(x\otimes1),
\]
such that $W$ is a Hodge--Lefschetz module of weight $m$ polarized by $\mathcal C+\operatorname{Amp}(Z)$ and
\[
S_W(jx,jy)=S(x,\Gamma y)\qquad(x,y\in V).
\]
Moreover, setting $U=j(V)$, the map $j(x)\mapsto\Gamma x$ induces a canonical isometry $U/(U\cap U^\perp)\simeq(\Gamma V,S_\Gamma)$ compatible with all operators in $\operatorname{span}(\mathcal C)$.
\end{proposition}

\begin{proof}
For each $1\le i\le a$, choose a partial flag variety $T_i=\operatorname{GL}_{N_i}/P_i$ of the marked quotient ranks of $F^{(i)}_\bullet$ with ambient dimension $N_i$ large enough that the complementary Chern relations (starting in degree $N_i-r+1$, where $r$ is the top marked rank) vanish on $V$ and $w_i\in S_{N_i}$. Put $t_i=\dim T_i$, fix a Borel subgroup $B_i\subseteq P_i$, and compose the Bott--Samelson resolution in $\operatorname{GL}_{N_i}/B_i$ of the minimal coset representative of $\sigma_{w_i}$ \cite[Theorem~3.5]{Hansen}, \cite[\S\S3.11--3.12]{Dem} with the projection $\operatorname{GL}_{N_i}/B_i\to T_i$ (which is an isomorphism on open Schubert cells) to obtain a projective resolution $b_i\colon Z_i\to T_i$ of the Schubert variety dual to $\sigma_{w_i}$, with $\dim Z_i=t_i-k_i$. Set $Z=\prod_{i=1}^a Z_i$.

Denoting by $S_{i,j}$ the kernel of the $j$th universal quotient on $T_i$, the pullbacks $b_i^*S_{i,j}^\vee$ form a quotient flag on $Z_i$ with globally generated first stage. Suppressing pullbacks to $Z$, form the formal incidence bundle
\[
K_i=\sum_j\bigl(F^{(i)}_j-F^{(i)}_{j-1}\bigr)\boxtimes b_i^*S_{i,j}^\vee,\qquad \rk K_i=t_i
\]
(taking the operator to be $1$ whenever $t_i=0$). By Lemma~\ref{lem:partial-flag-diagonal}, $c_{t_i}(K_i)$ equals the pullback of $[\Delta_{T_i}]$ specialized on the first factor at $F^{(i)}_\bullet$, so integrating over $Z_i$ evaluates the second factor against $(b_i)_*[Z_i]$ and yields $\mathfrak S_{w_i}(F^{(i)}_\bullet)$ by Schubert duality \cite[Proposition~5.1]{Pra96}. Setting $\beta=\prod_{i=1}^a c_{t_i}(K_i)$, integration over $Z$ satisfies the operator identity $(1\otimes\int_Z)\beta=\Gamma$.

By Proposition~\ref{prop:inc-module} (or \cite[Corollary~8.12]{LZ}), $W=\beta(V\otimes H^*(Z,\R))$ is polarized of weight $d+\dim Z-\sum_{i=1}^a t_i=m$ by $\mathcal C+\operatorname{Amp}(Z)$. Evaluating the tensor-product pairing on $x\otimes1$ and $\beta(y\otimes1)$ integrates out the $Z$-factor and yields $S_W(jx,jy)=S(x,\Gamma y)$. Finally, nondegeneracy of $S$ implies that the radical $U\cap U^\perp$ of $S_W|_U$ equals $j(\ker\Gamma)$, so $j(x)\mapsto\Gamma x$ induces the stated isometry.
\end{proof}

\begin{theorem}[Flags on Schubert incidence modules]\label{thm:schubert-incidence-preservation}
Let $W=\beta(V\otimes H^*(Z,\R))$ be the Schubert incidence module of Proposition~\ref{prop:schubert-incidence-module}, polarized by $\mathcal C_Z=\mathcal C+\operatorname{Amp}(Z)$. Every formally filtered nef flag $G_\bullet$ on $(V,\mathcal C)$ whose Chern operators commute with those of $F^{(1)}_\bullet,\ldots,F^{(a)}_\bullet$ induces a formally filtered nef flag on $(W,\mathcal C_Z)$.

In particular, when $\Gamma$ is a product of top Chern operators of commuting formally nef bundles on $V$, every formally filtered nef flag commuting with those bundles remains formally filtered nef on $(\Gamma V,\mathcal C)$ equipped with $S_\Gamma(\Gamma u,\Gamma v)=S(u,\Gamma v)$.
\end{theorem}

\begin{proof}
Set $V_Z=V\otimes H^*(Z,\R)$, write $\beta=\prod_i c_{t_i}(K_i)$, and let $U=G_j$ be any positive-rank total quotient of $G_\bullet$ with projective variable $\xi$. Realizing $U$ over $Z$ as the incidence bundle of the constant flag equal to $\OO_Z$ up to stage $j$ and $0$ thereafter, Proposition~\ref{prop:inc-module} shows that each $K_i$ is formally nef over the projective bundle module $W(U)$ on $V_Z$ equipped with $\mathcal C_Z+\R_{>0}\xi$, so simultaneous top Chern descent \cite[Corollary~8.12]{LZ} polarizes $\beta W(U)$ by $\mathcal C_Z+\R_{>0}\xi$.

Since $\beta$ commutes with the Chern operators of $U$, its image in the free basis $1,\xi,\ldots,\xi^{\rk U-1}$ decomposes as $\bigoplus_{c=0}^{\rk U-1}\xi^c W$, identifying $\beta W(U)$ with the projective bundle module of the induced operators on $W$. Indeed, for $u,v\in V_Z$ and $l=c+c'-\rk U+1$, the descended pairing on $\xi^c\beta u$ and $\xi^{c'}\beta v$ equals
\[
(-1)^c S_Z\bigl(u,\beta s_l(U)v\bigr)=(-1)^c S_W\bigl(\beta u,s_l(U)\beta v\bigr).
\]
Thus every positive-rank total quotient of $G_\bullet$ is formally nef on $W$, and Theorem~\ref{thm:abstract-stages} shows that the induced flag is formally filtered nef. Taking $Z$ to be a point with $K_i$ equal to the specified bundles on $V$ yields the final assertion.
\end{proof}

\subsection{Signature bounds on the base}\label{subsec:pairing-base}
Because $j\colon V\to W$ commutes with operators in $\overline{\mathcal C}$ and satisfies $S(x,\Gamma y)=S_W(jx,jy)$, restricting the mixed Hodge--Riemann form on $W$ bounds the negative index of the $\Gamma$-twisted pairing on $V^{p,q}$.

\begin{corollary}\label{cor:schubert-index-module}
Let $p+q\le m$ and $\eta_1,\ldots,\eta_{m-p-q}\in\overline{\mathcal C}$, and put $L=\prod_{j=1}^{m-p-q}\eta_j$ (with $L=1$ when $p+q=m$). Then the Hermitian form
\begin{equation}\label{eq:schubert-hermitian}
(x,y)\longmapsto \ii^{\,q-p}S(x,\Gamma L\overline y)\qquad(x,y\in V^{p,q})
\end{equation}
has at most $\dim W^{p-1,q-1}$ negative eigenvalues. Moreover, if $V^{p-2,q-2}=0$, this bound may be replaced by $\dim V^{p-1,q-1}$. In particular, if $V^{p-1,q-1}=0$, the form \eqref{eq:schubert-hermitian} is positive semidefinite.
\end{corollary}

\begin{proof}
Perturbing each $\eta_j\otimes1$ by small positive multiples of a base polarizer and an ample class on $Z$ yields polarizers on $W$, for which the mixed Hodge--Riemann form on $W^{p,q}$ is positive definite on a primitive subspace of codimension $\dim W^{p-1,q-1}$. Restricting along $j$ and passing to the limit preserves the upper bound $\dim W^{p-1,q-1}$ on the number of negative eigenvalues. Since each Bott--Samelson variety $Z_i$ is an iterated $\PP^1$-bundle over a point, its cohomology is concentrated on the diagonal and
\[
\dim W^{p-1,q-1}\le\sum_{b\ge0}\dim V^{p-1-b,q-1-b}\,\dim H^{b,b}(Z).
\]
Because any base polarizer acts injectively on $V$ below middle degree and $p+q\le m\le d$, the vanishing $V^{p-2,q-2}=0$ forces $V^{p-1-b,q-1-b}=0$ for all $b\ge1$, leaving only $\dim V^{p-1,q-1}\dim H^{0,0}(Z)=\dim V^{p-1,q-1}$ by connectedness of $Z$ (and $0$ when $V^{p-1,q-1}=0$).
\end{proof}

When $V$ is effective, Corollary~\ref{cor:schubert-index-module} implies positive semidefiniteness in bidegree $(0,0)$, corresponding on a connected Poincar\'e duality algebra to the nonnegativity of Schubert intersection numbers. In bidegree $(1,1)$ with $\dim V^{0,0}=1$, the negative of \eqref{eq:schubert-hermitian} is the real intersection form, which therefore has at most one positive eigenvalue.

By contrast, the quotient $\Gamma V\simeq U/(U\cap U^\perp)$ in Proposition~\ref{prop:schubert-incidence-module} need not satisfy the Hodge--Riemann relations in higher bidegrees. Indeed, on $X=(\PP^2)^3$ with factor hyperplane classes $x,y,z$ and $E=\OO(1,0,0)\oplus\OO(0,1,0)\oplus\OO(0,0,1)$, put $\Gamma=c_2(E)=xy+xz+yz$ and $h=x+y+z$. Inside $A=\R[x,y,z]/(x^3,y^3,z^3)$, the class
\[
v=2\Gamma-5(x^2+y^2+z^2)
\]
satisfies $\Gamma hv=0$ and $\int_X\Gamma v^2=-36$, defining a nonzero primitive $(2,2)$-class in $A/\Ann(\Gamma)$ whose primitive square has the opposite sign to that required by the Hodge--Riemann relations (compare Ross and Toma \cite[Example~9.2]{RT}, where an ample real twist of the same bundle yields a degenerate middle pairing).

\subsection{A Schubert positivity criterion for two blocks}\label{subsec:two-block}
Although successive kernels in a quotient flag can have negative first Chern class, for two marked ranks we can characterize exactly when multiplying a Schur class of the smaller quotient by a power of the relative first Chern class produces a Schubert-positive polynomial. Let $1\le s<r$, set $t=r-s$ and $L=x_{s+1}+\cdots+x_r$, and write $s_\lambda(x_1,\ldots,x_s)$ in this subsection for the symmetric Schur polynomial in the roots (padding a partition $\lambda$ of length at most $s$ with trailing zeros up to $\lambda_s$).

\begin{proposition}\label{prop:block-criterion}
For every integer $q\ge1$,
\begin{equation}\label{eq:relative-schubert-power}
s_\lambda(x_1,\ldots,x_s)L^q\text{ is Schubert-positive}\quad\Longleftrightarrow\quad \lambda_s\ge q+t-1.
\end{equation}
Consequently, for any formally filtered nef flag $F_1,F_2$ of quotient ranks $(s,r)$ on an effective Hodge--Lefschetz module $V$ with $\dim V^{0,0}=1$, the operator $s_\lambda(F_1)\bigl(c_1(F_2)-c_1(F_1)\bigr)^q$ has nonnegative intersection against complementary products of nef divisors whenever $\lambda_s\ge q+t-1$.
\end{proposition}

\begin{proof}
For sufficiency, expanding $L^q$ by Pieri's rule expresses it as a nonnegative combination of Schur polynomials $s_\mu(x_{s+1},\ldots,x_r)$ indexed by partitions $\mu\vdash q$ of length at most $t$. When $\lambda_s\ge q+t-1\ge\mu_1+t-1$, the concatenated exponent sequence
\[
(\lambda_1+s-1,\ldots,\lambda_s,\,\mu_1+t-1,\ldots,\mu_t)
\]
is weakly decreasing, so the corresponding monomial is a dominant Schubert polynomial. Applying the longest divided-difference operators separately in the first $s$ variables and in the next $t$ variables produces $s_\lambda(x_1,\ldots,x_s)s_\mu(x_{s+1},\ldots,x_r)$ by the bialternant formula. Because divided differences send Schubert polynomials to Schubert polynomials or zero \cite{LS82}, this nonzero product is a single Schubert polynomial with descents contained in $\{s,r\}$ (in particular, when $q=1$ and $\lambda_s\ge t$, the polynomial $s_\lambda(x_1,\ldots,x_s)L$ itself is a single Schubert polynomial).

Conversely, let $w$ be the Grassmannian permutation of $s_\lambda(x_1,\ldots,x_s)$, whose first $s$ values increase starting at $p+1$ where $p=\lambda_s$, and whose remaining values appear in increasing order. By Monk's formula \cite[Theorem~3]{Monk} in the polynomial form of \cite[(4.15)]{Macdonald}, multiplication by $L=\mathfrak S_{s_r}-\mathfrak S_{s_s}$ is a signed sum of Bruhat covers interchanging positions $i<j$: steps with $s<i\le r<j$ appear with sign $+1$ and preserve $\sum_{a=1}^s w(a)$, whereas steps with $i\le s<j\le r$ appear with sign $-1$ and strictly increase $\sum_{a=1}^s w(a)$.

Suppose now that $p<q+t-1$. Among values outside the first $s$ positions, the smallest $p$ values are $1,\ldots,p$. Denoting the $(p+1)$st value there by $v$, we have $v-1$ among the first $s$ positions. If $p<t$, then $v$ already lies in positions $s+1,\ldots,r$. If $p\ge t$, applying $p-t+1\le q-1$ positive covers raises the entry at position $r$ successively from $t$ through $p$ to $v$. Next interchange $v-1$ and $v$, a negative cover increasing the sum of the first $s$ values by exactly $1$, and complete the remaining steps to $q$ via positive covers exchanging the largest entry in positions $s+1,\ldots,r$ with the smallest larger entry to the right of $r$. Every $q$-step Monk path reaching this final permutation must contain a single negative step because $\sum_{a=1}^s w(a)$ has increased by $1$, so its Schubert coefficient in $s_\lambda L^q$ is strictly negative. Combining \eqref{eq:relative-schubert-power} with Corollary~\ref{cor:schubert-index-module} in bidegree $(0,0)$ yields the numerical statement. When $\lambda_s<q+t-1$, Fulton's test \cite[\S3]{Ful95} constructs an ample filtered bundle on which $s_\lambda(F_1)(c_1(F_2)-c_1(F_1))^q$ has negative degree.
\end{proof}

\section{Lorentzian polynomials and translated Schubert classes}\label{sec:ineq-lorentzian}

By Br\"and\'en and Huh \cite[\S2.4]{BH20}, a nonzero homogeneous polynomial of degree $m\ge2$ with nonnegative coefficients is Lorentzian if and only if its support is $M$-convex (meaning that for any exponent vectors $u,v$ in the support and any index $i$ with $u_i>v_i$, there exists $j$ with $u_j<v_j$ such that both $u-e_i+e_j$ and $v+e_i-e_j$ lie in the support) and every partial derivative of order $m-2$ has Hessian with at most one positive eigenvalue. In degrees $0$ and $1$, Lorentzianity reduces to coefficient nonnegativity. Throughout this section, $(V,S)$ denotes an effective Hodge--Lefschetz module of weight $d$ with $V^{0,0}=\C e$ for a chosen real generator $e$ (taking $e=1$ on cohomology algebras), we write $\int_V T=S(e,Te)$ for operators $T$ of bidegree $(d,d)$, and we call operators in $\overline{\mathcal C}$ nef.

\begin{lemma}\label{lem:ineq-volume}
Let $(U,S_U)$ be an effective polarized Hodge--Lefschetz module of weight $m\ge0$ with $\dim U^{0,0}\le1$, let $a\in U^{0,0}$ be real, and let $D_1,\ldots,D_s$ be commuting nef operators in the closure of a common polarizing cone on $U$. Then
\[
\frac1{m!}S_U\bigl(a,(t_1D_1+\cdots+t_sD_s)^m a\bigr)
\]
is Lorentzian or zero.
\end{lemma}

\begin{proof}
Assuming $a\ne0$ and taking $D_1,\ldots,D_s$ initially in the open polarizing cone, the mixed Hodge--Riemann relations in bidegree $(0,0)$ force every degree-$m$ monomial coefficient to be strictly positive (settling $m\in\{0,1\}$ and ensuring full simplex support). When $m\ge2$ and $P$ is any product of $m-2$ factors from $\{D_1,\ldots,D_s\}$, the Hessian $\bigl(S_U(a,D_iD_jPa)\bigr)_{1\le i,j\le s}$ is the pullback along $(b_i)\mapsto\sum_i b_i D_i a$ of the real bilinear form $-S_U(v,Pw)$ on $U^{1,1}_{\R}$, which is negative definite on the mixed primitive hyperplane and has at most one positive direction coming from $U^{0,0}$. Perturbing general nef operators as $D_i+\varepsilon H$ for a polarizer $H$ and letting $\varepsilon\downarrow0$ concludes the proof by closedness of the space of Lorentzian polynomials \cite[Theorem~2.25]{BH20}.
\end{proof}

\begin{theorem}\label{thm:ineq-schubert-volume}
Let $\mathcal F_1,\ldots,\mathcal F_r$ be formally filtered nef flags on $(V,\mathcal C)$ whose Chern operators commute with one another and with $\mathcal C$. For compatible permutations $w_1,\ldots,w_r$, set
\[
\Gamma=\prod_{i=1}^r\mathfrak S_{w_i}(\mathcal F_i),\qquad m=d-\sum_{i=1}^r\ell(w_i)\ge0.
\]
Then for all commuting nef divisors $D_1,\ldots,D_s\in\overline{\mathcal C}$, the polynomial
\[
\frac1{m!}\int_V\Gamma(t_1D_1+\cdots+t_sD_s)^m
\]
is Lorentzian or zero.
\end{theorem}

\begin{proof}
Let $j\colon V\to W$ be the Schubert incidence module of Proposition~\ref{prop:schubert-incidence-module}, and put $a=j(e)\in W^{0,0}$. Because $V$ is effective and $Z$ is connected, $W$ is effective with $\dim W^{0,0}\le1$ spanned by $a$, the operators $D_i\otimes1$ lie in $\overline{\mathcal C_Z}$, and $S_W(a,D_{i_1}\cdots D_{i_m}a)=\int_V\Gamma D_{i_1}\cdots D_{i_m}$. Applying Lemma~\ref{lem:ineq-volume} to $W$ yields the result.
\end{proof}

When $m\ge2$ and $M$ is a product of $m-2$ nef divisors, differentiating once in each factor of $M$ and setting its variable to zero leaves a binary quadratic form in $H,K\in\overline{\mathcal C}$ with nonpositive determinant, yielding the Hodge-index inequality
\[
\left(\int_V\Gamma MHK\right)^2\ge\left(\int_V\Gamma MH^2\right)\left(\int_V\Gamma MK^2\right)
\]
for any fixed product $\Gamma$ (while nonnegative linear combinations of such products still retain nonnegative coefficients).

To compare Schubert classes of varying degrees, define the derived polynomials $g^{[j]}$ of a polynomial $g(x_1,\ldots,x_e)$ in the roots of a flag via the uniform root shift
\[
g(x_1+u,\ldots,x_e+u)=\sum_{j\ge0}g^{[j]}(x_1,\ldots,x_e)u^j,
\]
so $g^{[j]}$ is the $j$th divided derivative along $(1,\ldots,1)$ and vanishes for $j>\deg g$.

\begin{theorem}[Lorentzianity of derived Schubert classes]\label{thm:ineq-luzheng}
With $(V,\mathcal C)$ and $\mathcal F_1,\ldots,\mathcal F_r$ as in Theorem~\ref{thm:ineq-schubert-volume}, put $g_i=\mathfrak S_{w_i}$ and $k_i=\ell(w_i)$. For any nonnegative integers $m,n_1,\ldots,n_r$ satisfying $d=m+\sum_{i=1}^r(k_i-n_i)$, the polynomial
\begin{equation}\label{eq:ineq-luzheng}
\sum_{\substack{a_1+\cdots+a_r=m\\0\le a_i\le n_i}} \frac{x_1^{a_1}\cdots x_r^{a_r}}{a_1!\cdots a_r!}\int_V\prod_{i=1}^r g_i^{[n_i-a_i]}(\mathcal F_i)
\end{equation}
is Lorentzian or zero.
\end{theorem}

\begin{proof}
Set $Y=\prod_{i=1}^r\PP^{n_i}$ with factor hyperplane classes $h_i$ normalized by $\int_Y\prod_i h_i^{n_i}=1$ (taking $h_i=0$ when $n_i=0$), and twist $\mathcal F_i$ over $V\otimes H^*(Y,\R)$ to $\mathcal F_i\langle h_i\rangle$. In any auxiliary test with incidence bundle $K_i$ over a base $B$, each root of the incidence class involves a single root from $\mathcal F_i$, so twisting replaces $K_i$ by $K_i\langle h_i\rangle$ and shifts its projective relation from $P_i(z)$ to $P_i(z-h_i)$, identifying the new hyperplane operator with $\xi_i+h_i$. For $N\in\mathcal C$, $b\in\operatorname{Amp}(B)$, $c\in\operatorname{Amp}(Y)$, and $t>0$, rewriting $N+b+c+t(\xi_i+h_i)=(N+b+t\xi_i)+(c+th_i)$ expresses it as the sum of a polarizer on $W(K_i)$ and an ample class on $Y$, which polarizes $W(K_i\langle h_i\rangle)$ by \cite[Proposition~2.11]{LZ}.

Consequently, Theorem~\ref{thm:ineq-schubert-volume} applies over $V\otimes H^*(Y,\R)$ to the twisted flags $\mathcal F_i\langle h_i\rangle$ and the nef classes $h_1,\ldots,h_r$, in residual weight $d+\sum_i n_i-\sum_i k_i=m$, showing that
\[
\frac1{m!}\int_{V\otimes H^*(Y)}\prod_{i=1}^r g_i\bigl(\mathcal F_i\langle h_i\rangle\bigr)\left(\sum_{i=1}^r x_i h_i\right)^m
\]
is Lorentzian or zero. Expanding the multinomial and integrating against $h_1^{a_1}\cdots h_r^{a_r}$ over $Y$ extracts the $u^{n_i-a_i}$ coefficient $g_i^{[n_i-a_i]}(\mathcal F_i)$ in each factor, recovering \eqref{eq:ineq-luzheng}.
\end{proof}

To upgrade Lorentzianity to strict Lorentzianity, we rule out Hessian kernels in all $(m-2)$nd derivatives by perturbing the flags negatively and differentiating with respect to the translation parameters.

\begin{theorem}\label{thm:strict-filtered}
In Theorem~\ref{thm:ineq-luzheng}, assume that $r\ge1$, $m\le n_i\le k_i$, and $g_i(1,\ldots,1)>0$ for every $1\le i\le r$. If there exist $H\in\mathcal C$ and $\epsilon>0$ such that each negatively twisted flag $\mathcal G_i=\mathcal F_i\langle-\epsilon H\rangle$ is formally filtered nef, then \eqref{eq:ineq-luzheng} is strictly Lorentzian.
\end{theorem}

\begin{proof}
For $t=(t_1,\ldots,t_r)\in\R_{\ge0}^r$, set $\mathcal F_i(t)=\mathcal G_i\langle t_i H\rangle$. Convexity of $\mathcal C$ ensures that shifting the hyperplane operator to $\xi+t_i H$ preserves formal filtered nefness. Adjoining an extra variable $y$ (indexed by $0$) for $H$, define for $a=(a_0,a_1,\ldots,a_r)$ of total weight $|a|=d$ the coefficients
\[
c_a(t)=\int_V H^{a_0}\prod_{i=1}^r g_i^{[k_i-a_i]}\bigl(\mathcal F_i(t)\bigr)
\]
(setting $c_a(t)=0$ if $a_i>k_i$ for some $i\ge1$). Applying the construction of Theorem~\ref{thm:ineq-luzheng} on $Y=\prod_{i=1}^r\PP^{k_i}$ with divisor $yH+\sum_{i=1}^r x_i h_i$ shows that $P_t=\sum_{|a|=d}c_a(t)y^{a_0}\prod_{i=1}^r x_i^{a_i}/(a_0!\prod_{i=1}^r a_i!)$ is Lorentzian or zero. In particular, $c_b(0)\ge0$ for all $b$. Because $c_{(d,0,\ldots,0)}(0)=\int_V H^d\prod_{i=1}^r g_i(1,\ldots,1)>0$ and expanding $c_a(t)$ about $t=0$ expresses it as a nonnegative combination of the $c_b(0)$ in which $c_{(d,0,\ldots,0)}(0)$ appears with coefficient $\prod_{i=1}^r\binom{k_i}{a_i}t_i^{a_i}$, we have $c_a(t)>0$ whenever $0\le a_i\le k_i$ and all $t_i>0$. Evaluating at $t_i=\epsilon$, the polynomial \eqref{eq:ineq-luzheng} equals
\[
\left.\partial_{x_1}^{k_1-n_1}\cdots\partial_{x_r}^{k_r-n_r}P_t\right|_{y=0},
\]
whose degree-$m$ coefficients are all strictly positive, settling $m\in\{0,1\}$.

Assume now that $m\ge2$, fix $t\in\R_{>0}^r$, and let $\alpha=(\alpha_0,\alpha_1,\ldots,\alpha_r)$ satisfy $|\alpha|=d-2$ and $\alpha_i\le k_i-2$ for $1\le i\le r$. The Hessian $A_\alpha=(c_{\alpha+e_j+e_h}(t))_{0\le j,h\le r}$ of $\partial^\alpha P_t$ has strictly positive entries and at most one positive eigenvalue. We claim that $\ker A_\alpha=0$. For any $u\in\ker A_\alpha$ and $1\le i\le r$, the function $u^T A_\alpha(s)u-\bigl(u^T A_\alpha(s)e_0\bigr)^2/(A_\alpha(s))_{00}$ is nonpositive near $s=t$ (since $e_0^\perp$ is negative semidefinite for $A_\alpha(s)$) and vanishes at $s=t$, so differentiating at $s=t$ forces $u^T(\partial_{t_i}A_\alpha)u=0$.

Meanwhile, differentiating $c_a(t)$ directly yields
\[
\partial_{t_i}c_a=\begin{cases}
(k_i-a_i+1)c_{a-e_i+e_0},&a_i>0,\\
0,&a_i=0.
\end{cases}
\]
Writing $A=A_\alpha$, suppose first that $\alpha_i>0$ and set $B=A_{\alpha-e_i+e_0}$. Since $Be_i=Ae_0$, entrywise differentiation yields the matrix identity
\[
\partial_{t_i}A=(k_i-\alpha_i+1)B-e_i(e_0^T A)-(Ae_0)e_i^T,
\]
whence $u^TBu=0$ and $u^TBe_i=0$. Because $B_{ii}=A_{i0}>0$ and $B$ has at most one positive eigenvalue, $B$ is negative semidefinite on $(Be_i)^\perp$, forcing $Bu=0$. Iterating this reduction lowers $\alpha_i$ to $0$ while keeping $u$ in the kernel. At $\alpha_i=0$, entrywise differentiation reads
\[
\partial_{t_i}A=k_i\bigl(e_i(e_0^T A)+(Ae_0)e_i^T\bigr)-(k_i+1)A_{i0}\,e_i e_i^T,
\]
so evaluating on $(u,u)$ with $Au=0$ and $A_{i0}>0$ forces $u_i=0$. Carrying out this argument independently for each $1\le i\le r$ yields $u_1=\cdots=u_r=0$, and the zeroth row of $A_\alpha u=0$ then forces $u_0=0$.

Finally, any $(m-2)$nd partial derivative $\partial^\gamma$ of \eqref{eq:ineq-luzheng} in $x_1,\ldots,x_r$ has Hessian equal to the lower-right $r\times r$ block of $A_\alpha$ for $\alpha=(0,k_1-n_1+\gamma_1,\ldots,k_r-n_r+\gamma_r)$, which satisfies $\alpha_i\le k_i-2$. Since $r\ge1$ and the diagonal entry at $e_1$ is positive while $A_\alpha$ is nonsingular of signature $(1,r)$, this $r\times r$ block is nonsingular of signature $(1,r-1)$, so \eqref{eq:ineq-luzheng} is strictly Lorentzian by \cite[Definition~2.1]{BH20}.
\end{proof}

We call a positive-rank formal bundle $F$ \emph{formally ample} over $(V,\mathcal C)$ if its hyperplane operator $\xi_F$ lies in the open polarizing component of $W(F)$ inside $\operatorname{span}(\mathcal C,\xi_F)$ containing $N+t\xi_F$ for $N\in\mathcal C$ and small $t>0$.

\begin{corollary}\label{cor:stagewise-strict}
In Theorem~\ref{thm:ineq-luzheng}, if $r\ge1$, $m\le n_i\le k_i$, $g_i(1,\ldots,1)>0$, and every positive-rank total quotient of each $\mathcal F_i$ is formally ample, then \eqref{eq:ineq-luzheng} is strictly Lorentzian.
\end{corollary}

\begin{proof}
Fix $H\in\mathcal C$. By openness of the polarizing cone, there is a common $\epsilon>0$ such that $\xi_F-\epsilon H$ remains in the distinguished polarizing component for every positive-rank total quotient $F$ of $\mathcal F_1,\ldots,\mathcal F_r$. Since $\mathcal C$ lies in the closure of that component by \cite[Lemma~8.2]{LZ}, convexity ensures that $N+t(\xi_F-\epsilon H)$ polarizes $W(F)$ for all $N\in\mathcal C$ and $t>0$, meaning that $F\langle-\epsilon H\rangle$ is formally nef (conversely, formal nefness of $F\langle-\epsilon H\rangle$ places both $\xi_F=\epsilon H+(\xi_F-\epsilon H)$ and $H+c\xi_F=(1+c\epsilon)H+c(\xi_F-\epsilon H)$ in the same connected polarizing cone, so the two descriptions of formal ampleness are equivalent). By Theorem~\ref{thm:abstract-stages}, each $\mathcal F_i\langle-\epsilon H\rangle$ is formally filtered nef, so Theorem~\ref{thm:strict-filtered} applies.
\end{proof}

For one-stage flags on a smooth complex projective variety, Theorems~\ref{thm:ineq-luzheng} and~\ref{thm:strict-filtered} recover the projective case of Lu and Zheng's theorem \cite[Theorem~1.9]{LuZhengHR}, where strictness was proved via mixed Hodge--Riemann nondegeneracy \cite[Theorem~6.43 and Lemma~6.44]{LuZhengHR} alongside weak Hodge--Riemann relations for derived Schubert classes under $H^{p-2,q-2}(X)=0$ \cite[Corollary~5.27]{LuZhengHR}. When $r=2$, setting $b_j=\int_V g_1^{[n_1-j]}(\mathcal F_1)g_2^{[n_2-m+j]}(\mathcal F_2)$ ($0\le j\le m$) identifies $b_j/(j!(m-j)!)$ with the coefficient of $x_1^j x_2^{m-j}$, yielding log-concavity $b_j^2\ge b_{j-1}b_{j+1}$ ($0<j<m$) with interval support \cite[Example~2.26]{BH20} (and strict log-concavity under Theorem~\ref{thm:strict-filtered}). For Schur polynomials of two nef vector bundles $E,F$ on a smooth projective $d$-fold with $n_1=n_2=m=k_1+k_2-d\ge0$, this recovers Ross and Toma's log-concavity theorem \cite[Theorem~10.5 and Corollary~10.7]{RTnef}.

\section{The Hodge--Lefschetz package for flag matroids}\label{sec:matroids}

Let $X_E$ be the permutohedral variety on a finite nonempty ground set $E$. Associated with any matroid $M$ on $E$ is the dual tautological class $F_M=\mathcal S_M^\vee\in K^0(X_E)$ \cite[Definition~3.9]{BEST}. When $\rk M>0$, Larson and Partida proved that its formal projective bundle ring satisfies hard Lefschetz and the Hodge--Riemann relations for $h+t\xi$ whenever $h\in\operatorname{Amp}(X_E)$ and $t>0$ \cite[Example~1.1 and Theorem~1.4]{LP}.

Recall that $N\preceq M$ denotes a matroid quotient, meaning that every circuit of $M$ is a union of circuits of $N$ \cite[\S1.7.3]{BGW}. By \cite[Theorem~1.7.1]{BGW} (see also \cite[\S3.3]{BGW00}), chains $M_1\preceq\cdots\preceq M_s$ of matroid quotients on $E$ of strictly increasing ranks correspond bijectively to flag matroids on $E$, whose flags are nested basis tuples $B_1\subsetneq\cdots\subsetneq B_s$ ($B_i\in\mathcal B(M_i)$). Repeated constituents contribute zero in $K^0(X_E)$ and may be omitted.

\begin{theorem}\label{thm:mat}
Let $M_1\preceq\cdots\preceq M_s$ be the constituents of a flag matroid on $E$, and set $F_0=0$ and $F_i=\mathcal S_{M_i}^\vee\in K^0(X_E)$ for $1\le i\le s$. For every morphism $g\colon Y\to X_E$ from a smooth complex projective variety $Y$ and every flag of quotient vector bundles $W_1\twoheadrightarrow\cdots\twoheadrightarrow W_s$ on $Y$ with $W_1$ globally generated, the $K$-theory class
\[
K=\sum_{i=1}^s g^*(F_i-F_{i-1})\otimes W_i
\]
defines a formally nef bundle on $Y$ whenever it has positive rank.
\end{theorem}

\subsection{Line-bundle representatives for successive differences}
To invoke Theorem~\ref{thm:filtered}, we represent each $K$-theoretic difference $F_i-F_{i-1}$ on $X_E$ by an explicit direct sum of toric line bundles. For any linear order $\sigma$ on $E$, let $B_\sigma(M)$ be the greedy basis of $M$ formed by scanning $E$ from largest to smallest element and retaining each element that increases the rank of the prefix. Writing $t_e$ for the $e$th coordinate character of $(\C^*)^E$, the localization of $F_M=\mathcal S_M^\vee$ at the torus fixed point indexed by $\sigma$ has character $\sum_{e\in B_\sigma(M)}t_e$ \cite[Proposition~3.8 and Definition~3.9]{BEST}. When $N\preceq M$, the closure inclusion $\operatorname{cl}_M(A)\subseteq\operatorname{cl}_N(A)$ on every prefix $A$ implies $B_\sigma(N)\subseteq B_\sigma(M)$ \cite[Corollary~1.7.2]{BGW}.

Listing the elements of $B_\sigma(M)\setminus B_\sigma(N)$ in scan order, fix an index $1\le k\le\rk M-\rk N$. Across any wall between adjacent braid chambers where consecutive elements $e,f\in E$ are swapped, all greedy selections outside $\{e,f\}$ remain unchanged and the cardinality selected from $\{e,f\}$ stays fixed because prefix ranks before and after $\{e,f\}$ agree. Hence the $k$th character either stays constant or switches between $t_e$ and $t_f$, which agree along the wall $t_e=t_f$. Dividing by $t_b$ for a chosen base element $b\in E$ produces wall-compatible characters of the effective torus $(\C^*)^E/\C^*$, defining a toric line bundle on $X_E$ \cite[Theorem~4.2.12]{CLS} whose $(\C^*)^E$-linearization after twisting back by $t_b$ has fixed-point weight equal to the $k$th character in $B_\sigma(M)\setminus B_\sigma(N)$.

Summing these line bundles over $1\le k\le\rk M_i-\rk M_{i-1}$ produces a vector bundle $G_i$ on $X_E$ sharing all fixed-point characters with $F_i-F_{i-1}$. By injectivity of fixed-point restriction in equivariant $K$-theory \cite[Theorem~6.2]{VV}, \cite[Theorem~2.1(a)]{BEST}, we have $[G_i]=F_i-F_{i-1}$ and $[\bigoplus_{a\le i}G_a]=F_i$ in $K^0(X_E)$. Since each positive-rank $F_i$ is formally nef on $X_E$ \cite[Example~1.1 and Theorem~1.4]{LP}, Theorem~\ref{thm:filtered} establishes Theorem~\ref{thm:mat} when $g=\mathrm{id}_{X_E}$.

\subsection{Pullbacks of tautological classes}
To pass from $X_E$ to an arbitrary smooth projective variety $Y$, we verify that formal nefness on all torus orbit closures survives arbitrary pullbacks.

\begin{lemma}\label{lem:toric}
Let $X$ be a product of permutohedral varieties, and let $D$ be a vector bundle of positive rank on $X$ whose restriction to every torus orbit closure (including $X$ itself) is formally nef. Then $g^*D$ is formally nef for every morphism $g\colon Y\to X$ from a smooth complex projective variety $Y$.
\end{lemma}

\begin{proof}
Working on each connected component of $Y$, we induct on $d=\dim X$, the case $d=0$ being trivial. Because every toric boundary divisor of $X$ is itself a product of permutohedral varieties whose torus orbits are orbits of $X$, the induction hypothesis applies whenever $g(Y)$ lies in the toric boundary. We may therefore assume that $U=g^{-1}(T)$ is dense in $Y$, where $T$ is the open torus of $X$.

On $U\times T$, the fiber of $(y,z)\mapsto z/g(y)$ over $c\in T$ is the graph of $c\cdot g|_U$. Fixing a reference element $b$ in each permutohedral factor, the coordinate ratios $z_i/z_b$ ($i\ne b$) define torus coordinates $\chi_1,\ldots,\chi_d$ extending to morphisms $X\to\PP^1$ (since coordinate differences have constant weak sign on every braid cone), so $\beta_i(y,z)=\chi_i(z)/\chi_i(g(y))$ defines morphisms $Y\times T\to\PP^1$. Representing each pencil $\beta_i$ by two sections of a line bundle on $Y\times X$ without common divisorial components and successively principalizing their base ideals \cite[Theorem~2.5(i),(ii)]{Res} via blowups along smooth centers disjoint from $Y\times T$, we obtain a smooth projective blowup $W$ of $Y\times X$ equipped with projections $p_Y\colon W\to Y$, $p_X\colon W\to X$, and a morphism $\beta\colon W\to(\PP^1)^d$.

Pulling $D$ back to $Y\times X$ preserves formal nefness by Lemma~\ref{lem:grass}. At each smooth blowup step, every component of the blowup center maps into a boundary divisor of $X$, so $D$ restricts to a formally nef bundle on the center by the induction hypothesis. Since projectivization commutes with blowups and the blowup decomposition theorem has only the ambient projectivization and the projectivization over the center as constant-coefficient supports, Lemma~\ref{lem:support} preserves formal nefness at each step, showing that $p_X^*D$ is formally nef on $W$.

By generic smoothness, for general $c\in T$ the fiber $L=\beta^{-1}(c)\subseteq W$ is smooth of pure dimension $\dim Y$, cut out regularly by the globally generated bundle $\bigoplus_{i=1}^d\beta_i^*\OO_{\PP^1}(1)$. Because $W\setminus(U\times T)$ has dimension at most $\dim Y+d-1$, a general fiber meets the boundary in dimension at most $\dim Y-1$. Hence the open graph $\{z=c\cdot g(y)\}\subseteq U\times T$ is dense in $L$, the projection $p_Y|_L\colon L\to Y$ is birational, and $p_X|_L=(c\cdot g)\circ(p_Y|_L)$ holds on all of $L$. Since torus translation by $c\in T$ preserves Chern classes, $p_X^*D|_L$ and $(p_Y|_L)^*g^*D$ share the same rank and Chern classes, so $g^*D$ is formally nef by Lemma~\ref{lem:section}.
\end{proof}

\begin{proof}[Completion of the proof of Theorem~\ref{thm:mat}]
On any torus orbit closure of $X_E$, indexed by an ordered partition with prefix unions $\varnothing=A_0\subsetneq A_1\subsetneq\cdots\subsetneq A_k=E$, the restriction of $F_i$ splits in $K$-theory into the direct sum of the dual tautological classes of the minors $(M_i|A_j)/A_{j-1}$ pulled back from the permutohedral factors \cite[Proposition~5.3]{BEST}. Each positive-rank summand is formally nef by \cite[Theorem~1.4]{LP} and Lemma~\ref{lem:grass}, and their direct sum is formally nef by \cite[Proposition~8.14]{LZ}. Lemma~\ref{lem:toric} therefore shows that $g^*E_i$ is formally nef on $Y$ for every positive-rank partial sum $E_i=\bigoplus_{a\le i}G_a$, and Theorem~\ref{thm:filtered} applied to $g^*G_1,\ldots,g^*G_s$ yields formal nefness of $K$.
\end{proof}

\subsection{An illustrative two-stage example}
Consider the uniform flag $U_{1,2}\preceq U_{2,2}$ on a two-element set $E$. On $X_E=\PP^1$, we have $F_1=[\OO(1)]$ and $F_2=[\OO^{\oplus2}]$, so $F_2-F_1=[\OO(-1)]$. Pairing against the quotient flag $\OO^{\oplus2}\twoheadrightarrow\OO(1)$ on a second copy of $\PP^1$ yields the bundle representative
\[
K=\OO(1,0)^{\oplus2}\oplus\OO(-1,1)
\]
on $\PP^1\times\PP^1$. Writing $x,y$ for the two hyperplane classes, $c_1(K)=x+y$ and $c_2(K)=2xy$, and Theorem~\ref{thm:mat} guarantees that
\[
\R[x,y,\xi]\big/\bigl(x^2,\,y^2,\,\xi^3-(x+y)\xi^2+2xy\xi\bigr)
\]
satisfies mixed hard Lefschetz and the mixed Hodge--Riemann relations for all $ax+by+t\xi$ with $a,b,t>0$, even though $K$ admits a quotient of degree $-1$ along fibers of the second projection and therefore fails geometric nefness.

\subsection{Flag matroids on general fans}\label{subsec:fan-flags}
Whereas an ordering of $E$ selects greedy bases across cones of the braid fan, on a general rational fan $\Sigma$ we govern all constituents simultaneously via a single integral linear map $f$ from the cocharacter lattice of $\Sigma$ to $\mathbb Z^E$. Denote by $\mathrm{PP}(\Sigma)$ the graded real algebra of conewise polynomials on $\Sigma$ (with linear functions in degree one) \cite[\S2.1]{LP}. Ordinary Chow classes are obtained modulo global linear functions (restricting first to a complement of the lineality space when $\Sigma$ is not pointed).

Let $M_1\preceq\cdots\preceq M_s$ be a flag matroid on $E$ of positive ranks $r_1<\cdots<r_s$, and assume that $f$ maps every cone of $\Sigma$ into a cone of the normal fan of $P(M_i)$ for each $i$ (equivalently, into a cone of the normal fan of $\sum_{i=1}^s P(M_i)$). Writing $\mathfrak c_j(M_i)$ on the normal cone of a basis $B\in\mathcal B(M_i)$ for the $j$th elementary symmetric polynomial in $\{t_e : e\in B\}$ \cite[\S1]{LP}, define for the dual classes
\[
C_i(u)=\sum_{j=0}^{r_i}(-1)^j f^*\mathfrak c_j(M_i)u^j\qquad(1\le i\le s),\qquad C_0(u)=1.
\]

\begin{theorem}\label{thm:fan-flags}
The polynomials $C_1(u),\ldots,C_s(u)$ are the total Chern polynomials of a formally filtered nef flag in both of the following settings:
\begin{enumerate}
\item On the Chow ring $A^*(\Sigma)$ of a simplicial Lefschetz fan $\Sigma$ possessing a nonempty cone of strictly convex classes.
\item On the intersection cohomology module $\IH^*(\Sigma)$ of any rational projective fan $\Sigma$, polarized by its cone of strictly convex classes, where the coefficients of $C_i(u)$ act via conewise polynomials (matching {\rm(a)} when $\Sigma$ is also simplicial).
\end{enumerate}
\end{theorem}

\begin{proof}
To check that each ratio $C_i(u)/C_{i-1}(u)\in\mathrm{PP}(\Sigma)[[u]]$ is a polynomial of degree at most $r_i-r_{i-1}$, subdivide $\Sigma$ along the preimages of the braid chambers. On any refined cone, refining the weight preorder by a single linear order on $E$ breaks ties uniformly across all stages and yields nested greedy bases $B_1\subseteq\cdots\subseteq B_s$ \cite[Corollary~1.7.2]{BGW}, so
\[
\frac{C_i(u)}{C_{i-1}(u)}=\prod_{e\in B_i\setminus B_{i-1}}(1-f^*t_e\,u)
\]
vanishes in degrees exceeding $r_i-r_{i-1}$ on every refined cone and hence in $\mathrm{PP}(\Sigma)$.

In case (a), the $i$th projective bundle module has defining relation
\[
\xi^{r_i}+f^*\mathfrak c_1(M_i)\xi^{r_i-1}+\cdots+f^*\mathfrak c_{r_i}(M_i)=0
\]
and is polarized by $h+t\xi$ for every strictly convex class $h$ and $t>0$ by Larson and Partida \cite[Theorem~3.9]{LP}. Theorem~\ref{thm:abstract-stages} then proves formal filtered nefness.

In case (b), equip $V^j=\IH^{2j}(X_\Sigma,\R)$ with bidegree $(j,j)$. Realizing $V$ as global sections of the minimal extension sheaf on $\Sigma$ modulo global linear functions \cite{BBFK,BL03} equips $V$ with a natural action of $\mathrm{PP}(\Sigma)$ balanced with respect to the intersection pairing $\langle\,\cdot\,,\,\cdot\,\rangle_\Sigma$ \cite[Theorems~4.2 and~5.2]{BBFK05} (for a top-degree Schubert operator $T$, its degree is $\langle 1,T1\rangle_\Sigma$ where $1\in\IH^0(X_\Sigma,\R)$), and setting $S(u,v)=(-1)^j\langle u,v\rangle_\Sigma$ for $u\in V^j$ makes every strictly convex class $h$ polarize $(V,S)$ by Karu's theorem \cite[Theorems~0.1 and~0.2]{Karu}.

Choosing a projective simplicial subdivision $\pi\colon\Sigma'\to\Sigma$, the decomposition of minimal extension sheaves \cite[Theorem~2.6 and Corollary~2.7]{BL03} splits $\IH^*(\Sigma)$ off as a graded $\mathrm{PP}(\Sigma)$-linear direct summand of $A^*(\Sigma')$, and this splitting extends coefficientwise in $1,\xi,\ldots,\xi^{r_i-1}$ to the $i$th projective bundle modules with invariance under $\pi^*h$ and $\xi$. Over $A^*(\Sigma')$, adding positive multiples of the nef class $\pi^*h$ and $\xi$ to $a+\xi$ ($a\in\operatorname{Amp}(\Sigma')$) remains inside the polarizing cone of \cite[Theorem~1.4]{LP}, so constancy of the monodromy filtration \cite[Theorem~3.3(i)]{CK}, \cite[Proposition~2.14]{LZ} shows that the monodromy filtration of $\pi^*h+t\xi$---and hence of its restriction $h+t\xi$ over $\IH^*(\Sigma)$---is independent of $t>0$. Since $h+t\xi$ polarizes the projective bundle module over $\IH^*(\Sigma)$ for small $t>0$ by \cite[Lemma~8.2]{LZ}, hard Lefschetz and the positive definiteness of primitive pairings persist for all $t>0$, so each total quotient is formally nef and Theorem~\ref{thm:abstract-stages} applies.
\end{proof}

\begin{remark}[Integral polymatroids]\label{rem:polymatroids}
For a chain $P_1,\ldots,P_s$ of integral polymatroid quotients of strictly increasing positive ranks on $E$, replacing each $e\in E$ by a block of at least $\max_i\rk_{P_i}(\{e\})$ parallel copies produces a ground set $\widetilde E$ on which the multisymmetric matroid lifts $M_i$ \cite[Definition~3.1]{EL24} form a flag matroid \cite[Proposition~4.10]{BLS}. Pulling back along the diagonal weight repetition map $\delta\colon X_E\to X_{\widetilde E}$ via Theorem~\ref{thm:fan-flags} shows that $\delta^*(\mathcal S_{M_1}^\vee),\ldots,\delta^*(\mathcal S_{M_s}^\vee)$ form the total quotients of a formally filtered nef flag on $X_E$.
\end{remark}

\subsection{Ordinary and augmented matroid Chow rings}
Specializing Theorem~\ref{thm:fan-flags} to the ordinary Bergman fan (whose Chow ring was introduced by Feichtner--Yuzvinsky \cite[Definition~3]{FY} and studied by Adiprasito--Huh--Katz \cite[\S1]{AHK}) and to the augmented Bergman fan of Braden, Huh, Matherne, Proudfoot, and Wang \cite[\S1.2]{BHMPW22} yields filtered packages over any loopless base matroid $N$.

\begin{corollary}\label{cor:both-chow-rings}
Let $N$ be an arbitrary loopless matroid on $E$, with ordinary and augmented Chow rings $A_{\rm ord}(N)$ and $A_{\rm aug}(N)$. For any flag matroid $M_1\preceq\cdots\preceq M_s$ on $E$ of positive ranks $r_1<\cdots<r_s$ (chosen independently of $N$), the restricted classes
\[
\begin{array}{c|c}
A_{\rm ord}(N) & \mathcal S_{M_1}^{\vee},\ldots,\mathcal S_{M_s}^{\vee}\\[2pt]
A_{\rm aug}(N) & \mathcal Q^{\rm aug}_{M_1^*},\ldots,\mathcal Q^{\rm aug}_{M_s^*}
\end{array}
\]
form the total quotients of formally filtered nef flags of ranks $r_1<\cdots<r_s$, where $M^*$ denotes the dual matroid of $M$.
\end{corollary}

\begin{proof}
Both Bergman fans satisfy Poincar\'e duality, hard Lefschetz, and the Hodge--Riemann relations for strictly convex classes \cite[Theorem~1.6 and Remark~1.7]{BHMPW22}, with top Minkowski weights generated by the constant weight $1$ \cite[Proposition~2.8 and Remark~2.14]{BHMPW22} and ample cones nonempty by restriction from the Boolean fans \cite[Proposition~2.6]{BHMPW22}. The star Lefschetz condition holds on the ordinary Bergman fan by \cite[Proposition~2.2(3)]{LP}, and follows on the augmented Bergman fan by induction on rank using the product decomposition of stars at proper flats $F$ into the ordinary fan of $N/F$ and the augmented fan of $N|F$ (and the augmented contraction fan along coordinate rays) \cite[Proposition~2.7]{BHMPW22}, \cite[Proposition~2.2(2)]{LP}.

On the ordinary Bergman fan inside the permutohedral fan, taking $f=-\mathrm{id}$ yields the restrictions of $\mathcal S_{M_i}^\vee$ \cite[Examples~1.1 and~3.10]{LP}. On the augmented Bergman fan inside the stellahedral fan \cite[Definition~5.10]{EHL}, replace $E$ by $E\sqcup E'$ and form the parallel doubling $D(M_i)$ (whose rank on $J\subseteq E\sqcup E'$ is $\rk_{M_i}$ of the projection of $J$ to $E$, preserving quotient relations). Setting $f(v)=(v,0)$ assigns coordinate functions $t_e$ to $E$ and $0$ to $E'$. Dualizing negates the nonzero roots to $-t_e$, which matches the Chern roots of $\mathcal Q^{\rm aug}_{M_i^*}$ in \cite[Proposition~4.4]{EHL}. Indeed, on a maximal stellahedral cone indexed by $I$ and a flag, the minimum-weight basis of $D(M_i)$ selects the minimum-weight basis of $M_i|_{E\setminus I}$ from $E$ (whose complement in $E\setminus I$ is the oppositely ordered greedy basis of $M_i^*/I$) together with $r_i-\rk_{M_i}(E\setminus I)$ zero-weight elements from $E'$, as in \cite[Example~8.6 and \S8.3]{EHL} and \cite[Example~1.2]{LP}. Because projecting $P(D(M_i))$ onto $\R^E$ produces the independence polytope of $M_i$, whose normal fan is refined by the stellahedral fan \cite[\S3.3]{EHL}, $f$ is a single compatible weight map across all stages, and Theorem~\ref{thm:fan-flags}(a) applies.
\end{proof}

Notice also that for loopless $N$, the ordinary Bergman fan is the star of the ray $\varnothing$ in the augmented Bergman fan \cite[Remark~2.5]{BHMPW22}, along which $\mathcal Q^{\rm aug}_{M_i^*}$ restricts to the ordinary tautological quotient $\mathcal Q_{M_i^*}$ \cite[\S4.1]{EHL} (matching $f=\mathrm{id}$ in \cite[Example~1.1]{LP}). In the augmented setting, using $\mathcal Q^{\rm aug}_{M^*}$ rather than $(\mathcal S^{\rm aug}_M)^\vee$ is essential even for $s=1$: for $U_{1,1}$ on a singleton, the stellahedral variety is $\PP^1$ and $\mathcal S^{\rm aug}_{U_{1,1}}=\OO(1)$ \cite[Remark~4.3]{EHL}, whose dual $\OO(-1)$ fails formal nefness.

\subsection{Intersection cohomology of matroids and pullbacks from toric varieties}
Braden, Huh, Matherne, Proudfoot, and Wang realize the intersection cohomology module $\IH(N)$ of $N$ as a graded module summand of $A_{\rm aug}(N)$ over the graded M\"obius algebra $H(N)$, polarized by strictly positive combinations $h=\sum_{e\in E}a_e y_e$ ($a_e>0$) of the nef coordinate divisors $y_e$ \cite[Theorem~1.6]{BHMPW26}, \cite[Definition~3.7]{EHL}. Whenever multiplication by every Chern operator $c_j(F_i)$ of the flag on $A_{\rm aug}(N)$ in Corollary~\ref{cor:both-chow-rings} preserves the subspace $\IH(N)$---for example whenever every $c_j(F_i)$ lies in $H(N)$---self-adjointness and nondegeneracy of the degree pairing on $\IH(N)$ \cite{BHMPW26} preserve $\IH(N)^\perp$ as well, extending the orthogonal splitting to each projective bundle module so that the monodromy and small-perturbation argument of Theorem~\ref{thm:fan-flags} makes the restricted flag formally filtered nef on $\IH(N)$.

Such invariance does not hold automatically, however: for the Boolean matroid $N=U_{2,2}$, the stellahedral surface is the blowup $q\colon X\to\PP^1\times\PP^1$ at one torus fixed point with exceptional divisor $D$, and $\IH(N)=H(N)$ has degree-one piece $\IH^1(N)$ spanned by the two ruling pullbacks $h_1,h_2$ \cite{BHMPW26}. Taking $M=U_{1,2}=M^*$, the class $c_1(\mathcal Q^{\rm aug}_M)$ corresponds to the triangle $\operatorname{conv}(0,e_1,e_2)$ \cite[Proposition~4.6]{EHL}, whose linear system $\operatorname{span}(1,z_1,z_2)$ of bidegree $(1,1)$ has a simple base point at the blown-up point. Hence $c_1(\mathcal Q^{\rm aug}_M)=h_1+h_2-[D]$ pairs to $1$ with $D$ (whereas $h_1\cdot D=h_2\cdot D=0$) and carries $1\in\IH^0(N)$ outside $\IH^1(N)$.

We conclude this section by descending tautological flags on fans along morphisms from smooth projective varieties to projective toric varieties.

\begin{corollary}\label{cor:fan-geom}
Let $\Sigma$ be a rational projective fan equipped with a compatible weight map $f$ for $M_1\preceq\cdots\preceq M_s$ as in Theorem~\ref{thm:fan-flags}, let $g\colon Y\to X_\Sigma$ be a morphism from a smooth complex projective variety, and let $F_i$ be the formal bundle on $Y$ obtained by pulling back the ordinary operational Chern classes of $C_i(u)$ (with $F_0=0$). For every flag of quotient vector bundles $W_1\twoheadrightarrow\cdots\twoheadrightarrow W_s$ on $Y$ with $W_1$ globally generated, the formal bundle $\sum_{i=1}^s(F_i-F_{i-1})\otimes W_i$ is formally nef whenever its rank is positive.
\end{corollary}

\begin{proof}
Working on each connected component, assume that $Y$ is irreducible. Let $\Delta$ be the normal fan of $P=\sum_{i=1}^s P(M_i)$ with lineality space $L$, and put $X_P=X_{\Delta/L}$. Choosing a lattice complement to $L\cap\mathbb Z^E$ splits $\mathrm{PP}(\Delta)=\mathrm{PP}(\Delta/L)\otimes\operatorname{Sym}(L^*)$ over $\R$. Restricting the piecewise Chern polynomials to the complement discards multiples of global linear forms on $L$ and defines equivariant operational classes on $X_P$ \cite[Theorem~1]{Payne} whose ordinary classes pull back along $X_\Sigma\to X_P$ to the coefficients of $C_i(u)$, and along the proper surjective toric morphism $X_E\to X_P$ (induced by negation followed by quotienting by $L$) to the Chern classes of $\mathcal S_{M_i}^\vee$ \cite[Example~1.1]{LP}.

Choosing a closed point of $Y\times_{X_P}X_E$ over the function field of $Y$ and resolving its closure yields a generically finite surjection $p\colon Y'\to Y$ from a smooth projective variety $Y'$ equipped with a morphism $Y'\to X_E$ compatible with $g$. By Theorem~\ref{thm:mat}, the pullback of the incidence bundle to $Y'$ is formally nef. For any $h\in\operatorname{Amp}(Y)$, the identity $p_*p^*=(\deg p)\,\mathrm{id}$ splits $H^*(Y',\R)=p^*H^*(Y,\R)\oplus\ker p_*$ compatibly with all pulled-back Chern operators and with $\xi$, and since $p^*h$ is nef on $Y'$ while both $p^*h$ and $\xi$ lie in the closed polarizing cone of the pullback projective bundle, the monodromy-constancy and small-perturbation argument of Theorem~\ref{thm:fan-flags} shows that $h+t\xi$ polarizes the projective bundle module on $Y$ for all $t>0$.
\end{proof}

\section{Applications}\label{sec:applications}

We now apply our formal filtered Hodge--Lefschetz package in two settings. After recovering the classical geometric positivity criteria of Fulton--Lazarsfeld and Fulton (together with cycle-level positivity on projective toric varieties) in Subsection~\ref{sec:classical-pos}, we devote Subsections~\ref{subsec:comparison}--\ref{subsec:ordinary-truncations} to four applications in matroid theory. Taken together, these combinatorial applications explain what filtered Schubert positivity achieves beyond single-matroid Hodge theory: it produces intersection and finite-difference inequalities that detect whether one matroid is a quotient of another (Subsection~\ref{subsec:comparison}), endows the combinatorial $K$-class of a possibly nonrealizable flag matroid on a partial flag variety with the Schubert-positivity properties of a genuine torus orbit closure (Subsection~\ref{sec:orbit-chow}), extracts matroid morphisms and their canonical Higgs factorizations intrinsically from the supports of derived Schubert numbers on augmented Chow rings (Subsection~\ref{sec:schubert-quotients}), and establishes Lorentzianity for truncation invariants and reduced rank polynomials of arbitrary loopless matroids (Subsection~\ref{subsec:ordinary-truncations}).

\subsection{Classical numerical positivity}\label{sec:classical-pos}
On a smooth complex projective variety $X$, if $F_s\twoheadrightarrow\cdots\twoheadrightarrow F_1\twoheadrightarrow F_0=0$ is a quotient flag with $F_s$ nef, then every $F_i$ is nef \cite[Theorem~6.2.12]{Laz}. Since $\OO_{\PP(F_i)}(1)$ is nef and relatively ample, each positive-rank $F_i$ is formally nef, so $F_\bullet$ (as well as any real twist with nef top bundle) is formally filtered nef by Theorem~\ref{thm:abstract-stages}.

Fix marked ranks $0=r_0<r_1<\cdots<r_s=r$, let $P$ be a homogeneous polynomial of degree $d$ symmetric in each block of $r_i-r_{i-1}$ Chern roots of $\ker(F_i\to F_{i-1})$, and expand $P=\sum_{\ell(w)=d}a_w\mathfrak S_w$ over compatible permutations $w$.

\begin{corollary}[Fulton--Lazarsfeld, Fulton]\label{cor:classical-pos}
The intersection number $\int_X P(F_\bullet)$ is nonnegative for every quotient flag of marked ranks $r_1<\cdots<r_s$ with $F_s$ nef on a $d$-dimensional complex projective variety $X$ if and only if $a_w\ge0$ for all $w$. Likewise, $\int_X P(F_\bullet)>0$ whenever $F_s$ is ample on such $X$ if and only if $a_w\ge0$ for all $w$ and $P\ne0$.
\end{corollary}

For ample bundles, Corollary~\ref{cor:classical-pos} was proved by Fulton and Lazarsfeld for $s=1$ \cite[Theorem~I]{FL83} and by Fulton for general $s$ \cite[\S1, Theorem]{Ful95}. On compact K\"ahler manifolds, Demailly, Peternell, and Schneider proved nonnegativity for Schur classes of nef bundles against K\"ahler powers on analytic cycles \cite[Theorem~2.5]{DPS}, from which Fulton deduced the nef filtered implication via cone intersections on degeneracy cones \cite[Proposition~2.3]{DPS}, \cite[Theorem$'$]{Ful95}.

\begin{proof}
When $X$ is smooth and $F_s$ is nef, we have $\int_X\mathfrak S_w(F_\bullet)\ge0$ by Corollary~\ref{cor:schubert-index-module}. More generally, twisting by the hyperplane class on $X\times\PP^j$ and pairing against $H^j$ for any nef divisor $H$ on $X$ yields $\int_X H^j\mathfrak S_w^{[j]}(F_\bullet)\ge0$ upon integration over $\PP^j$. If $F_s$ is ample, pick $H\in\operatorname{Amp}(X)$ and $\epsilon>0$ small enough that $G_\bullet=F_\bullet\langle-\epsilon H\rangle$ has nef top bundle. Expanding $F_\bullet=G_\bullet\langle\epsilon H\rangle$ and retaining the $j=d$ term yields
\[
\int_X\mathfrak S_w(F_\bullet)\ge\epsilon^d\mathfrak S_w(1,\ldots,1)\int_X H^d>0,
\]
where $\mathfrak S_w(1,\ldots,1)>0$ by the compatible-word formula \cite[Theorem~1.1]{BJS} (for singular $X$, choose the ample twist prior to resolving singularities and apply the projection formula).

Conversely, on $T=\operatorname{Fl}(r_1,\ldots,r_s;N)$ with $N\ge r+d$, pulling the globally generated universal flag $Q_\bullet$ back to a resolution $Y_w$ of the $d$-dimensional Schubert variety dual to $\mathfrak S_w$ yields $\int_{Y_w}P(Q_\bullet)=a_w$ \cite[\S3]{Ful95}, forcing $a_w\ge0$ under the nef hypothesis (even if tested over compact K\"ahler manifolds). Under the ample hypothesis, choosing $H$ very ample on $Y_w$, passing to a Bloch--Gieseker finite cover on which $Q_\bullet\langle H/b\rangle$ becomes an ample vector-bundle flag \cite[Lemma~2.1]{BG}, dividing by the degree of the cover, and letting $b\to\infty$ likewise forces $a_w\ge0$ \cite[proof of Proposition~3.4]{FL83} (compare \cite[Remark~8.3.17]{Laz}), while strict positivity rules out $P=0$.
\end{proof}

Similarly, if $\Gamma$ is a product of compatible Schubert classes of tautological flags of flag matroids on a projective toric variety $X_\Sigma$ as in Corollary~\ref{cor:fan-geom}, pulling back to a projective resolution $g\colon\widetilde Z\to Z\hookrightarrow X_\Sigma$ of any irreducible subvariety $Z$ makes the pulled-back flags formally filtered nef on $\widetilde Z$ by Corollary~\ref{cor:fan-geom} and Theorem~\ref{thm:abstract-stages}. Hence $\deg(\Gamma\cap[Z])=\int_{\widetilde Z}g^*\Gamma\ge0$ whenever $\dim Z=\deg\Gamma$ (so $\Gamma$ pairs nonnegatively with every effective cycle of complementary dimension), while Theorem~\ref{thm:ineq-schubert-volume} yields Lorentzian divisor inequalities on $\widetilde Z$ whenever $\dim Z>\deg\Gamma$.

\subsection{Inequalities for matroid quotients}\label{subsec:comparison}\label{sec:flag-applications}\label{subsec:higher-differences}
Recent Hodge-theoretic results for tautological classes of matroids establish the K\"ahler package for projective bundle rings associated with a single matroid or an arbitrary tuple of matroids on $E$ \cite[Theorems~1.4 and~3.9 and Example~3.10]{LP}, as well as degree and Hodge-index inequalities for products of Schur classes of individual matroids \cite[Theorem~1.1 and \S5]{LS}. Because those statements require no compatibility among the matroids in a tuple, they apply only to polynomials built from the individual tautological bundles and cannot detect whether a chain $M_1\preceq\cdots\preceq M_s$ forms a flag matroid. By contrast, multi-descent Schubert polynomials depend essentially on the successive differences $F_i-F_{i-1}$: geometrically, they reflect the existence of compatible quotient maps between stages, yet Theorem~\ref{thm:mat} and Corollary~\ref{cor:both-chow-rings} show that their positivity persists for arbitrary nonrealizable flag matroids while failing as soon as the quotient relation is broken.

To see this concretely for two marked ranks $(1,3)$ on $n\ge4$ elements, let $N\preceq M$ be a rank-one quotient of a rank-three matroid on $E$, and set $d=n-1$, $a=c_1(F_N)$, and $b=c_1(F_M)$. Applying the divided-difference operator $\partial_2$ to the dominant monomial $\mathfrak S_{3412}=x_1^2x_2^2$ yields $\mathfrak S_{3142}=x_1^2(x_2+x_3)$, so Theorem~\ref{thm:mat} and Corollary~\ref{cor:schubert-index-module} imply that for all nef divisors $D_1,\ldots,D_{d-3}$ on $X_E$,
\[
\Gamma=\mathfrak S_{3142}(N,M)=a^2(b-a),\qquad \int_{X_E}\Gamma D_1\cdots D_{d-3}\ge0
\]
(and similarly on $A_{\rm ord}(N_0)$ or $A_{\rm aug}(N_0)$ of top degree $d\ge3$ in Corollary~\ref{cor:both-chow-rings}, with quadratic Hodge-index inequalities whenever $d\ge5$ by Theorem~\ref{thm:ineq-schubert-volume}).

The quotient assumption $N\preceq M$ is indispensable already for $n=4$ realizable matroids: taking $N=U_{1,4}$ and $M$ on $\{1,2,3,4\}$ with sole basis $\{1,2,3\}$, the element $4$ is a loop in $M$ but not in $N$, so $N\not\preceq M$. Because the coordinate $3$-plane is torus-fixed we have $b=0$, whereas $a$ is the hyperplane pullback along $X_{[4]}\to\PP^3$, yielding $\int_{X_{[4]}}a^2(b-a)=-1$. Algebraically, in ambient dimension at least $6$ where no degree-$3$ relations occur, Littlewood--Richardson expansions of products $x_1^i s_\lambda(x_1,x_2,x_3)$ ($i+|\lambda|=3$) of Schur classes of the two quotients always contain either $x_1^3$ (when $\lambda$ has at most one row) or $x_1x_2x_3$ (otherwise), neither of which appears in $\mathfrak S_{3142}=x_1^2x_2+x_1^2x_3$.

\subsubsection{Iterated differences of mixed degrees}
A classical consequence of the Khovanskii--Teissier inequalities is that for two nef divisor classes $a$ and $b$ on a $d$-dimensional base, the mixed degrees $m_j=\int_V a^{d-j}b^j$ ($0\le j\le d$) form a log-concave sequence. When $a=c_1(F_1)$ and $b=c_1(F_2)$ come from a formally filtered nef flag of ranks $(1,r)$ with $a,b$ nef, the quotient relation imposes higher-order monotonicity and log-concavity on the iterated forward differences $\Delta^q m_j$ ($\Delta m_j=m_{j+1}-m_j$). Indeed, for $q\ge1$ and $0\le j\le d-q$, the binomial expansion reads
\[
\Delta^q m_j=\int_V a^{d-j-q}(b-a)^q b^j,
\]
which is nonnegative whenever $j\le d-2q-r+2$ by Proposition~\ref{prop:block-criterion}. Moreover, when $k=d-rq\ge0$, Proposition~\ref{prop:block-criterion} identifies $G=a^{r-1}(b-a)$ with a single Schubert class, so Theorem~\ref{thm:ineq-schubert-volume} shows that
\[
\sum_{j=0}^k\frac{\Delta^q m_j}{(k-j)!j!}u^{k-j}v^j=\frac1{k!}\int_V G^q(ua+vb)^k
\]
is Lorentzian or zero. Hence $(\Delta^q m_j)_{0\le j\le k}$ is log-concave with interval support \cite[Example~2.26]{BH20}. In particular, for $(r,q)=(3,1)$,
\begin{equation}\label{eq:multidegree-differences}
\Delta m_j\ge0\quad(0\le j\le d-3),\qquad (\Delta m_j)^2\ge\Delta m_{j-1}\Delta m_{j+1}\quad(1\le j\le d-4).
\end{equation}

For example, take the nonrealizable Fano matroid $F_7$ (whose incidences force the three diagonal points of a projective frame to be collinear, forcing characteristic $2$) together with its rank-one truncation $U_{1,7}$, where $d=6$. Along the torus-invariant curve $C\subseteq X_E$ swapping the first two elements of an order, the greedy Fano basis retains both elements while the $U_{1,7}$ basis switches, so $(b-a)\cdot C=-1$ and $b-a$ is not nef. Even so, the rational map $\PP^6\dashrightarrow\PP^{\binom{7}{3}-1}$ defined by cubic basis monomials of $F_7$ has base locus equal to the seven coordinate $2$-planes dual to the Fano lines. General $j$-planes miss this base locus for $j\le3$, so $m_j=3^j$ for $0\le j\le3$. Meanwhile, a general $4$-plane meets the seven $2$-planes transversely in seven simple points (locally cut out by normal coordinates since two line coordinates plus one off-line coordinate form a Fano basis), so $m_4=3^4-7=74$ by B\'ezout's theorem, confirming $\Delta m_0,\ldots,\Delta m_3=2,6,18,47$ in \eqref{eq:multidegree-differences}. Likewise, expanding $(a+u)^2(b-a+2u)=\Gamma+2ab\,u+(b+3a)u^2+2u^3$ and applying Theorem~\ref{thm:ineq-luzheng} with $m=n_1=n_2=2$ against $x^{d-1}$ for any nef divisor $H$ shows that
\[
\int_{X_E}(b+3a)H^{d-1},\qquad 2(d-1)\int_{X_E}abH^{d-2},\qquad \binom{d-1}{2}\int_{X_E}\Gamma H^{d-3}
\]
is a nonnegative log-concave sequence for every rank-$(1,3)$ flag matroid.

Similarly, for marked ranks $(2,3)$ with $u=c_1(F_1)$, $v=c_2(F_1)$, and $b=c_1(F_2)$, Fulton's class $\mathfrak S_{1432}=uv+(u^2-v)(b-u)$ lies outside the cone generated by products of Schur classes of $F_1$ and $F_2$ \cite[\S4]{Ful95}, yet Theorem~\ref{thm:ineq-schubert-volume} applies to $\mathfrak S_{1432}$ even when two individually $\C$-realizable constituents admit no simultaneous $\C$-realization---for instance $L/9\preceq L\setminus9$ for the non-Pappus matroid $L$ labeled as in \cite[\S1.7.5, Example~7]{BGW} (any common $\C$-realization of $L/9\preceq L\setminus9$ would realize $L$ over $\C$ by taking the kernel line as element $9$), yielding degree-four Lorentzian polynomials on $X_{[8]}$ whenever nonzero.

\subsection{The Chow class of a flag matroid}\label{sec:orbit-chow}
Our second matroidal application explains the geometric meaning of flag matroid Schubert numbers inside the Chow ring of an ambient partial flag variety $T=\operatorname{Fl}(r_1,\ldots,r_s;n)$ with $n=|E|$ (omitting a rank-$n$ stage). When a flag matroid $\mathcal M=(M_1,\ldots,M_s)$ comes from a genuine linear flag in $\C^E$, the closure $Z\subseteq T$ of its $(\C^*)^E$-orbit is a projective subvariety of dimension $d=\dim\bigl(\sum_{i=1}^s P(M_i)\bigr)$ \cite[Corollary~1.13.5]{BGW}. Its fundamental cycle $[Z]\in A_d(T)$ expands with nonnegative integer coefficients in the Schubert basis, and by Riemann--Roch $[Z]$ appears as the highest-dimensional term of $\tau_T([\OO_Z])=\operatorname{ch}([\OO_Z])\operatorname{td}(T_T)\cap[T]$. For an arbitrary (possibly nonrealizable) flag matroid $\mathcal M$, there is no torus orbit closure $Z\subseteq T$, yet Cameron, Dinu, Micha\l ek, and Seynnaeve \cite[Definition~84 and Proposition~85]{CDMS} constructed a combinatorial $K$-theory class $y(\mathcal M)\in K^0(T)$ extending the Grassmannian $K$-classes of Speyer \cite[Appendix]{Speyer} and Fink--Speyer \cite[Proposition~3.2]{FS11}, \cite[Proposition~3.2]{FS}. We show that the top-dimensional component of $\tau_T(y(\mathcal M))$ behaves cohomologically like the fundamental cycle of a virtual orbit closure of dimension $d$.

To account for disconnectedness when $d<n-1$, partition $E$ into the finest blocks $E_1,\ldots,E_c$ over which every $M_i$ decomposes as a direct sum. Basis exchanges inside circuits show that the linear span of differences in $P=\sum_{i=1}^s P(M_i)$ is $L_{\R}$, where
\[
L=\left\{u\in\mathbb Z^E:\ \sum_{e\in E_j}u_e=0\text{ for all }1\le j\le c\right\},
\]
so $d=\dim P=n-c$. Following the product reduction suggested for a single matroid in \cite[proof of Corollary~1.10]{LS}, set $Y=\prod_{j=1}^c X_{E_j}$ of dimension $d$, whose fan refines the normal fan of $P$ in $L_{\R}^\vee$ (with singleton blocks contributing points). Ordering the blocks embeds $Y$ as a torus orbit closure in $X_E$ on which the $i$th total quotient splits as the sum of pullbacks of $\mathcal S_{M_i|E_j}^\vee$ \cite[Proposition~5.3]{BEST}. Appending the free matroid $U_{n,E}$ supplies a rank-$(n-r_s)$ bundle representing $[\OO^n]-F_s$, so specialization defines a ring homomorphism $\Phi\colon A^*(T)\to A^*(Y)$. Dualizing against $D\in A^j(Y)$ via
\[
\deg_T(\alpha\cap C_D)=\int_Y\Phi(\alpha)D\qquad\bigl(\alpha\in A^{d-j}(T)\bigr)
\]
produces cycles $C_D\in A_{d-j}(T)$ whose Schubert coefficients are nonnegative (and integral when $D$ is an integral product of nef divisors) by Theorem~\ref{thm:mat} and Corollary~\ref{cor:schubert-index-module}.

\begin{proposition}\label{prop:orbit-rr}
For every flag matroid $\mathcal M$ on $E$, the highest-dimensional nonzero homogeneous component of $\tau_T(y(\mathcal M))$ is $C_1\in A_d(T)$ in dimension $d=\dim P$, and its coefficients in the Schubert basis are nonnegative integers.
\end{proposition}

\begin{proof}
Let the coordinate torus $(\C^*)^E$ act on $(\C^E)^\vee$ by $t_e^{-1}$, so that at a coordinate flag $B_\bullet$ the $i$th universal quotient has character $\sum_{e\in B_i}t_e$ \cite[\S2.2]{DES}. Choosing a linear order $\sigma_j$ on each block $E_j$ determines a torus fixed point $\sigma$ of $Y$ and a flag $B_\sigma$ formed by the unions of the greedy bases on the blocks. Writing $v_B=\sum_{i=1}^s\mathbf 1_{B_i}$ and $C_B=\operatorname{cone}(P-v_B)$, the equivariant Euler class of $T^*_\sigma Y$ is \cite[\S2.3]{BEST}
\[
q_\sigma=\lambda_{-1}(T^*_\sigma Y)=\prod_{j=1}^c\prod_{a=1}^{|E_j|-1}\bigl(1-t_{\sigma_j(a+1)}/t_{\sigma_j(a)}\bigr).
\]
We claim that the rational lattice-point generating function of $C_B$ satisfies
\begin{equation}\label{eq:orbit-cone}
\operatorname{Hilb}(C_B)=\sum_{\sigma:B_\sigma=B}q_\sigma^{-1}.
\end{equation}
Indeed, take the Minkowski sum of $C_B$ with the product of standard permutohedra on the blocks $E_j$: in each decreasing-weight chamber, the greedy vertex maximizes the pairing, so the vertices of the Minkowski sum correspond precisely to those $\sigma$ with $B_\sigma=B$ (where the weight is nonpositive on $C_B$). At the vertex indexed by $\sigma$, the inward primitive edge directions $e_{\sigma_j(a+1)}-e_{\sigma_j(a)}$ form a $\mathbb Z$-basis of $L$ whose tangent cone at the origin contributes $q_\sigma^{-1}$. Since the recession cone of the Minkowski sum is $C_B$ and $L=L_{\R}\cap\mathbb Z^E$, Ishida's recession-cone identity \cite[Theorem~2.3]{Ish} yields \eqref{eq:orbit-cone} (compare \cite[\S2]{Brion}, \cite[Theorem~4.2]{DES}, and \cite[proof of Proposition~1]{Eur} for two stages).

Now let $A$ be any polynomial in exterior powers of the universal quotients on $T$, with specialization $A'$ at the bundles $F_i$ on $Y$ under the same $(\C^*)^E$-linearization. Since $A'|_\sigma=A|_{B_\sigma}$, by equivariant localization \cite[\S4.7]{Nielsen} and \eqref{eq:orbit-cone} we have
\[
\chi_Y^{\rm eq}(A')=\sum_B A|_B\operatorname{Hilb}(C_B)=\chi_T^{\rm eq}\bigl(A\,y(\mathcal M)^{\rm eq}\bigr),
\]
where in the second equality the cotangent factor of \cite[Definition~2.9]{DES} cancels the localization denominator of \cite[Theorem~2.8(3)]{DES}. Passing to the non-equivariant limit and applying Hirzebruch--Riemann--Roch \cite[Corollary~15.2.1]{Ful98} yields
\[
\tau_T\bigl(y(\mathcal M)\bigr)=\Phi_*\bigl(\operatorname{td}(T_Y)\cap[Y]\bigr),
\]
where $\Phi_*$ is integration-adjoint to $\Phi$: indeed, $\operatorname{ch}(A')=\Phi(\operatorname{ch}(A))$, and because $\operatorname{ch}(\gamma^j(Q-\rk Q))$ starts in degree $j$ with $c_j(Q)$ for each universal quotient $Q$, products of such $\gamma$-operations realize all Chern monomials as leading terms and hence span $A^*(T)_{\mathbb Q}$ by descending induction on degree.

Thus the dimension-$j$ component of $\tau_T(y(\mathcal M))$ equals $C_{\operatorname{td}_{d-j}(T_Y)}$ (vanishing for $j>d$), with $j=d$ component $C_1$ because $\operatorname{td}_0(T_Y)=1$. Finally, $C_1\ne0$ because the tensor product of the Pl\"ucker line bundles specializes on $Y$ to the toric line bundle of $P$ (with fixed-point weights $v_{B_\sigma}$), whose top self-intersection equals the positive normalized volume of $P$ in $L$ \cite[Theorem~13.4.3]{CLS}.
\end{proof}

When $\mathcal M$ is realizable with orbit closure $Z\subseteq T$, we have $y(\mathcal M)=[\OO_Z]$ \cite[(4.7)]{CDMS}, \cite[Remark~2.10]{DES} and hence $C_1=[Z]$ \cite[Theorem~18.3(5)]{Ful98} (for Grassmannian orbit polytopes see also \cite[\S2.4]{GGMS} and \cite[\S3]{FS}, and for $s=1$ compare \cite[Corollary~1.10]{LS}). Writing $Z_j\in A_j(T)$ for the dimension-$j$ component of $\tau_T(y(\mathcal M))$ and $H_i$ for the $i$th Pl\"ucker divisor on $T$ (which specializes to the nef divisor of $P(M_i)$ on $Y$), we obtain Lorentzian multidegree polynomials in the top two dimensions as well as Schubert positivity in small codimensions.

\begin{corollary}\label{cor:orbit-multidegrees}
Let $j\in\{d,d-1\}$, let $\sigma_w$ be a compatible Schubert class with $m=j-\ell(w)\ge0$, and set $b_\nu=\deg_T\bigl(\sigma_w\prod_{i=1}^s H_i^{\nu_i}\cap Z_j\bigr)$ for $|\nu|=m$. Then $\sum_{|\nu|=m}b_\nu t^\nu/\nu!$ is Lorentzian or zero. In particular, every $b_\nu\ge0$, the support $\{\nu : b_\nu>0\}$ is $M$-convex when nonempty, and $b_\nu^2\ge b_{\nu+e_p-e_q}b_{\nu-e_p+e_q}$ whenever $p\ne q$ and $\nu_p,\nu_q>0$.
\end{corollary}

\begin{proof}
For $j=d$, apply Theorem~\ref{thm:ineq-schubert-volume} directly on $Y$ (when $w=1$, the coefficients $b_\nu$ are the mixed volumes of $P(M_1),\ldots,P(M_s)$). For $j=d-1$, we have $Z_{d-1}=C_{c_1(T_Y)/2}$. Since on each factor $X_{E_j}$ the anticanonical polytope has convex support function $u\mapsto\max_e u_e-\min_e u_e$ equal to $1$ on every primitive braid ray, $c_1(T_Y)$ is nef, so differentiating the $(m+1)$-variable Lorentzian polynomial once along $c_1(T_Y)/2$ and setting its variable to zero preserves Lorentzianity, while $M$-convexity and the exchange inequality follow from \cite[Theorem~2.25]{BH20}.
\end{proof}

\begin{corollary}\label{cor:orbit-todd}
The Chow components $Z_{d-k}$ have nonnegative rational Schubert coefficients for $0\le k\le\min(3,d)$, as do $Z_0$ and (when $d\ge1$) $Z_1$. Furthermore, if $|E_j|\le7$ for every block $E_j$, then $Z_j$ is Schubert-positive for all $0\le j\le d$.
\end{corollary}

\begin{proof}
Put $c_i=c_i(T_Y)$. Proposition~\ref{prop:orbit-rr} handles $k=0$, while for $k=1,2,3$ the Todd polynomials are $c_1/2$, $(c_1^2+c_2)/12$, and $c_1c_2/24$ (compare \cite[Lemma~3.8 and Theorem~1.5(2)]{CL18} for the corresponding Berline--Vergne coefficients). By the generalized Euler sequence on smooth toric varieties \cite[\S8.1 and \S12.5]{CLS}, $c_2$ is the sum of the fundamental classes of all codimension-two torus orbit closures. Restricting the matroid flag and the nef class $c_1$ to those closures (and pairing $c_1^2$ on $Y$ itself) yields nonnegative pairings against every complementary Schubert class by Theorem~\ref{thm:mat} and Corollary~\ref{cor:schubert-index-module}.

In dimensions $0$ and $1$, Demazure vanishing \cite[Theorem~9.2.3]{CLS} and Riemann--Roch yield the Ehrhart expansion
\[
\#\bigl(tP(M_i)\cap\mathbb Z^E\bigr)=\sum_{j=0}^d\frac{t^j}{j!}\deg_T(H_i^j\cap Z_j)\qquad(t\in\mathbb Z_{\ge0}),
\]
whose constant term is $1$ and whose linear Ehrhart coefficient $\deg_T(H_i\cap Z_1)$ is nonnegative by Castillo--Liu \cite[Theorem~7.5]{CL21} and Jochemko--Ravichandran \cite[Theorem~1.2]{JR} (vanishing when $P(M_i)$ is a point). Since $H_1,\ldots,H_s$ form the dual basis to the dimension-one Schubert cycles, $Z_1$ is Schubert-positive. When every block satisfies $|E_j|\le7$, each factor $X_{E_j}$ has dimension at most $6$, where every Todd component is an effective combination of torus orbit closures \cite[Theorem~1.5(1)]{CL18}, \cite[\S2]{CL21}, so multiplicativity of $\operatorname{td}(T_Y)$ and orbit-closure restriction yield Schubert positivity in all dimensions.
\end{proof}

In higher dimensions, Schubert positivity of lower Todd components can fail even for $s=1$: Ferroni's matroid with negative quadratic Ehrhart coefficient \cite[Theorem~5.3]{Fer} satisfies $\deg_T(H^2\cap Z_2)<0$, which forces at least one Schubert coefficient of $Z_2$ to be negative because $H^2$ is Schubert-positive.

\subsection{Matroid quotients and Higgs lifts from Schubert numbers}\label{sec:schubert-quotients}\label{subsec:schubert-minors}
Our third matroidal application shows how Schubert intersection numbers in an augmented Chow ring $A_{\rm aug}(N)$ construct new matroids, matroid morphisms, and their Higgs factorizations. Let $N$ be a loopless matroid of rank $d$ on $E$. In $A_{\rm aug}(N)$, each coordinate divisor $y_e$ ($e\in E$) is pulled back from a $\PP^1$ factor and therefore is nef with $y_e^2=0$ \cite[Definition~3.7]{EHL}. Consequently, when we pair a product $\Gamma$ of compatible Schubert classes of the flags on $A_{\rm aug}(N)$ in Corollary~\ref{cor:both-chow-rings} of residual degree $m=d-\deg\Gamma\ge0$ against powers of $\sum_{e\in E}x_e y_e$, all non-squarefree monomials drop out, leaving a multiaffine homogeneous polynomial whose coefficients $w_I=\deg_N\bigl(\Gamma\prod_{e\in I}y_e\bigr)$ ($|I|=m$) record minor-contraction degrees on $N/I$. By Br\"and\'en and Huh's characterization of multiaffine Lorentzian supports \cite[Theorem~2.25]{BH20}, nonvanishing of this polynomial immediately implies that $\{I : w_I>0\}$ is the set of bases of a rank-$m$ matroid on $E$.

\begin{proposition}\label{prop:schubert-minors}
The multiaffine polynomial $\sum_{|I|=m}w_I x^I$ (where $x^I=\prod_{e\in I}x_e$) is Lorentzian or zero. Whenever it is nonzero, its support $\{I : w_I>0\}$ is the set of bases of a matroid on $E$ contained in the independent sets of $N$. Moreover, for each $N$-independent $m$-set $I$, the coefficient $w_I$ depends only on $\operatorname{cl}_N(I)$ and equals $\deg_{N/I}\Gamma_I$, where $\Gamma_I$ is obtained from $\Gamma$ by replacing each constituent $\mathcal Q^{\rm aug}_{M^*}$ with $\mathcal Q^{\rm aug}_{(M\setminus I)^*}\oplus\OO^{\,\rk M-\rk(M\setminus I)}$.
\end{proposition}

\begin{proof}
Because $y_e^2=0$, we have $\sum_{|I|=m}w_I x^I=\frac1{m!}\deg_N\bigl(\Gamma(\sum_{e\in E}x_e y_e)^m\bigr)$, which is Lorentzian or zero by Theorem~\ref{thm:ineq-schubert-volume}. When nonzero, $M$-convexity of its squarefree support is equivalent to the basis-exchange axiom \cite[Theorem~2.25]{BH20}. In $A_{\rm aug}(N)$, the product $\prod_{e\in I}y_e$ vanishes when $I$ is dependent in $N$ and equals $y_{\operatorname{cl}_N(I)}$ when $I$ is independent \cite[Lemma~2.11 and Proposition~2.18]{BHMPW22}. Iterating \cite[Proposition~2.7]{BHMPW22} identifies the star of $I$ with the augmented Bergman fan of the loopless part of $N/I$, compatible with the degree map \cite[Proposition~2.29 and Remark~2.30]{BHMPW22}, and restricting $\mathcal Q^{\rm aug}_{M^*}$ along this star yields $\mathcal Q^{\rm aug}_{M^*/I}\oplus\OO^{\,|I|-\rk_{M^*}(I)}$ \cite[Proposition~4.8]{EHL}. By matroid duality, $M^*/I=(M\setminus I)^*$ and $|I|-\rk_{M^*}(I)=\rk M-\rk(M\setminus I)$, which yields the stated formula while keeping the formal ranks unchanged in $\Gamma$.
\end{proof}

Notice also that for any $(m-2)$-subset $J\subseteq E$ and distinct $i,j,k,l\in E\setminus J$, writing $w_{Jij}=w_{J\cup\{i,j\}}$, the quantities $A=\sqrt{w_{Jij}w_{Jkl}}$, $B=\sqrt{w_{Jik}w_{Jjl}}$, and $C=\sqrt{w_{Jil}w_{Jjk}}$ satisfy all three triangle inequalities: restricting the Hessian of the $J$-derivative to $\{i,j,k,l\}$ produces a $4\times 4$ zero-diagonal nonnegative symmetric matrix with at most one positive eigenvalue (attained on $(1,1,1,1)^T$ whenever nonzero), so its determinant $-(A+B+C)(-A+B+C)(A-B+C)(A+B-C)$ must be nonpositive (in particular, $C=0$ forces $w_{Jij}w_{Jkl}=w_{Jik}w_{Jjl}$).

\subsubsection{Higgs lifts}\label{subsec:derived-higgs}
Suppose now that we vary the codimension simultaneously across all derived layers of $\Gamma$ by shifting every Chern root of every flag by $u$ and expanding $\Gamma(u)=\sum_{j=0}^k G_j u^{k-j}$, where $k=\deg\Gamma$. Here $G_k=\Gamma$, while $G_0>0$ is the positive integer obtained by evaluating each Schubert polynomial at $(1,\ldots,1)$. For $0\le j\le\min(k,d)$ and $|I|=d-j$, set
\[
b_j(I)=\deg_N\Bigl(G_j\prod_{e\in I}y_e\Bigr),\qquad \mathcal B_j=\bigl\{I\subseteq E:\ |I|=d-j,\ b_j(I)>0\bigr\},
\]
so $\mathcal B_0=\mathcal B(N)$ and $b_k(I)=w_I$ when $k\le d$. Recall that if $K$ is a quotient of $N$ of corank $q=\rk N-\rk K$, Higgs \cite{Higgs} (see also Las Vergnas \cite[(2.4)]{LV99}) proved that the matroid morphism $N\to K$ factors canonically into $q$ elementary rank-one quotients $N=K_0\twoheadrightarrow K_1\twoheadrightarrow\cdots\twoheadrightarrow K_q=K$, called the \emph{Higgs lifts} of $K$ relative to $N$. Whereas Eur and Huh \cite[Definition~1.1]{EH20} started from a prescribed morphism $N\to K$ to study the bases across its Higgs lifts, the next theorem shows that the last nonempty support $\mathcal B_q$ of the derived Schubert classes automatically defines a quotient $K_q$ of $N$ whose intermediate Higgs lifts are precisely the intermediate supports $\mathcal B_j$.

\begin{theorem}\label{thm:derived-higgs}
Let $q\le\min(k,d)$ be the largest index for which $\mathcal B_q\ne\varnothing$. For every $0\le j\le q$, the family $\mathcal B_j$ is the set of bases of a rank-$(d-j)$ matroid $K_j$ on $E$. Moreover, $K_q$ is a quotient of $N$, and $K_0,\ldots,K_q$ are the Higgs lifts of $K_q$ relative to $N$, with rank functions
\[
\rk_{K_j}(A)=\min\bigl\{\rk_N(A),\,\rk_{K_q}(A)+q-j\bigr\}\qquad(A\subseteq E).
\]
Equivalently, a $(d-j)$-subset $I\subseteq E$ belongs to $\mathcal B_j$ if and only if $I$ is independent in $N$ and spans $K_q$. In particular, whenever $\mathcal B_k\ne\varnothing$, the support of $\Gamma$ itself determines all derived supports $\mathcal B_0,\ldots,\mathcal B_k$.
\end{theorem}

\begin{proof}
Let $h$ be the hyperplane class on $\PP^k$. Twisting all flags by $h$ over $A_{\rm aug}(N)\otimes H^*(\PP^k,\R)$ preserves formal filtered nefness as in the proof of Theorem~\ref{thm:ineq-luzheng}, so Theorem~\ref{thm:ineq-schubert-volume} implies that
\[
\frac1{d!}\int\Gamma(h)\Bigl(zh+\sum_{e\in E}x_e y_e\Bigr)^d=\sum_{j=0}^{\min(k,d)}\frac{z^j}{j!}\sum_{|I|=d-j}b_j(I)x^I
\]
is Lorentzian (at $z=0$ it equals $G_0\sum_{B\in\mathcal B(N)}x^B\ne0$, and integrating $h^j\cdot G_j h^{k-j}$ over $\PP^k$ picks out the $z^j/j!$ layer). Polarizing $z^j$ to $e_j\bigl((z_s)_{s\in S}\bigr)/\binom{q}{j}$ in $q$ fresh variables indexed by a disjoint $q$-element set $S$ preserves Lorentzianity \cite[Proposition~3.1]{BH20}, so its multiaffine support is the basis family of a rank-$d$ matroid $\widehat N$ on $E\sqcup S$ \cite[Theorem~2.25]{BH20}.

Inspecting the $j=0$ and $j=q$ layers shows that $E$ spans $\widehat N$, $\widehat N\setminus S=N$, and $S$ is independent in $\widehat N$. For any $j$-subset $S_j\subseteq S$, contracting $S_j$ leaves $\rk_{\widehat N/S_j}(E)=\rk_{\widehat N}(E\cup S_j)-j=d-j=\rk(\widehat N/S_j)$, so deleting $S\setminus S_j$ preserves rank and defines a rank-$(d-j)$ matroid $K_j=(\widehat N/S_j)\setminus(S\setminus S_j)$ with basis set $\{I\subseteq E : I\cup S_j\in\mathcal B(\widehat N)\}=\mathcal B_j$. In particular, $K_q=\widehat N/S$ is a quotient of $N=\widehat N\setminus S$. Because $\widehat N$ is invariant under all permutations of $S$, for every $A\subseteq E$ the restriction $(\widehat N/A)|S$ is a uniform matroid of rank $u=q+\rk_{K_q}(A)-\rk_N(A)$, and contraction yields
\[
\rk_{K_j}(A)=\rk_N(A)+\min(j,u)-j=\min\bigl\{\rk_N(A),\,\rk_{K_q}(A)+q-j\bigr\},
\]
matching Higgs's rank formula \cite{Higgs}, \cite[(2.4)]{LV99}. Evaluating on $|I|=d-j$ shows that $\rk_{K_j}(I)=d-j$ holds if and only if $\rk_N(I)=d-j$ and $\rk_{K_q}(I)=d-q$.
\end{proof}

When the flags pull back from $N\setminus e$ along the deletion morphism $\theta_e\colon A_{\rm aug}(N\setminus e)\to A_{\rm aug}(N)$ of \cite[\S3.1]{BHMPW22}, the layer polynomials $B_{j,N}=\sum_{|I|=d-j}b_j(I)x^I$ (putting $B_{j,N}=0$ out of range) satisfy the deletion recurrences
\[
B_{j,N}|_{x_e=0}=B_{j,N\setminus e}\quad(e\text{ not a coloop}),\qquad B_{j,N}=x_e B_{j,N\setminus e}\quad(e\text{ a coloop}).
\]
Indeed, when $e$ is not a coloop, $\theta_e$ preserves the remaining coordinate divisors and top degree \cite[Lemma~3.3]{BHMPW22}. When $e$ is a coloop, top-degree products omitting $y_e$ vanish by the rank drop while intersecting with $y_e$ restricts to the star of $e$, which projects isomorphically onto the augmented fan of $N/e=N\setminus e$ as in Proposition~\ref{prop:schubert-minors}. Since $E\setminus e$ spans $N$ (and hence every quotient $K_j$) when $e$ is not a coloop, every nonempty layer retains a basis avoiding $e$, so deletion preserves $q$ and sends each $K_j$ to $K_j\setminus e$.

\subsection{Principal truncations and the reduced rank polynomial}\label{subsec:ordinary-truncations}
Our final application connects Schubert intersection numbers on ordinary Chow rings to principal truncations of matroids and deduces Lorentzianity for a purely combinatorial rank generating polynomial. Recall first that in the augmented setting, a rank-one constituent supported on a nonempty subset $S\subseteq E$ implements the principal truncation of $N$ at $\operatorname{cl}_N(S)$. Indeed, for any $N$-independent $(d-1)$-set $I$, the minimum support function of the $S$-simplex on the rank-one star $N/I$ vanishes on positive coordinate rays and equals $-1$ on the remaining ray precisely when $S\not\subseteq\operatorname{cl}_N(I)$. The coordinate degree is therefore $1$ if $S\not\subseteq\operatorname{cl}_N(I)$ and $0$ otherwise \cite[Theorem~1.6, Definition~5.10, and Corollary~5.11]{EHL}, \cite[Proposition~3.1.5 and Lemmas~3.2.4 and~3.3.6]{BES}.

Analogously, inside the ordinary Chow ring $A_{\rm ord}(N)$, Backman, Eur, and Simpson \cite[Definition~3.2.1 and Theorem~3.2.3]{BES} showed that capping the Bergman Minkowski weight with the nef simplicial generator $h_S$ (which depends only on $\operatorname{cl}_N(S)$) performs the principal truncation of $N$ at $\operatorname{cl}_N(S)$. Iterating this cap product as an intersection of Minkowski weights on the permutohedral fan \cite[Theorem~3.1]{FuSt} and applying the projection formula, any product $\Gamma$ of compatible Schubert classes of degree $k\le d-1$ from Corollary~\ref{cor:both-chow-rings} and any $m=d-1-k$ nonempty subsets $S_1,\ldots,S_m\subseteq E$ satisfy
\[
\deg_{A_{\rm ord}(N)}(\Gamma h_{S_1}\cdots h_{S_m})=\deg_{A_{\rm ord}(N_{\boldsymbol S})}(\Gamma),
\]
where $N_{\boldsymbol S}$ is obtained by successively truncating at the closures of $S_1,\ldots,S_m$ (with the degree vanishing whenever any truncating flat has rank at most $1$). Holding $S_1,\ldots,S_{m-2}$ fixed and varying the final pair $(A,B)$, nefness of $h_S$ and Theorem~\ref{thm:ineq-schubert-volume} yield $w_{AB}^2\ge w_{AA}w_{BB}$, extending the $\Gamma=1$ truncation inequalities of \cite[Theorems~5.2.4 and~5.3.1]{BES}.

Specializing the two-stage class $a^{r-1}(b-a)$ to the flag $U_{1,E}\preceq M$ on the augmented Chow ring of the Boolean matroid produces a remarkably simple Lorentzian polynomial directly from the rank function of $M$.

\begin{corollary}\label{cor:reduced-rank-polynomial}
For every loopless matroid $M$ of rank $r\ge2$ on $E$, the polynomials
\[
\sum_{|T|\ge r}\bigl(\rk_M(T)-1\bigr)\frac{z^{|T|-r}x^{E\setminus T}}{(|T|-r)!}\qquad\text{and}\qquad \sum_{|T|=t}\bigl(\rk_M(T)-1\bigr)x^{E\setminus T}\quad(2\le t\le|E|)
\]
are Lorentzian.
\end{corollary}

\begin{proof}
On the augmented Boolean Chow ring over $E$, consider the flag $U_{1,E}\preceq M$ with first Chern classes $a$ and $b$. By Proposition~\ref{prop:block-criterion}, $a^{r-1}(b-a)$ is a single Schubert class, so Theorem~\ref{thm:ineq-schubert-volume} applied with the nef divisors $a$ and $(y_e)_{e\in E}$ produces a Lorentzian normalized volume polynomial in $z$ and $(x_e)_{e\in E}$ in which the coefficient of $z^{t-r}x^{E\setminus T}/(t-r)!$ for $t=|T|\ge r$ equals the degree of $a^{t-1}b-a^t$ on the star of $E\setminus T$. On that star, $a$ and $b$ correspond to the standard simplex and the independence polytope of $M|T$ \cite[Proposition~4.6]{EHL}. Over $\PP^t$, their linear systems are generated by hyperplanes and by the monomials $Z_0^{s-|B|}\prod_{e\in B}Z_e$ for $M|T$-independent sets $B$, where $s=\rk_M(T)$. Because these monomials include both $Z_0^s$ and a basis monomial, they have no common factor, so a general line avoids the base locus and the intersection degree equals $s-1=\rk_M(T)-1$ (which is positive at $T=E$ since $r\ge2$).

For $t\ge r$, differentiating $t-r$ times in $z$ and setting $z=0$ yields the homogeneous degree-$(|E|-t)$ polynomial. For $2\le t<r$, replacing $M$ by its rank-$t$ truncation preserves $\rk_M(T)$ on all $t$-subsets (and remains nonzero since $M$ is loopless of rank at least $2$, so some $t$-set has rank at least $2$).
\end{proof}

For instance, when $M=F_7$ is the Fano matroid ($r=3$) and $t=3$, the class $\mathfrak S_{3142}=a^2(b-a)$ produces the degree-$4$ Lorentzian polynomial assigning weight $1$ to the complements of the $7$ Fano lines and weight $2$ to the complements of the $28$ non-collinear triples.


\begin{thebibliography}{BHM+26}

\bibitem[AHK18]{AHK}
K. Adiprasito, J. Huh, and E. Katz, \emph{Hodge theory for combinatorial geometries}, Ann. of Math. (2) \textbf{188} (2018), no.~2, 381--452, [DOI 10.4007/annals.2018.188.2.1].

\bibitem[AHL25]{AHL}
O. Amini, J. Huh, and M. Larson, \emph{A decomposition theorem for Lefschetz modules} (2025). [arXiv:2511.02026v2].

\bibitem[AP26]{AP}
O. Amini and M. Piquerez, \emph{Hodge theory for tropical fans}, Bull. Soc. Math. France \textbf{154} (2026), 713--809, [DOI 10.24033/bsmf.2931].

\bibitem[ADH23]{ADH}
F. Ardila, G. Denham, and J. Huh, \emph{Lagrangian geometry of matroids}, J. Amer. Math. Soc. \textbf{36} (2023), no.~3, 727--794, [DOI 10.1090/jams/1009].

\bibitem[AK06]{AK}
F. Ardila and C. J. Klivans, \emph{The Bergman complex of a matroid and phylogenetic trees}, J. Combin. Theory Ser. B \textbf{96} (2006), no.~1, 38--49, [DOI 10.1016/j.jctb.2005.06.004].

\bibitem[BES24]{BES}
S. Backman, C. Eur, and C. Simpson, \emph{Simplicial generation of Chow rings of matroids}, J. Eur. Math. Soc. (JEMS) \textbf{26} (2024), no.~11, 4491--4535, [DOI 10.4171/JEMS/1350].

\bibitem[BBFK02]{BBFK}
G. Barthel, J.-P. Brasselet, K.-H. Fieseler, and L. Kaup, \emph{Combinatorial intersection cohomology for fans}, Tohoku Math. J. (2) \textbf{54} (2002), no.~1, 1--41, [DOI 10.2748/tmj/1113247177].

\bibitem[BBFK05]{BBFK05}
\bysame, \emph{Combinatorial duality and intersection product: a direct approach}, Tohoku Math. J. (2) \textbf{57} (2005), no.~2, 273--292, [jstage:tmj/57\_2\_273].

\bibitem[BBD82]{BBD}
A. A. Beilinson, J. Bernstein, and P. Deligne, \emph{Faisceaux pervers}, Analysis and topology on singular spaces, I (Luminy, 1981), Ast\'erisque, vol.~100, Soc. Math. France, Paris, 1982, pp.~5--171.

\bibitem[BEST23]{BEST}
A. Berget, C. Eur, H. Spink, and D. Tseng, \emph{Tautological classes of matroids}, Invent. Math. \textbf{233} (2023), no.~2, 951--1039, [DOI 10.1007/s00222-023-01194-5].

\bibitem[BJS93]{BJS}
S. C. Billey, W. Jockusch, and R. P. Stanley, \emph{Some combinatorial properties of Schubert polynomials}, J. Algebraic Combin. \textbf{2} (1993), no.~4, 345--374, [DOI 10.1023/A:1022419800503].

\bibitem[BG71]{BG}
S. Bloch and D. Gieseker, \emph{The positivity of the Chern classes of an ample vector bundle}, Invent. Math. \textbf{12} (1971), 112--117, [DOI 10.1007/BF01404655].

\bibitem[Bor53]{Borel}
A. Borel, \emph{Sur la cohomologie des espaces fibr\'es principaux et des espaces homog\`enes de groupes de Lie compacts}, Ann. of Math. (2) \textbf{57} (1953), 115--207, [DOI 10.2307/1969728].

\bibitem[BGW00]{BGW00}
A. V. Borovik, I. M. Gelfand, and N. White, \emph{Combinatorial flag varieties}, J. Combin. Theory Ser. A \textbf{91} (2000), no.~1--2, 111--136, [DOI 10.1006/jcta.2000.3091].

\bibitem[BGW03]{BGW}
\bysame, \emph{Coxeter matroids}, Progress in Mathematics, vol.~216, Birkh\"auser, Boston, 2003, [DOI 10.1007/978-1-4612-2066-4].

\bibitem[BHM+22]{BHMPW22}
T. Braden, J. Huh, J. P. Matherne, N. Proudfoot, and B. Wang, \emph{A semi-small decomposition of the Chow ring of a matroid}, Adv. Math. \textbf{409} (2022), Paper No.~108646, [DOI 10.1016/j.aim.2022.108646].

\bibitem[BHM+26]{BHMPW26}
\bysame, \emph{Singular Hodge theory for combinatorial geometries}, J. Amer. Math. Soc. (2026), published electronically, [DOI 10.1090/jams/1083].

\bibitem[BLS26]{BLS}
M.-C. Brandenburg, G. Loho, and B. Smith, \emph{Quotients of $M$-convex sets and $M$-convex functions}, Comb. Theory \textbf{6} (2026), no.~1, Paper No.~14, [DOI 10.5070/C66165701].

\bibitem[BH20]{BH20}
P. Br\"and\'en and J. Huh, \emph{Lorentzian polynomials}, Ann. of Math. (2) \textbf{192} (2020), no.~3, 821--891, [DOI 10.4007/annals.2020.192.3.4].

\bibitem[BEV05]{Res}
A. M. Bravo, S. Encinas, and O. Villamayor U., \emph{A simplified proof of desingularization and applications}, Rev. Mat. Iberoam. \textbf{21} (2005), no.~2, 349--458, [DOI 10.4171/RMI/425].

\bibitem[BL03]{BL03}
P. Bressler and V. A. Lunts, \emph{Intersection cohomology on nonrational polytopes}, Compos. Math. \textbf{135} (2003), no.~3, 245--278, [DOI 10.1023/A:1022232232018].

\bibitem[Bri88]{Brion}
M. Brion, \emph{Points entiers dans les poly\`edres convexes}, Ann. Sci. \'Ecole Norm. Sup. (4) \textbf{21} (1988), no.~4, 653--663.

\bibitem[Buc05]{Buch}
A. S. Buch, \emph{Quantum cohomology of partial flag manifolds}, Trans. Amer. Math. Soc. \textbf{357} (2005), no.~2, 443--458, [DOI 10.1090/S0002-9947-04-03655-4].

\bibitem[CDMS22]{CDMS}
A. Cameron, R. Dinu, M. Micha\l ek, and T. Seynnaeve, \emph{Flag matroids: algebra and geometry}, Interactions with Lattice Polytopes, Springer Proc. Math. Stat., vol.~386, Springer, Cham, 2022, pp.~73--114, [DOI 10.1007/978-3-030-98327-7\_4].

\bibitem[CL18]{CL18}
F. Castillo and F. Liu, \emph{Berline--Vergne valuation and generalized permutohedra}, Discrete Comput. Geom. \textbf{60} (2018), no.~4, 885--908, [DOI 10.1007/s00454-017-9950-3].

\bibitem[CL21]{CL21}
\bysame, \emph{On the Todd class of the permutohedral variety}, Algebr. Comb. \textbf{4} (2021), no.~3, 387--407, [DOI 10.5802/alco.157].

\bibitem[Cat08]{Cattani}
E. Cattani, \emph{Mixed Lefschetz theorems and Hodge--Riemann bilinear relations}, Int. Math. Res. Not. IMRN (2008), no.~10, Art. ID rnn025, [DOI 10.1093/imrn/rnn025].

\bibitem[CK82]{CK}
E. Cattani and A. Kaplan, \emph{Polarized mixed Hodge structures and the local monodromy of a variation of Hodge structure}, Invent. Math. \textbf{67} (1982), no.~1, 101--115.

\bibitem[CKS87]{CKS}
E. Cattani, A. Kaplan, and W. Schmid, \emph{$L^2$ and intersection cohomologies for a polarizable variation of Hodge structure}, Invent. Math. \textbf{87} (1987), no.~2, 217--252.

\bibitem[CLS11]{CLS}
D. A. Cox, J. B. Little, and H. K. Schenck, \emph{Toric varieties}, Graduate Studies in Mathematics, vol.~124, Amer. Math. Soc., Providence, RI, 2011.

\bibitem[dCM02]{dCM02}
M. A. A. de Cataldo and L. Migliorini, \emph{The hard Lefschetz theorem and the topology of semismall maps}, Ann. Sci. \'Ecole Norm. Sup. (4) \textbf{35} (2002), no.~5, 759--772, [DOI 10.1016/S0012-9593(02)01108-4].

\bibitem[dCM05]{dCM}
\bysame, \emph{The Hodge theory of algebraic maps}, Ann. Sci. \'Ecole Norm. Sup. (4) \textbf{38} (2005), no.~5, 693--750, [DOI 10.1016/j.ansens.2005.07.001].

\bibitem[dCM09]{dCMsurvey}
\bysame, \emph{The decomposition theorem, perverse sheaves and the topology of algebraic maps}, Bull. Amer. Math. Soc. (N.S.) \textbf{46} (2009), no.~4, 535--633, [DOI 10.1090/S0273-0979-09-01260-9].

\bibitem[DCP95]{DCP}
C. De Concini and C. Procesi, \emph{Wonderful models of subspace arrangements}, Selecta Math. (N.S.) \textbf{1} (1995), no.~3, 459--494, [DOI 10.1007/BF01589496].

\bibitem[DPS94]{DPS}
J.-P. Demailly, T. Peternell, and M. Schneider, \emph{Compact complex manifolds with numerically effective tangent bundles}, J. Algebraic Geom. \textbf{3} (1994), no.~2, 295--345.

\bibitem[Dem74]{Dem}
M. Demazure, \emph{D\'esingularisation des vari\'et\'es de Schubert g\'en\'eralis\'ees}, Ann. Sci. \'Ecole Norm. Sup. (4) \textbf{7} (1974), 53--88, [DOI 10.24033/asens.1261].

\bibitem[DES21]{DES}
R. Dinu, C. Eur, and T. Seynnaeve, \emph{$K$-theoretic Tutte polynomials of morphisms of matroids}, J. Combin. Theory Ser. A \textbf{181} (2021), Paper No.~105414, [DOI 10.1016/j.jcta.2021.105414].

\bibitem[EW14]{EW14}
B. Elias and G. Williamson, \emph{The Hodge theory of Soergel bimodules}, Ann. of Math. (2) \textbf{180} (2014), no.~3, 1089--1136, [DOI 10.4007/annals.2014.180.3.6].

\bibitem[EW21]{EW21}
\bysame, \emph{Relative hard Lefschetz for Soergel bimodules}, J. Eur. Math. Soc. (JEMS) \textbf{23} (2021), no.~8, 2549--2581, [DOI 10.4171/JEMS/1061].

\bibitem[Eur21]{Eur}
C. Eur, \emph{Erratum/amplification for ``$K$-theoretic Tutte polynomials of morphisms of matroids''}, author's note, 2~pp., [Eur\_KTutteFlag\_ErrAmp.pdf].

\bibitem[EH20]{EH20}
C. Eur and J. Huh, \emph{Logarithmic concavity for morphisms of matroids}, Adv. Math. \textbf{367} (2020), Paper No.~107094, [DOI 10.1016/j.aim.2020.107094].

\bibitem[EHL23]{EHL}
C. Eur, J. Huh, and M. Larson, \emph{Stellahedral geometry of matroids}, Forum Math. Pi \textbf{11} (2023), Paper No.~e24, [DOI 10.1017/fmp.2023.24].

\bibitem[EL24]{EL24}
C. Eur and M. Larson, \emph{Intersection theory of polymatroids}, Int. Math. Res. Not. IMRN \textbf{2024} (2024), no.~5, 4207--4241, [DOI 10.1093/imrn/rnad213].

\bibitem[FeS05]{FeSt}
E. M. Feichtner and B. Sturmfels, \emph{Matroid polytopes, nested sets and Bergman fans}, Port. Math. (N.S.) \textbf{62} (2005), no.~4, 437--468, [emis:52803].

\bibitem[FY04]{FY}
E. M. Feichtner and S. Yuzvinsky, \emph{Chow rings of toric varieties defined by atomic lattices}, Invent. Math. \textbf{155} (2004), no.~3, 515--536, [DOI 10.1007/s00222-003-0327-2].

\bibitem[Fer22]{Fer}
L. Ferroni, \emph{Matroids are not Ehrhart positive}, Adv. Math. \textbf{402} (2022), Paper No.~108337, [DOI 10.1016/j.aim.2022.108337].

\bibitem[FS11]{FS11}
A. Fink and D. E. Speyer, \emph{$K$-classes for matroids and equivariant localization}, Discrete Math. Theor. Comput. Sci. Proc. \textbf{AO} (2011), 339--350, [DOI 10.46298/dmtcs.2915].

\bibitem[FS12]{FS}
\bysame, \emph{$K$-classes for matroids and equivariant localization}, Duke Math. J. \textbf{161} (2012), no.~14, 2699--2723, [DOI 10.1215/00127094-1813296].

\bibitem[Ful92]{Ful92}
W. Fulton, \emph{Flags, Schubert polynomials, degeneracy loci, and determinantal formulas}, Duke Math. J. \textbf{65} (1992), no.~3, 381--420, [DOI 10.1215/S0012-7094-92-06516-1].

\bibitem[Ful95]{Ful95}
\bysame, \emph{Positive polynomials for filtered ample vector bundles}, Amer. J. Math. \textbf{117} (1995), no.~3, 627--633, [DOI 10.2307/2375083].

\bibitem[Ful98]{Ful98}
\bysame, \emph{Intersection theory}, 2nd ed., Ergebnisse der Mathematik und ihrer Grenzgebiete (3), vol.~2, Springer-Verlag, Berlin, 1998, [DOI 10.1007/978-1-4612-1700-8].

\bibitem[FL83]{FL83}
W. Fulton and R. Lazarsfeld, \emph{Positive polynomials for ample vector bundles}, Ann. of Math. (2) \textbf{118} (1983), no.~1, 35--60, [DOI 10.2307/2006953].

\bibitem[FuS97]{FuSt}
W. Fulton and B. Sturmfels, \emph{Intersection theory on toric varieties}, Topology \textbf{36} (1997), no.~2, 335--353, [DOI 10.1016/0040-9383(96)00016-X].

\bibitem[GGMS87]{GGMS}
I. M. Gel'fand, R. M. Goresky, R. D. MacPherson, and V. V. Serganova, \emph{Combinatorial geometries, convex polyhedra, and Schubert cells}, Adv. Math. \textbf{63} (1987), no.~3, 301--316, [DOI 10.1016/0001-8708(87)90059-4].

\bibitem[Gro58]{Gro}
A. Grothendieck, \emph{La th\'eorie des classes de Chern}, Bull. Soc. Math. France \textbf{86} (1958), 137--154, [DOI 10.24033/bsmf.1501].

\bibitem[Han73]{Hansen}
H. C. Hansen, \emph{On cycles in flag manifolds}, Math. Scand. \textbf{33} (1973), 269--274, [DOI 10.7146/math.scand.a-11489].

\bibitem[Hig68]{Higgs}
D. A. Higgs, \emph{Strong maps of geometries}, J. Combinatorial Theory \textbf{5} (1968), 185--191, [DOI 10.1016/S0021-9800(68)80054-7].

\bibitem[Ish90]{Ish}
M.-N. Ishida, \emph{Polyhedral Laurent series and Brion's equalities}, Internat. J. Math. \textbf{1} (1990), no.~3, 251--265, [DOI 10.1142/S0129167X90000150].

\bibitem[JR22]{JR}
K. Jochemko and M. Ravichandran, \emph{Generalized permutahedra: Minkowski linear functionals and Ehrhart positivity}, Mathematika \textbf{68} (2022), no.~1, 217--236, [DOI 10.1112/mtk.12122].

\bibitem[KaP15]{KaPr}
S. Kaji and P. Pragacz, \emph{Diagonals of flag bundles} (2015). [arXiv:1503.03217v5].

\bibitem[Kar04]{Karu}
K. Karu, \emph{Hard Lefschetz theorem for nonrational polytopes}, Invent. Math. \textbf{157} (2004), no.~2, 419--447, [DOI 10.1007/s00222-004-0358-3].

\bibitem[Kar19]{Karu19}
\bysame, \emph{Relative hard Lefschetz theorem for fans}, Adv. Math. \textbf{347} (2019), 859--903, [DOI 10.1016/j.aim.2019.03.006].

\bibitem[KK87]{KK}
M. Kashiwara and T. Kawai, \emph{The Poincar\'e lemma for variations of polarized Hodge structure}, Publ. Res. Inst. Math. Sci. \textbf{23} (1987), no.~2, 345--407.

\bibitem[KP87]{KP87}
W. Kra\'skiewicz and P. Pragacz, \emph{Foncteurs de Schubert}, C. R. Acad. Sci. Paris S\'er. I Math. \textbf{304} (1987), no.~9, 209--211.

\bibitem[KP04]{KP}
\bysame, \emph{Schubert functors and Schubert polynomials}, European J. Combin. \textbf{25} (2004), no.~8, 1327--1344, [DOI 10.1016/j.ejc.2003.09.016].

\bibitem[LP26]{LP}
M. Larson and E. Partida, \emph{Hodge theory for combinatorial projective bundles} (2026). [arXiv:2604.21925v1].

\bibitem[LS26]{LS}
M. Larson and A. Stapledon, \emph{Positivity properties of Schur classes} (2026). [arXiv:2609.13025v2].

\bibitem[LS82]{LS82}
A. Lascoux and M.-P. Sch\"utzenberger, \emph{Polyn\^omes de Schubert}, C. R. Acad. Sci. Paris S\'er. I Math. \textbf{294} (1982), no.~13, 447--450.

\bibitem[LV99]{LV99}
M. Las Vergnas, \emph{The Tutte polynomial of a morphism of matroids. I. Set-pointed matroids and matroid perspectives}, Ann. Inst. Fourier (Grenoble) \textbf{49} (1999), no.~3, 973--1015, [DOI 10.5802/aif.1702].

\bibitem[Laz04]{Laz}
R. Lazarsfeld, \emph{Positivity in algebraic geometry. II}, Ergebnisse der Mathematik und ihrer Grenzgebiete (3), vol.~49, Springer-Verlag, Berlin, 2004.

\bibitem[LiZ26]{LZ}
E. Li and W. Zheng, \emph{Hodge--Riemann relations for nef bundles} (2026). [arXiv:2610.37650v2].

\bibitem[LL97]{LL97}
E. Looijenga and V. A. Lunts, \emph{A Lie algebra attached to a projective variety}, Invent. Math. \textbf{129} (1997), no.~2, 361--412, [DOI 10.1007/s002220050166].

\bibitem[LuZ26]{LuZhengHR}
Q. Lu and W. Zheng, \emph{Hodge--Riemann polynomials} (2026). [arXiv:2506.16992v7].

\bibitem[Mac91]{Macdonald}
I. G. Macdonald, \emph{Notes on Schubert polynomials}, Publ. LaCIM, vol.~6, Universit\'e du Qu\'ebec \`a Montr\'eal, Montreal, 1991, [macdonaldschubert.pdf].

\bibitem[McM93]{McM}
P. McMullen, \emph{On simple polytopes}, Invent. Math. \textbf{113} (1993), no.~2, 419--444, [DOI 10.1007/BF01244313].

\bibitem[Mon59]{Monk}
D. Monk, \emph{The geometry of flag manifolds}, Proc. London Math. Soc. (3) \textbf{9} (1959), 253--286, [DOI 10.1112/plms/s3-9.2.253].

\bibitem[Nie74]{Nielsen}
H. A. Nielsen, \emph{Diagonalizably linearized coherent sheaves}, Bull. Soc. Math. France \textbf{102} (1974), 85--97, [DOI 10.24033/bsmf.1771].

\bibitem[Pay06]{Payne}
S. Payne, \emph{Equivariant Chow cohomology of toric varieties}, Math. Res. Lett. \textbf{13} (2006), no.~1, 29--41, [Equivariant-Chow-cohomology.pdf].

\bibitem[Pra96]{Pra96}
P. Pragacz, \emph{Symmetric polynomials and divided differences in formulas of intersection theory}, Parameter Spaces (Warsaw, 1994), Banach Center Publ., vol.~36, Polish Acad. Sci., Warsaw, 1996, pp.~125--177, [DOI 10.4064/-36-1-125-177].

\bibitem[Pra07]{Pra08}
\bysame, \emph{Miscellany on the zero schemes of sections of vector bundles} (notes by \"O.~\"Ozt\"urk), Algebraic Cycles, Sheaves, Shtukas, and Moduli, Trends Math., Birkh\"auser, Basel, 2007, pp.~105--116, [DOI 10.1007/978-3-7643-8537-8\_6].

\bibitem[Pra18]{PragaczKP}
\bysame, \emph{On a certain family of $U(\mathfrak b)$-modules}, Schubert Varieties, Equivariant Cohomology and Characteristic Classes---IMPANGA~15, EMS Ser. Congr. Rep., Eur. Math. Soc., Z\"urich, 2018, pp.~203--224, [DOI 10.4171/182-1/9].

\bibitem[RT23a]{RT}
J. Ross and M. Toma, \emph{Hodge--Riemann bilinear relations for Schur classes of ample vector bundles}, Ann. Sci. \'Ecole Norm. Sup. (4) \textbf{56} (2023), no.~1, 197--241, [DOI 10.24033/asens.2531].

\bibitem[RT23b]{RTnef}
\bysame, \emph{On Hodge--Riemann cohomology classes}, Birational Geometry, K\"ahler--Einstein Metrics and Degenerations, Springer Proc. Math. Stat., vol.~409, Springer, Cham, 2023, pp.~763--793, [DOI 10.1007/978-3-031-17859-7\_39].

\bibitem[Sai88]{Saito88}
M. Saito, \emph{Modules de Hodge polarisables}, Publ. Res. Inst. Math. Sci. \textbf{24} (1988), no.~6, 849--995, [DOI 10.2977/PRIMS/1195173930].

\bibitem[Spe09]{Speyer}
D. E. Speyer, \emph{A matroid invariant via the $K$-theory of the Grassmannian}, Adv. Math. \textbf{221} (2009), no.~3, 882--913, [DOI 10.1016/j.aim.2009.01.010].

\bibitem[Tim99]{Tim}
V. A. Timorin, \emph{An analogue of the Hodge--Riemann relations for simple convex polytopes}, Russian Math. Surveys \textbf{54} (1999), no.~2, 381--426, [DOI 10.1070/RM1999v054n02ABEH000134].

\bibitem[VV03]{VV}
G. Vezzosi and A. Vistoli, \emph{Higher algebraic $K$-theory for actions of diagonalizable groups}, Invent. Math. \textbf{153} (2003), no.~1, 1--44, [DOI 10.1007/s00222-002-0275-2].

\end{thebibliography}
\end{document}